\documentclass[reqno,11pt]{amsart}
\usepackage{amssymb, amsmath,latexsym,amsfonts,amsbsy, amsthm,mathtools}
\usepackage{rotating}
\usepackage{tikz}
\usepackage{calc}
\usepackage{hyperref}
\usepackage{mathabx}
\usepackage{enumitem}  
\usepackage{pifont}    
\usetikzlibrary{calc}
\usepackage{cases}
\usepackage{tikz-cd}
\usepackage{mathabx}
\usepackage{tikz-3dplot}
\usetikzlibrary{calc}
\usetikzlibrary{decorations.markings,arrows}
\usetikzlibrary{shapes.geometric}
\usepackage{color}
\usepackage{hyperref}
\usepackage{extarrows}
\usepackage{graphicx}
\allowdisplaybreaks[4]
\let\pa=\partial
\let\al=\alpha

\let\s=\sigma
\let\f=\frac

\let\om=\omega
\let\G=\Gamma

\let\ve=\varepsilon
\let\pa=\partial

\newcommand{\beq}{\begin{equation}}
\newcommand{\eeq}{\end{equation}}
\newcommand{\ben}{\begin{eqnarray}}
\newcommand{\een}{\end{eqnarray}}
\newcommand{\beno}{\begin{eqnarray*}}
\newcommand{\eeno}{\end{eqnarray*}}

\renewcommand{\theequation}{\thesection.\arabic{equation}}

\newtheorem{theorem}{Theorem}[section]

\newtheorem{lemma}[theorem]{Lemma}

\newtheorem{Theorem}{Theorem}[section]

\newtheorem{Remark}[Theorem]{Remark}

\newcommand{\ud}{\mathrm{d}}

\newcommand{\BS}{\mathbb{S}}

\newcommand{\BR}{\mathbb{R}}

\begin{document}
\title[Vector-valued Allen-Cahn dynamics]{Asymptotic limit Of a vector-valued Allen-Cahn equation
for phase transition dynamics}

\author{Huan Dong}
\address{School of Mathematical Sciences, Peking University, Beijing 100871, China}
\email{huandong@math.pku.edu.cn}

\author{Wei Wang}
\address{School of Mathematical Sciences, Zhejiang University, Hangzhou 310058, China}
\email{wangw07@zju.edu.cn}

\renewcommand{\theequation}{\thesection.\arabic{equation}}
\setcounter{equation}{0}


\begin{abstract}
We study the sharp interface limit  for a vector-valued Allen-Cahn equation.
The limiting system is shown to be a two-phase flow in which the interface evolves by mean curvature flow, while in the bulk regions the solution follows a harmonic map heat flow into $\mathbb{S}^{n-1}$, with a mixed boundary condition across the interface.
 This result extends the work of Bronsard-Stoth \cite{Bron} from the two-dimensional radially symmetric case to a higher-dimensional general setting.
\end{abstract}

\date{\today}
\maketitle
\numberwithin{equation}{section}

\section{Introduction}
Phase transition problems have attracted considerable interest in both analysis and
applications. The simplest diffuse interface model for phase transitions is the scalar Allen-Cahn
equation \cite{AL}:
\begin{equation}\label{allencahn}
    \pa_t u=\Delta u-\frac{1}{\ve^2}F'(u),
\end{equation}
where $x\in\Omega \subseteq \mathbb{R}^m,\ u(x)\in \mathbb{R}$ is a scalar order parameter and $F:\mathbb{R}\to \mathbb{R}$ is a potential function with double wells (e.g., $F(u)=\frac{(u^2-1)^2}{4}$). When the small parameter $\ve\to 0$, the domain $\Omega$ separates into two regions where $u$ tends to $\pm 1$, respectively. The interface between these two regions
evolves according to mean curvature flow.
This asymptotic behavior, referred to as the sharp interface limit, has been rigorously justified in various settings: for the static case see \cite{ML,  ML77}, and for the dynamical case see \cite{Brondy, C1992, deS, evans, il, Ste}.

In \cite{RSKf,RSKr}, Rubinstein-Sternberg-Keller introduced the vector-valued Allen-Cahn equation:
\begin{align}\label{equation:main}
  \partial_t u=\Delta u-\frac{1}{\ve^2}\partial_uF(u),
\end{align}
where $x\in\Omega\subseteq \BR^m,\ u(x,t)\in\BR^n$ is a phase-indicator function and $F:\BR^n\to \BR$ is a smooth potential function which vanishes on two disjoint connected submanifolds $N^\pm\subset\BR^n$. Formal asymptotic analysis as $\ve\to 0$ suggests that the interface evolves by mean curvature flow, while away from the interface $u$ follows the harmonic map heat flow into $N^\pm$. The sharp interface limit for \eqref{equation:main}
has been studied extensively when $F$ has minima at isolated points; see, e.g. \cite{Bethuel, fischer2024, Fonseca, Laux18, Moser90}.
However, when the zero set of $F$ consists of two disjoint submanifolds, the rigorous justification of the limit $\ve\to 0$ for \eqref{equation:main} is a challenging problem, which is known as the Keller-Rubinstein-Sternberg problem.

In \cite{Bron}, Bronsard-Stoth rigorously justified the sharp interface limit of \eqref{equation:main} for the special potential:
\begin{align}\label{potential-2D}
  F(u)=(|u|-a)^2(b-|u|)^2, \ u \in \mathbb{R}^2,
\end{align}
under the assumption of radially symmetry (i.e., $u(x)=u(|x|)$). This model is introduced as a toy model for the general case.
In \cite{LPW2012}, Lin-Pan-Wang studied a general double-well potential $F(u): \mathbb{R}^n \to \mathbb{R}_+$ which attains its global minima on two disjoint, connected compact submanifolds $N^\pm \subset\mathbb{R}^n$. They analyzed the asymptotic behavior, as $\ve\to 0$, of energy minimizers $u_\ve$ for the functional
$$E(u)=\int_\Omega |\nabla u|^2+\frac{1}{\ve^2}F(u)\ud x.$$
The regularity of the limiting minimizing map was later investigated by Lin-Wang \cite{lw2019}.

An important physical case of the Keller-Rubinstein-Sternberg problem is the isotropic-nematic transition of liquid crystals with isotropic elasticity.  Using the matched asymptotic expansions, Fei-Wang-Zhang-Zhang \cite{FWZZ2018} rigorously justified the sharp interface limit as $\ve\to 0$ without the radial symmetry assumption (the radially symmetric case was studied in \cite{symliquid}). More recently, this problem is studied in \cite{ll2021} via the related entropy method which is recently developed in \cite{FLS2020}. We also refer the reader to \cite{lw23,liu24,LWZ,DRW} for analyses of the isotropic–nematic transition problem with anisotropic elasticity.

Another important example of the Keller-Rubinstein-Sternberg problem is the $O(n)$ model introduced by Osting-Wang \cite{OW}:
\begin{equation}\label{matrixallencahn}
\pa_t\mathbf{A}=\Delta \mathbf{A}-\frac{1}{\ve^2}(\mathbf{A}\mathbf{A}^\mathsf{T}\mathbf{A}-\mathbf{A}),
\end{equation}
where $\mathbf{A}:\Omega \subseteq \mathbb{R}^m\to \mathbb{R}^{n\times n}$ is a matrix-valued function. The corresponding potential is $F(\mathbf{A})=\frac{1}{4}\|\mathbf{A}\mathbf{A}^\mathsf{T}-\mathbf{I}\|^2$. The asymptotic limit as $\ve\to 0$ was rigorously justified by Fei-Lin-Wang-Zhang \cite{FLWZ} in the way of matched asymptotic expansion. Compared with the isotropic-nematic problem, a key difference is that the system \eqref{matrixallencahn} is  {\it partially minimally paired} (see \cite[Definition 3.5]{FLWZ}),  which introduces new difficulties in constructing approximate solutions and in establishing a uniform spectral lower bound for the linearized operator around them. We refer to \cite{FLWZ} for details.

Several methods have been developed to analyze sharp interface limits of phase transition systems. The first is the  matched asymptotic expansion, which is introduced by \cite{deS, ABC1994}, and later applied to more complicated systems \cite{Abelsliu2018,  FWZZ2015,  FWZZ2018, FLWZ}.
Another approach uses geometric measure theory and varifold flows for the scalar Allen-Cahn equation, originally developed for the scalar Allen-Cahn equation based on a monotonicity formula involving the ``discrepancy function"\cite{il, PP}; it is unclear whether this method extends to vectorial systems.
The third approach is the relative entropy (or modified energy) method, introduced by Fischer-Laux-Simon \cite{FLS2020} ( see also \cite{jerrard, fischerj}). This method has been successfully applied to various models \cite{linenergy, ll2021, AHFacns, liucpam}. Particularly, Liu \cite{liucpam} recently used it to justify the sharp interface limit of \eqref{equation:main} under the assumption that, near the vacuum manifolds $N^\pm$, the potential $F$ depends only on the distance to $N^\pm$.

In this paper, we study the vector-valued Allen-Cahn system
\begin{align}\label{equation:main1}
  \partial_t u^\ve=\Delta u^\ve-\frac{1}{\ve^2}(\partial_{u} F)(u^\ve),
\end{align}
where $x\in\Omega \subseteq \BR^m,\ u(x,t)\in\BR^n$ and the potential $F:\mathbb{R}^n\to\mathbb{R}$ is given by
\begin{align}\label{Fu}
F(u)=\frac{(|u|^2-a^2)^2(|u|^2-b^2)^2}{4}\text{ with } 0<a<b.
\end{align}
For simplicity, we take $\Omega = \mathbb{T}^m$.
This choice of $F$ generalizes \eqref{potential-2D} in \cite{Bron} to the arbitrary dimension $n\ge 2$. Clearly, $F$ attains its minimum on two concentric spheres $N^-=\{|u|=a\}$ and $N^+=\{|u|=b\}$. We denote
\begin{align}\label{fex}
  f(u):=\partial_uF(u)=(|u|^2-a^2)(|u|^2-b^2)(2|u|^2-a^2-b^2)u.
\end{align}

We decompose the solution as $u^\ve(x,t)=\rho(x,t)\omega(x,t)$ where $\rho=|u|\ge 0$ and $\omega\in\mathbb{S}^{n-1}$ is a unit vector. As $\ve\to0$, the domain $\Omega$ is divided into two regions $U_\pm=\bigcup_{t\in[0,T]}\Omega_t^\pm$ in which $u\to u^\pm=\rho^\pm(x,t)\omega^\pm(x,t)$ with
\begin{align*}
  \rho^+(x,t)=b,~\text{ for } x\in\Omega_t^+, \ \text{and}\  \rho^-(x,t)=a, ~\text{ for }x\in\Omega_t^-.
\end{align*}
 The sharp interface limit system of \eqref{equation:main1} reads
\begin{subnumcases}{\label{model1}}
\partial_t\omega^\pm=\Delta\omega^\pm+|\nabla\omega^\pm|^2\omega^\pm,  \quad (x,t)\in\Omega_t^\pm,  \label{harmonic1}\\
V=\kappa,  \qquad\qquad\qquad\qquad\quad\ \ (x,t)\in\Gamma_t,\label{meancu11}\\
\omega^+=\omega^-,  \qquad\qquad\qquad\qquad\ (x,t)\in\Gamma_t,  \\
b^2\partial_\nu\omega^+ =a^2\partial_\nu\omega^-,   \qquad\qquad\ \ (x,t)\in \Gamma_t,  \label{jumpcondition11}
\end{subnumcases}
where $\nu$, $V$ and $\kappa$ are the unit normal vector, the normal velocity, and the mean curvature of $\Gamma$ respectively.  System \eqref{model1} describes a coupled dynamics: the first equation is the harmonic map heat flow into $\mathbb{S}^{n-1}$ in the bulk phases; the second governs the motion of $\Gamma$ by mean curvature flow; the third and fourth equations provide, respectively, continuity and a flux balance for $\omega^\pm$ across the interface.

The goal of this paper is to rigorously establish the convergence of \eqref{equation:main1} to the sharp interface limit system \eqref{model1} as $\ve \to 0$.
We employ the matched asymptotic expansion methods which is developed in \cite{ABC1994}.
Let $(\Gamma_t,\omega^\pm)$ be a smooth solution of  \eqref{model1} on $[0,T]$, and let $d(x,\Gamma_t)$ be the signed distance function from $x$ to $\Gamma_t$, with $d(x,\Gamma_t)<0$ in $\Omega_t^-$ and $d(x,\Gamma_t)>0$ in $\Omega_t^+$. Define the tubular neighborhood
\begin{align*}
\Gamma_t(\delta)&=\{x:|d(x,  \Gamma_t)|<\delta\}\text{ and }\G(\delta)=\bigcup_{t\in[0,T]}(\G_t(\delta)\times\{t\}).
\end{align*}

Our first main result guarantees the existence of high-order approximate solutions near a given smooth solution $(\Gamma_t,\omega^\pm)$ of \eqref{model1}.
\begin{Theorem}[Construction of approximate solutions]\label{thK-1}
Let  $(\Gamma_t,\omega^\pm)$ be a smooth solution of \eqref{model1} on $[0,T]$. For any $K\in \mathbb{Z}^+$, there exists an approximate solution $u^K$  satisfying
\begin{align}\label{uKapproximate}
\partial_tu^K-\Delta u^K+\ve^{-2}f(u^K)=\mathcal{R}^K,
\end{align}
where $\mathcal{R}^K=O(\ve^{K-1})$. Moreover,
\begin{align*}
  |u^K-b\omega^+|\to 0\quad x\in\Omega_t^+,\qquad
  |u^K-a\omega^-|\to 0\quad x\in\Omega_t^-,
\end{align*}
as $\ve\to 0$.
\end{Theorem}
 Our second main result provides a uniform spectral lower bound for the linearized operator around the approximate solution $u^K$.
For convenience, we set
\begin{align}
F(u)=\frac{(|u|^2-a^2)^2(|u|^2-b^2)^2}{4}:=\frac{1}{2}G(|u|^2),
\end{align}
where $G(s)=\f12 (s-a^2)^2(s-b^2)^2.$
Hence $\partial_uF(u)=G'(|u|^2)u$, and \eqref{equation:main1} can be rewritten as
\begin{align*}
\partial_t u^\ve =\Delta u^\ve -\frac{1}{\ve^2}G'(|u^\ve|^2)u^\ve.
\end{align*}
\begin{Theorem}[Spectral lower bound estimate]\label{th:uA}
Let $u^K$ be an approximate solution constructed as in Theorem \ref{thK-1}. Then for any $u\in H^1(\Omega;\mathbb{R}^n)$ and $t\in [0,T]$,
\begin{align}\label{th2ineq}
\int_{\Omega}\Big\{|\nabla u|^2+\frac{1}{\ve^2}\Big(G'(|u^K|^2)|u|^2+2G''(|u^K|^2)(u^K\cdot u)^2\Big)\Big\}\ud x\geq -C\int_{\Omega}|u|^2\ud x,
\end{align}
where $C>0$ is independent of $\ve$. Here $$G'(s)=(s-a^2)(s-b^2)(2s-a^2-b^2),\ G''(s)=(2s-a^2-b^2)^2+2(s-a^2)(s-b^2).$$
\end{Theorem}
Combining Theorems \ref{thK-1} and \ref{th:uA}, we can estimate the error between the exact solution and the approximate solution  via an energy method.
\begin{Theorem}[Nonlinear stability]\label{th:error}
Let $k=3([\frac{m}{2}]+1)+3$ and $K=k+1$. Let $u^K$ be an approximate solution as in Theorem \ref{thK-1}, and $u^\ve$ denotes the solution of \eqref{equation:main1}.  Assume that the initial error satisfies
\begin{align}\label{errorestimate0}
\mathcal{E}(u^\ve(\cdot,0)-u^K(\cdot,0))\leq C_0\ve^{2k},
\end{align}
where
$$\mathcal{E}(u)=\sum_{i=0}^{[\frac{m}{2}]+1}\ve^{6i}\int_{\Omega}|\partial^iu|^2\ud x.$$
Then there exist $\ve_0,\ C_1>0$ such that for all $\ve \leq \ve_0$,
\begin{align}\label{errorestimatet}
\mathcal{E}(u^\ve(\cdot,t)-u^K(\cdot,t))\leq C_1\ve^{2k},\ \text{for } \forall \ t\in [0,T].
\end{align}
\end{Theorem}

The above results extend the analysis of Bronsard-Stoth \cite{Bron} for two-dimensional vector-valued problems to the general $n$-dimensional case. Our proof follows the framework of \cite{ABC1994} and \cite{FLWZ} in two steps: constructing the approximate solution and establishing a uniform spectral lower bound estimate \eqref{th2ineq}. This approach allows us to treat the general case without radial symmetry.

For $u^+=b\omega^+\in N^+, u^-=a\omega^-\in N^-$, the infimum of the connecting energy
\begin{align*}
  \inf_{u(z)\in H^1_{loc}(\mathbb{R})\atop
  u(\pm\infty)=u^\pm} \int_\BR \Big\{ \frac{1}{2}|u'|^2+ F(u)\Big\}\ud z
\end{align*}
can be attained only when $\omega^+=\omega^-$. Hence the problem is {\it partially minimally paired} (cf. \cite{FLWZ}), which introduces extra difficulties in the construction of approximate solutions. Following \cite{FLWZ}, we employ the idea of constructing ``quasi-minimal connected orbits". However, due to the asymmetry between $U^+$ and $U^-$ (i.e. $a\neq b$), the resulting orbit $u_0(z,x,t)$ satisfies only a first-order (in $d_0$) degenerate remainder estimates near the interface:
\begin{align}
   |\partial_z^2u_0-f(u_0)|=O(d_0e^{-\al_0|z|}),\qquad \text{ for }|d_0(x,\Gamma_t)|<\delta.
\end{align}
This weak degeneracy introduces additional difficulties in analyzing the system of equations at each order of the expansion, especially the subleading-order system.

On the other hand, to establish the spectral lower bound of the linearized operator
around $u^K$ (constructed in Theorem \ref{thK-1}),
we introduce the local decomposition near the interface:
$$u(x,t)=\sum_{i=0}^{n-1}a_i(x,t)E_i(x,t)$$
where $E_0=\frac{u^K}{|u^K|}$, $E_i\in \BS^{n-1}$, $E_i\cdot E_j=\delta_i^j$ for $(0\leq i,j\leq n-1)$ and $\pa_rE_j=-(\pa_rE_0\cdot E_j)E_0\text{ for }1\leq j\leq n-1$. This yields
\begin{align*}
|\nabla u|^2&+\frac{1}{\ve^2}(G'(|u^K|^2)|u|^2+2G''(|u^K|^2)(u^K\cdot u)^2)\\
\geq&|\nabla_\G u|^2
+(\partial_ra_0)^2+\frac1{\ve^2}\Big(G'(|u^K|^2)+2G''(|u^K|^2)|u^K|^2\Big)a_0^2\\
&+\sum_{i=1}^{n-1}\Big\{(\partial_ra_i)^2+\frac{1}{\ve^2}G'(|u^K|^2)a_i^2\Big\}+\underbrace{2\sum_{i=1}^{n-1}(\partial_ra_0a_i-a_0\partial_ra_i)E_0\cdot\partial_rE_i}_{\text{cross terms}}.
\end{align*}
Here $\nabla_\Gamma$ and $\partial_r$ denote the tangential and normal derivatives, respectively.
By applying the decomposition trick in \cite{FLWZ}, it is expected that $\|a_0\|_{L^2}^2$ and $\|a_i\|_{L^2}^2(1\le i\le n-1)$ can be bounded by
\begin{align*}
\int_{\Gamma(\delta)} \Big\{ (\partial_ra_0)^2+\frac1{\ve^2}\Big(G'(|u^K|^2)+2G''(|u^K|^2)\Big)a_0^2\Big\} \quad\text{and}\quad  \int_{\Gamma(\delta)}\Big\{ (\partial_ra_i)^2+\frac{1}{\ve^2}G'(|u^K|^2)a_i^2\Big\} ,
\end{align*}
respectively. However, it cannot be expected that the normal derivatives $\partial_ra_0, \partial_ra_i$ can be bounded by them.
Moreover, since $\pa_\nu\bar{\omega}|_\G\neq 0$, one has $\pa_r E_0$ is only of order $O(1)$ near the interface (this similar term is of order $O(\ve)$ in \cite{FLWZ}, see \cite[Lemma 7.10]{FLWZ}). This produces new difficulty in estimating the cross terms such as
\begin{align}\label{crossworst}
(\partial_ra_0 a_i-a_0\partial_ra_i)E_0\cdot \pa_rE_i\text{ for }i=1,\cdots,n-1.
\end{align}
A crucial ingredient of our proof is that, if we set $a_0=\rho_0'(\frac{r}{\ve})\bar{a}_0$ and $a_i=\rho_0(\frac{r}{\ve})\bar{a}_i$,
the most singular part of \eqref{crossworst} can be written as
\begin{align}
  \frac{\rho_0'(\frac{r}{\ve})}{\rho_0(\frac{r}{\ve})}\bar{a}_0\bar{a}_i\pa_r\Big(\rho_0^2\Big(\frac{r}{\ve}\Big)\pa_rE_0\Big)\cdot E_i,
\end{align}
near $\Gamma$. By exploiting the jump condition \eqref{jumpcondition11} in the construction of the approximate solution in the inner expansion, we can choose $\rho_0$ and $E_0$ (i.e. $u^K$) such that
$\pa_r\Big(\rho_0^2\Big(\frac{r}{\ve}\Big)\pa_rE_0\Big)\sim O(1)$. This reduction of singularity in $\ve$ is crucial for our estimates; details are given in Lemmas \ref{par} and \ref{cross0ac}.

Finally, we note that, in studying asymptotic limits for phase transition systems, understanding the role played by boundary conditions in rigorous sharp-interface analysis is of fundamental importance. Our present work addresses this question for a specific yet important class of problems, and the approach developed here also sheds light on the analysis of other phase transition models.

The paper is organized as follows. The construction of the approximate solution (Theorem \ref{thK-1}) is carried out via matched asymptotic expansions: Section \ref{section2} treats the  outer expansion away from the interface, and Section \ref{section3} deals with inner expansion near the interface. Section 4 describes the recursive procedure for solving the coupled outer-inner expansion system at each order (see Figure \ref{fig:expansion} for a schematic). In Subsection \ref{section4.4}, the approximate solutions on the whole domain are obtained by gluing the inner and outer expansions. Section \ref{spectralestimate} presents the proof of the uniform-in-$\ve$ spectral lower bound of the linearized operator (Theorem \ref{th:uA}).

\subsection*{Notations}
\begin{itemize}
\item Let $d(x,\Gamma_t)$ be the signed distance from $x$ to $\Gamma_t$, and $\nu=\nabla d|_{\Gamma_t}$ be the unit outer normal of $\Omega_t^-$. We denote
\begin{align*}
 &   \Gamma(\delta)=\{(x,t)\in\Omega\times[0,T]: |d(x,\Gamma_t)|<\delta \},\\
 & U_\pm =\{(x,t)\in\Omega\times[0,T]: d(x,\Gamma_t)\gtrless 0\}.
\end{align*}
We also write $\Gamma=\{(x,t)\in\Omega\times[0,T]:x\in\Gamma_t\}$ for simplicity.
    \item For any two vectors $m=(m_1, m_2, \cdots,m_n)$ and $l = (l_1, l_2,\cdots,l_n)$, we denote $$m\otimes l = (m_il_j )_{1\leq i,j\leq n}\text{ and }m\cdot l=\sum_{i=1}^nm_il_i.$$
     \item Let $\al_0>0$. By $f(z,x,t) = O(e^{-\alpha_0|z|})$ as $z \to \pm\infty$ we mean that there exist constants $C > 0$ and $k \ge 0$ such that
       \begin{align}
       |f(z,x,t)|\leq C|z|^ke^{-\al_0|z|},\;\; \;\text{as}\;z\to\pm\infty.
\end{align}
In this paper, we take
$$0<\alpha_0 <\min \{\sqrt{2}(b^2-a^2)b,\sqrt{2}(b^2-a^2)a\}.$$
\item We use $f|_{\G}$ and $f|_{\G(\delta)}$ to denote
$f(z,x,t)|_{(x,t)\in\G}$ and $f(z,x,t)|_{(x,t)\in\G(\delta)}$ respectively.
\item For convenience, if $f(x,t)=O(d_0)$, we note
\begin{align*}
\frac{f(x,t)}{d_0}:=\begin{cases}
    \frac{f(x,t)}{d_0},\quad \ \ \ \ d_0\neq 0,\\
    \pa_\nu f(x,t),\quad d_0=0.
\end{cases}
\end{align*}
\item We define
$$\langle f(z),g(z)\rangle=\int_\mathbb{R}f(z)g(z)\ud z.$$
\end{itemize}
\section{Outer expansion}\label{section2}
\subsection{Formal outer expansion}
We perform an outer expansion in $U_\pm$ by using the Hilbert expansion method as in \cite{ABC1994,WZZ,WZZ2}. Assume that
\begin{equation}\label{outerexlater}
u^\ve(x,t)=u_0^\pm(x,t)+\ve u_1^\pm (x,t)+\cdots+\ve^k u_k^\pm (x,t)+\cdots.
\end{equation}
Substituting the expansion \eqref{outerexlater} into the equation \eqref{equation:main1}, one has
\begin{align*}
\partial_t\big(&u_0^\pm+\ve u_1^\pm+\cdots+\ve^k u_k^\pm+\cdots\big)=\Delta \big(u_0^\pm+\ve u_1^\pm+\cdots+\ve^k u_k^\pm+\cdots\big)\\
&-\frac{1}{\ve^2}\Big(f(u_0^\pm)+Df(u_0^\pm)\sum_{k=1}^{+\infty}\ve^k u_k^\pm+\sum_{k=1}^{+\infty}\ve^kf^{(k-1)}(u_0^\pm,u_1^\pm,\cdots,u_{k-1}^\pm)\Big).
\end{align*}
Thus we can obtain the system of $u_k^\pm(x,t)$ for $k\geq 0$:
 \begin{align}
O(\ve^{-2})\  \text{term}:&\ f(u_0^\pm)=0;\label{ve-2}\\
O(\ve^{-1})\  \text{term}:&\ Df(u_0^\pm)u_1^\pm=0;\label{u1pmequation}\\
O(\ve^{k-2})\  \text{term}:&\ \partial_t u_{k-2}^\pm=\Delta u_{k-2}^\pm-Df(u_0^\pm)u_{k}^\pm-f^{(k-1)}(u_0^\pm,\cdots,u_{k-1}^\pm)\ \text{for}\ k\geq 2.\label{outuk-2}
\end{align}
The boundary conditions for $u_k^\pm(x,t)$ on $\Gamma$ will be determined to ensure the solvability of  expansions in inner regions (particularly on $\Gamma$).
Moreover, we will use the value of $u_k^\pm(x,t)$ not only in $ U_\pm$, but also in $U_\pm\cup\Gamma(\delta)$. Therefore, we  smoothly extend $u_k^\pm(x,t)$ from $U_\pm$ to $U_\pm\cup\Gamma(\delta)$.
For $f$, we have
\begin{lemma}\label{fex}
Let $v^\ve=\sum_{i=0}^{+\infty}\ve^i v_i$. Then
\begin{align}
f(v^\ve)&=(|v^\ve|^2-a^2)(|v^\ve|^2-b^2)(2|v^\ve|^2-a^2-b^2)v^\ve\notag\\
&=f(v_0)+Df(v_0)\sum_{k=1}^{+\infty}\ve^k v_k+\sum_{k=1}^{+\infty}\ve^kf^{(k-1)}(v_0,v_1,\cdots,v_{k-1}),\label{fexpansion}
\end{align}
where
\begin{align}
Df(v_0)=&(|v_0|^2-a^2)(|v_0|^2-b^2)(2|v_0|^2-a^2-b^2)I\nonumber\\
  &+2\Big((2|v_0|^2-a^2-b^2)^2+2(|v_0|^2-a^2)(|v_0|^2-b^2)\Big)v_0 \otimes v_0,\label{dfu}\\
  f^{(0)}(v_0)=&0,\notag\\
f^{(1)}(v_0,v_1)=&2(v_0\cdot v_1)\big[(2|v_0|^2-a^2-b^2)^2+2(|v_0|^2-a^2)(|v_0|^2-b^2)\big]v_1\nonumber\\
&+(v_1\cdot v_1)\big[(2|v_0|^2-a^2-b^2)^2+2(|v_0|^2-a^2)(|v_0|^2-b^2)\big]v_0\nonumber\\
&+12(v_0\cdot v_1)^2(2|v_0|^2-a^2-b^2)v_0,\label{de2}
\end{align}
and for $k\geq 3$,
\begin{align}
f^{(k-1)}(v_0,\cdots,v_{k-1})=&f^{(k-1)}_1(v_0,v_1,v_{k-1})v_0+f^{(k-1)}_2(v_0,v_{k-1})v_1\nonumber\\
&+f_3^{(k-1)}(v_0,v_1)v_{k-1}+f^{(k-1)}_4(v_0,\cdots,v_{k-2}),\label{fk-1}
\end{align}
with
\begin{align*}
f_1^{(k-1)}(v_0,v_1,v_{k-1})=&2(v_1\cdot v_{k-1})\big[(2|v_0|^2-a^2-b^2)^2+2(|v_0|^2-a^2)(|v_0|^2-b^2)\big]\\
&+24(v_0\cdot v_{k-1})(v_0\cdot v_1)(2|v_0|^2-a^2-b^2),\\
f_2^{(k-1)}(v_0,v_{k-1})=&2(v_0\cdot v_{k-1})\big[(2|v_0|^2-a^2-b^2)^2+2(|v_0|^2-a^2)(|v_0|^2-b^2)\big],\\
f_{3}^{(k-1)}(v_0,v_1)=&2(v_0\cdot v_1)\big[(2|v_0|^2-a^2-b^2)^2+2(|v_0|^2-a^2)(|v_0|^2-b^2)\big],
\end{align*}
 and $f^{(k-1)}_4$ only depends on $v_i$ for $i<k-1$.
 \begin{proof}
These conclusions follow from the formal expansion \eqref{fexpansion} directly.
 \end{proof}
\end{lemma}
\subsection{Solving $u_0^\pm$}
For $k=0$, it yields from \eqref{ve-2} that
$$f(u_0^\pm)=0.$$
In view of \eqref{Fu} and \eqref{fex}, this shows that $u_0^\pm$ are critical points of the bulk energy $F(u).$
Without loss of generality, we take
\begin{equation}
u_0^+=b\omega^+,\ u_0^-=a\omega^- \text{ with } |\omega^\pm|=1.\label{u0pm}
\end{equation}
\subsection{Solving $\rho_1^\pm$}
For $u_1^\pm$, by \eqref{u1pmequation}, we have
\begin{align}\label{u1pmeq}
Df(u_0^\pm)u_1^\pm=0,
\end{align}
which gives that $u_1^\pm$ belongs to the kernel of the operator $Df(u_0^\pm)$.

To proceed, we assume the decomposition:
$$u_1^\pm=\rho_1^\pm\omega^\pm+\sum_{\al}\sigma_{1,\al}^\pm\xi_{\al}^\pm,$$
where $$\{\omega^\pm(x,t),\xi_1^\pm(x,t),\cdots,\xi_{n-1}^\pm(x,t)\}$$ forms a group of complete standard orthogonal basis of $\mathbb{R}^n$,  which will
be determined later. Now, we  derive equations for $\rho_1^\pm$ and $\sigma_{1,\beta}^\pm.$ Taking
$v_0=u_0^+$  in \eqref{dfu}, we deduce from \eqref{u0pm} and \eqref{u1pmeq} that
$$Df(u_0^+)u_1^+=2(b^2-a^2)^2b^2\rho_1^+\omega^+ =0,$$
showing that
\begin{equation}\label{rho1+}
\rho_1^+=0.
\end{equation}
Similarly, we have
\begin{equation}\label{rho1_}
\rho_1^-=0.
\end{equation}
Note that by \eqref{u0pm}, \eqref{rho1+} and \eqref{rho1_}, we obtain that
\begin{equation}\label{u0u1cdot}
    u_0^\pm\cdot u_1^\pm=0.
\end{equation}
\subsection{Solving $\rho_2^\pm$ and $\omega^\pm$}
The equation \eqref{outuk-2} for $k = 2$ gives us that
\begin{align*}
\partial_t u_0^\pm=\Delta u_0^\pm-D f(u_0^\pm)u_2^\pm-f^{(1)}(u_0^\pm,u_1^\pm).
\end{align*}
Applying \eqref{u0pm} and the decomposition $$u_2^\pm=\rho_2^\pm\omega^\pm+\sum_{\alpha}\sigma_{2,\alpha}^\pm\xi_{\alpha}^\pm,$$ we deduce that
\begin{align*}
\rho_0^\pm\partial_t\omega^\pm=\rho_0^\pm\Delta \omega^\pm-2(\rho_0^\pm)^2(b^2-a^2)^2\rho_2^\pm\omega^\pm-f^{(1)}(u_0^\pm,u_1^\pm).
\end{align*}
By taking $v_0=u_0^\pm$ and $v_1=u_1^\pm$ in \eqref{de2}, it is clear that $f^{(1)}(u_0^\pm,u_1^\pm)$ is parallel to $\omega^\pm$.
Thus, for the above equality,  taking the inner product with $\omega^\pm$ and $\xi_\beta^\pm$ separately, by $|\omega^\pm|=1$, we get
\begin{align}
&\rho_0^\pm\Delta \omega^\pm\cdot \omega^\pm-2(\rho_0^\pm)^2(b^2-a^2)^2\rho_2^\pm-f^{(1)}(u_0^\pm,u_1^\pm)\cdot\omega^\pm=0,\label{rho2equation}\\
&(\partial_t\omega^\pm-\Delta \omega^\pm)\cdot\xi_{\beta}^\pm=0\label{0termxi}.
\end{align}
Thus, $\rho_2^\pm$ can be expressed explicitly by \eqref{rho2equation}:
\begin{equation*}
  \rho_2^\pm=\frac{\rho_0^\pm\Delta \omega^\pm\cdot \omega^\pm-f^{(1)}(u_0^\pm,u_1^\pm)\cdot \omega^\pm}{2(\rho_0^\pm)^2(b^2-a^2)^2}.
\end{equation*}
From the fact \eqref{0termxi} and
$$(\pa_t \omega^\pm -\Delta \omega^\pm)\cdot \omega^\pm=|\nabla \omega^\pm|^2,$$
we arrive at $$\partial_t\omega^\pm-\Delta \omega^\pm=|\nabla \omega^\pm|^2\omega^\pm,$$
which is independent of the choice of $\{\xi_{\alpha}^\pm\}$. Indeed, it is just the heat flow of harmonic maps into $\mathbb{S}^{n-1}$.
\subsection{Solving $\rho_{k+1}^\pm$ and $\sigma_{k-1,\beta}^\pm(k\geq 2)$}
We now consider the term of order {\bf $O(\ve^{k-1})$} for $k\geq 2$:
\begin{align}\label{eq:order-ko}
\partial_t u_{k-1}^\pm=\Delta u_{k-1}^\pm-Df(u_0^\pm)u_{{k+1}}^\pm-f^{(k)}(u_0^\pm,\cdots,u_{k}^\pm).
\end{align}
As usual, we take \begin{equation}\label{compose}
u_{k+1}^\pm=\rho_{k+1}^\pm\omega^\pm+\sum_{\alpha}\sigma_{{k+1},\alpha}^\pm\xi_{\alpha}^\pm.
\end{equation}
After substituting the decomposition \eqref{compose} into \eqref{eq:order-ko}, we take separate inner products with $\omega^\pm$ and $\xi_\beta^\pm$ to obtain the equations:
\begin{align}\label{outerrhok}
\partial_t\rho_{k-1}^\pm&+\sum_{\alpha}\sigma^\pm_{k-1, \alpha}\partial_t\xi_\alpha^\pm\cdot \omega^\pm\notag\\
=&\Delta\rho_{k-1}^\pm+\rho_{k-1}^\pm\Delta\omega^\pm\cdot \omega^\pm+\sum_{\alpha}\sigma_{k-1, \alpha}^\pm\Delta\xi_\alpha^\pm\cdot\omega^\pm+2\sum_\al(\nabla\sigma_{k-1, \alpha}^\pm\cdot\nabla)\xi_\alpha^\pm\cdot\omega^\pm\notag\\
&-2(\rho_0^\pm)^2(b^2-a^2)^2\rho_{k+1}^\pm-f^{(k)}(u_0^\pm, \cdots, u_{k}^\pm)\cdot\omega^\pm,
\end{align}
and
\begin{align}\label{outersigmak}
\partial_t\sigma_{k-1, \beta}^\pm&-\Delta\sigma_{k-1, \beta}^\pm-\sigma_{k-1, \beta}^\pm\Delta \xi_\beta^\pm\cdot \xi_\beta^\pm\notag\\
=&-\rho_{k-1}^\pm\partial_t\omega^\pm\cdot\xi_{\beta}^\pm+\rho_{k-1}^\pm\Delta\omega^\pm\cdot\xi_{\beta}^\pm+2(\nabla\rho_{k-1}^\pm\cdot\nabla)\omega^\pm\cdot\xi_{\beta}^\pm\nonumber\\
&+\sum_{\alpha\neq \beta}\sigma_{k-1, \alpha}^\pm\Delta\xi_\alpha^\pm\cdot\xi_{\beta}^\pm+2\sum_{\al\neq\beta}(\nabla\sigma_{k-1, \alpha}^\pm\cdot\nabla)\xi_\alpha^\pm\cdot\xi_\beta^\pm\nonumber\\
&-\sum_{\alpha\neq \beta}\sigma _{k-1, \alpha}^\pm \partial_t\xi_\alpha^\pm \cdot \xi_{\beta}^\pm-f^{(k)}(u_0^\pm, \cdots, u_{k}^\pm)\cdot\xi_\beta^\pm.
\end{align}
The equation \eqref{outerrhok} directly gives us that
\begin{align}\label{outersolve}
\rho_{k+1}^\pm=&\frac{1}{2(\rho_0^\pm)^2(b^2-a^2)^2}\Big(\Delta\rho_{k-1}^\pm+\rho_{k-1}^\pm\Delta\omega^\pm\cdot \omega^\pm+\sum_{\alpha}\sigma_{k-1, \alpha}^\pm\Delta\xi_\alpha^\pm\cdot\omega^\pm\nonumber\\
&+2\sum_\al(\nabla\sigma_{{k-1}, \alpha}^\pm\cdot\nabla)\xi_\alpha^\pm\cdot\omega^\pm-\partial_t\rho_{k-1}^\pm-\sum_{\alpha}\sigma^\pm_{{k-1}, \alpha}\partial_t\xi_\alpha^\pm\cdot \omega^\pm\notag\\
&-f^{(k)}(u_0^\pm, \cdots, u_{k}^\pm)\cdot\omega^\pm\Big).
\end{align}
By \eqref{fk-1} and \eqref{u0u1cdot}, we obtain
\begin{align*}
f^{(k)}(u_0^\pm,\cdots,u_{k}^\pm)\cdot \xi_\beta^\pm=f_2^{(k)}\sigma_{1,\beta}^\pm+f_4^{(k)}\cdot \xi_\beta^\pm.
\end{align*}
We conclude that $f^{(k)}(u_0^\pm,\cdots,u_{k}^\pm)\cdot \xi_\beta^\pm$ depends only on
$u_i^\pm(0\leq i\leq k-1)$ and $\rho_{k}^\pm$.
Thus the right-hand side of \eqref{outersigmak} is independent of $\sigma_{k,\al}^\pm(1\leq \al\leq n-1).$ Hence,
$\sigma_{k-1,\beta}^\pm$ can be determined by \eqref{outersigmak} with boundary conditions which will be fixed in inner expansion.
\section{Inner expansion}\label{section3}

The aim of the inner expansion is to obtain a high-order approximation (up to any power of $\ve$) to the exact solution in the interfacial region $\Gamma(\delta)$. As the exact solution exhibits rapid variations near the interface, we introduce a rescaled variable $z=\frac{d^\ve(x,t)}{\ve}$ to capture this behavior. Our main strategy is to seek functions
\begin{align}\label{innerexpansion}
 u^\ve(z,x,t)=\sum_{k=0}^{+\infty}\ve^k u_k(z,x,t),\quad
  d^\ve(x,t) =\sum_{k=0}^{+\infty} \ve^kd_k(x,t),
\end{align}
for $(x,t)\in\Gamma(\delta)$, such that $\tilde{u}^\ve(x,t)=u^\ve(\frac{d_\ve(x,t)}{\ve},x,t)$
solves the original equation
\begin{align}\label{eq:main-inner}
 \partial_t {u}^\ve=\Delta {u}^\ve-\frac{1}{\ve^2} f({u}^\ve).
\end{align}
Here $d^\ve(x,t)$ is a signed distance function to a perturbative interface $\Gamma_t^\ve$. Then we have $|\nabla d^\ve|^2=1$, which implies
\begin{equation}\label{nablad0nabladk}
\nabla d_0\cdot \nabla d_k=\left\{
\begin{array}{ll}
1,   &k=0, \\
0,  &k=1, \\
-\frac{1}{2}\sum_{i=1}^{k-1}\nabla d_i\cdot\nabla d_{k-i},  &k\geq 2.
\end{array}
\right.
\end{equation}
Moreover, to match the inner and outer expansions, we require that, as $z\to\pm\infty$,
\begin{align}\label{matching-condition}
  D_x^mD_t^nD_z^l[u_k( z,x,t)-u_k^\pm(x,t)]=O(e^{-\alpha_0 |z|})\ \text{for}\ (x,t)\in\Gamma(\delta),
\end{align}
with $k,m,n,l\geq 0$, which are called  matching conditions.

Substituting \eqref{innerexpansion} into \eqref{eq:main-inner}, the leading $O(\ve^{-2})$ system is
\begin{align}\label{-2term}
\partial_{zz}u_0(z,x,t)-f(u_0(z,x,t))=0,
\end{align}
with boundary condition
\begin{align}\label{eq:bd-0}
u_0(+\infty, x, t)=u^+(x,t)(=b\omega^+(x,t)),\ u_0(-\infty,x,t)=u^-(x,t)(=a\omega^-(x,t)).
\end{align}

The equation \eqref{-2term} is the Euler-Lagrange equation to the 1-D connecting energy functional
\begin{align*}
 \int_\BR \Big\{ \frac{1}{2}|u'|^2+ F(u)\Big\}\ud z.
\end{align*}
For a given pair $(u^+, u^-)=(b\omega^+, a\omega^-)\in N^+\times N^-$, the infimum of  this energy with given boundary condition \eqref{eq:bd-0}
can be attained only when $\omega^+=\omega^-$. Therefore, the problem is {\it partially minimally paired} (in the sense of \cite[Definition 3.5]{FLWZ}), which introduces extra difficulties in the construction of approximate solutions. We refer to \cite[Subsection 3.2]{FLWZ} for more detailed explanations.

For  $(x,t)\in \G$,  it holds that $\omega^+=\omega^-$, then the equation \eqref{-2term} admits a solution
\begin{align*}
  u_0(z,x,t)=\rho_0(z)\omega^-(x,t),
\end{align*}
where $\rho_0(z)$ satisfies the ODE:
\begin{equation}\label{1D:profile}
\begin{cases}
  \partial_z^2 \rho_0=f(\rho_0),\\
  \rho_0(-\infty)=a,\ \rho_0(+\infty)=b,\\
  \end{cases}
\end{equation}
where
$$f(s)=(s^2-a^2)(s^2-b^2)(2s^2-a^2-b^2)s,$$
showing
\begin{equation*}
f(a)=f(b)=0, \ f'(a)>0, \ f'(b)>0, \ \int_{a}^{\rho}f(s)ds>0, \forall\rho\in (a, b).
\end{equation*}
\begin{lemma}\label{rerho}
The solutions of \eqref{1D:profile} satisfy
\begin{align}\label{pazrho0}
&\rho_0'=\frac{\sqrt{2}}{2}(\rho_0^2-a^2)(b^2-\rho_0^2)>0,\\
\label{rho0expression}
&b\ln\frac{\rho_0-a}{\rho_0+a}-a\ln\frac{b-\rho_0}{b+\rho_0}=\sqrt{2}ab(b^2-a^2)z+c_0,
\end{align}
where $c_0$ is an arbitrary constant. Moreover,
\begin{align}\label{dexpression}
\int_{\mathbb{R}}\big(\rho_0'(z)\big)^2\ud z=\frac{\sqrt{2}}{15}(b-a)(b^4+a^4+b^3a+a^3b-4a^2b^2)\triangleq e>0,
\end{align}
and there exists  positive constants $C_l$ for $l=1,2,\cdots$ such that
\begin{align}\label{basicestimate}
|\rho_0(z)^2-b^2|+|\partial_z^l\rho_0(z)|\leq C_le^{-\al_0 |z|}&,\  \forall z\in (0,+\infty) ,\nonumber\\
|\rho_0(z)^2-a^2|+|\partial_z^l\rho_0(z)|\leq C_le^{-\al_0 |z|}&,\  \forall z\in (-\infty,0),
\end{align}
where
$$0<\alpha_0 <\min \{\sqrt{2}(b^2-a^2)b,\sqrt{2}(b^2-a^2)a\}.$$
\proof
We place the proof in  Appendix \ref{rho0zpro}.
\qed
\end{lemma}
Without loss of generality, we choose $\rho_0$ to be the solution of
\eqref{rho0expression} with $c_0=0$.

However, for $(x,t)\in\Gamma(\delta)$, the relation $\omega^+=\omega^-$
does not hold in general. The reason is that $\omega^\pm|_{U_\pm\cap \Gamma(\delta)}$ is obtained by smooth extension of $\omega^\pm|_{U_\pm}$, which is solved by the outer expansion systems. However, the normal derivatives of $\partial_\nu\omega^\pm$ on $\Gamma$ do not coincide with each other due to the jump condition $a^2\partial_\nu\omega^-=b^2\partial_\nu\omega^+$. In this case, the
\eqref{-2term}-\eqref{eq:bd-0} has no solution in general, which motivates us to apply the idea of {\it quasi-minimal connecting orbit} in \cite{FLWZ}.

For  $\omega^+(x,t)\neq \omega^-(x,t)\in \mathbb{S}^{n-1}$, there is a
 geodesic ${\omega}_*(\tau,x,t)(\tau\in[0,1])$ on $\BS^{n-1}$ such that
\begin{align}\label{const.}
{\omega}_*(0,x,t)=\omega^-(x,t),\qquad {\omega}_*(1,x,t)&=\omega^+(x,t),\qquad
|\pa_\tau{\omega}_*(\tau,x,t)|=\text{const}.\text{ for }\tau\in (0,1).
\end{align}
Clearly, $\omega_*(\tau, x, t)$ are smooth in $[0,1]\times \Gamma(\delta)$. As $\omega_*(\tau, x, t)|_{\Gamma}\equiv\omega^\pm(x,t)$, one has
$$\partial_\tau\omega_*(\tau, x, t),\quad \partial_\tau^2\omega_*(\tau, x, t)=O(d_0).$$

Then we choose a smooth orthogonal basis of $\xi_{\alpha}^-$ on $\Gamma(\delta)$. Let $\xi_{*\alpha}(\tau, x, t)(0\le\tau\le 1)$ be the tangential vector which is generated by parallel transporting $\xi_\alpha^-(x, t)$ along the geodesic $\omega^*(\tau, x, t)(0\le \tau\le 1)$ on the unit sphere, namely, the covariant derivative
\begin{align*}
\frac{D}{D\tau}\xi_{*\alpha}(\tau, x, t)=0,\ \ \text{with }\ \ \xi_{*\alpha}(0, x, t)=\xi_\alpha^-(x, t).
\end{align*}
 Then we choose $\xi_{\alpha}^+(x, t)=\xi_{*\alpha}(1, x, t)$. Moreover, we have
\begin{align*}
  \xi_{*\alpha}(\tau, x, t)\cdot \xi_{*\beta}(\tau, x, t) =\delta_{\alpha\beta},\quad \xi_{*\alpha}(\tau, x, t)\cdot\omega_*(\tau, x, t)=0.
\end{align*}
Lemma \ref{lem:frame-product} gives us that $\partial_\tau\xi_{*\alpha}=\kappa_\alpha\omega_*(\tau,x,t)$, where $\kappa_\alpha=\kappa_\alpha(x,t)$ is independent of $\tau$. We also denote $\kappa_0(x,t)=|\partial_\tau\omega_*(\tau,x,t)|$.
In addition, we get from $\pa_\tau\omega_*\cdot \xi_{*\al}=-\pa_\tau \xi_{*\al}\cdot \omega_*(\tau,x,t)=-\kappa_\al$ that
\begin{align}\label{kappa0kappaal}
\kappa_0^2=|\pa_\tau \omega_*|^2=\Big|\sum_\al \kappa_\al \xi_{*\al}\Big|^2=\sum_{\al=1}^{n-1} \kappa_\al ^2.
\end{align}

Introduce
$$\bar\omega(z,x,t)=\omega_*(\eta(z),x,t),\quad \xi_{\alpha}(z, x,t)=\xi_{*\alpha}(\eta(z),x,t),$$
where $\eta$ is a smooth and monotone function such that $\eta(-\infty)=0$, $\eta(+\infty)=1$ and
\begin{align}\label{def:eta1}
  |\partial_z^k\eta(z)|\le C_k e^{\alpha_0z}, \text{ as }z\to-\infty;\quad
   |\partial_z^k(\eta(z)-1)|\le C_k e^{-\alpha_0z}, \text{ as }z\to+\infty.
\end{align}
For example, we can take $\eta(z)=\frac{\rho_0(z)-a}{b-a}$.
Then $\bar\omega(\pm\infty, x, t)=\omega^\pm(x,t)$ and $\xi_{\alpha}(\pm\infty, x, t)=\xi_{\alpha}^\pm(x,t)$. Moreover, ($\partial_z =\eta'(z)\partial_\tau$)
\begin{align}\label{pazw}
 & \bar{\omega}(z, x,t)=\omega^+(x,t)=\omega^-(x,t), \quad  \partial_z\bar{\omega}=0, \quad\text{for }(x,t)\in\Gamma,\\ \label{pazxigamma}
 & \xi_{\alpha}(z, x,t)=\xi_{\alpha}^+(x,t)=\xi_{\alpha}^-(x,t), \quad  \partial_z\xi_{\alpha}=0, \quad\text{for }(x,t)\in\Gamma,
 \end{align}
 and for $(x,t)\in\G(\delta)$,
 \begin{align}
\label{decay:bar-omega}
\partial_z^2\bar{\omega}(z, x,t), \ \partial_z\bar{\omega}(z, x,t),\  \partial_z^2\xi_{\alpha}(z, x,t), \ \partial_z\xi_{\alpha}(z, x,t)=O(d_0e^{-\al_0|z|}).
\end{align}
Here are some important equalities:
\begin{align}
\partial_z\bar\omega\cdot\xi_{\alpha}&=-\bar\omega\cdot\partial_z\xi_{\alpha}=-\bar\omega\cdot\partial_\tau\xi_{\alpha}\eta'=-\kappa_\alpha(x,t)\eta';\notag\\
\partial_z\bar\omega\cdot\partial_z\xi_{\alpha}&=\partial_\tau\bar\omega\cdot\partial_\tau\xi_{\alpha}(\eta')^2=0,\notag\\
\partial_z^2\bar\omega\cdot\bar\omega&=-\partial_z\bar\omega\cdot\partial_z\bar\omega=-\kappa_0^2(x,t)(\eta')^2,\notag\\
\partial_z^2\bar\omega\cdot\xi_{\alpha}&=\partial_z(\partial_z\bar\omega\cdot\xi_{\alpha})-\partial_z\bar\omega\cdot\partial_z\xi_{\alpha}=-\kappa_\alpha(x,t)\eta'',\label{basic computation}\\
\partial_z^2\xi_{\alpha}\cdot\bar\omega&=\partial_z(\partial_z\xi_{\alpha}\cdot\bar\omega)-\partial_z\bar\omega\cdot\partial_z\xi_{\alpha}
=\kappa_\alpha(x,t)\eta'',\notag\\
\partial_z\xi_{\alpha}\cdot\xi_{\beta}&=\kappa_\alpha\eta' \bar\omega \cdot\xi_\beta=0;\notag\\
\partial_z^2\xi_{\alpha}\cdot\xi_{\beta}&=\partial_z(\partial_z\xi_{\alpha}\cdot\xi_\beta)-\partial_z\xi_{\alpha}\cdot\partial_z\xi_{\beta}
=-\kappa_\alpha(x,t)\kappa_\beta(x,t)(\eta')^2.\notag
\end{align}
Define
\begin{align}\label{solution:profile}
  u_0(z,x,t)=\bar{\omega}(z,x,t) \rho_0(z).
\end{align}
Thus, one has for $(x,t)\in\G(\delta)$
\begin{align}\nonumber
\partial_{zz}u_0-f(u_0)&=\bar{\omega}\big(\partial_z^2 \rho_0-f(\rho_0)\big)+\partial_z^2\bar{\omega}\rho_0+2\partial_z\bar{\omega}\rho_0'\\
&=\partial_z^2\bar{\omega}\rho_0+2\partial_z\bar{\omega}\rho_0'\nonumber\\
&=O(d_0e^{-\al_0|z|}), \label{eq:0remainder}
\end{align}
due to \eqref{decay:bar-omega}.

The profile $u_0$ defined in \eqref{solution:profile} is called a {\it quasi-minimal connecting orbit}, which enables us to construct approximate solutions to a modified (but equivalent) system.
Namely, as in \cite{ABC1994, FLWZ}, we consider
\begin{align}\label{modify}
   \partial_t {u}^\ve=\Delta {u}^\ve-\frac{1}{\ve^2} f({u}^\ve)- (d^\ve -\ve z)g^\ve(z,x,t)-H^{+,\ve}(x,t)\eta_M^+(z)-H^{-,\ve}(x,t)\eta_M^-(z),
\end{align}
where $$g^\ve(z,x,t)=\sum_{k=0}^{+\infty}\ve^{k-2}g_k(z,x,t),\quad H^{\pm,\ve}=\sum_{k=0}^{
+\infty
}\ve^kH_k^{\pm}(x,t)$$
are vector-valued functions to be determined step by step. We choose a fixed smooth and nonnegative function
$\eta_0(z)$ satisfying: $\eta_0(z)=0$ if $z\le -1$, $\eta_0(z)=1$ if $z\ge 1$, $\eta_0'(z)\ge 0$. Set
$$\eta_M^\pm(z)=\eta_0(-M\pm z),\text{ with }M=\|d_1\|_{C^0(\G)}+2.$$
Then one has $\eta_0(-M\pm d^\ve(x,t)/\ve)=0$ for $(x,t)\in\G(\delta)\cap U_{\mp}.$ By \eqref{hkpm}, we obtain $H^{\pm,\ve}(x,t)=0$ for $(x,t)\in \Gamma(\delta)\cap U_\pm$. Thus,  one has
\begin{align*}
  H^{+,\ve}(x,t)\eta_M^+(z)+H^{-,\ve}(x,t)\eta_M^-(z)\Big|_{z=d^\ve(x,t)/\ve}=0,\quad \text{ for any }(x,t)\in\Gamma(\delta).
\end{align*}
As a result, all the modified terms vanish on $\{d^\ve =\ve z\}$, and will not change the system (\ref{eq:main-inner}).

Substituting \eqref{innerexpansion} into \eqref{modify}, we obtain
\begin{align*}
0=&\sum_{k=0}^{+\infty}\ve^k\Big(\partial_tu_k+\frac1{\ve}\partial_z u_k\partial_t \sum_{i=0}^{+\infty}\ve^i d_i\Big)\nonumber\\
&-\sum_{k=0}^{+\infty}\ve^k\Big\{\frac{1}{\ve^2}\partial_{zz}u_k+\frac{1}{\ve}\Big(2\partial_z\nabla u_k\cdot\nabla \sum_{i=0}^{+\infty}\ve^i d_i
+\partial_zu_k\Delta \sum_{i=0}^{+\infty}\ve^i d_i\Big)+\Delta u_k\Big\}\nonumber\\
&+\frac{1}{\ve^2}f\Big(\sum_{k=0}^{+\infty}\ve^ku_k\Big)+\Big(\sum_{k=0}^{+\infty}\ve^k d_k -\ve z\Big)\sum_{k=0}^{+\infty}\ve^{k-2}g_k+\sum_{k=0}^\infty \ve ^k\Big(H^{+}_k\eta_M^++H^{-}_k\eta_M^-\Big) .
\end{align*}
Then we arrive at the following systems of order $\ve^{k-2}(k\geq 0)$:
\begin{itemize}
\item The $O(\ve^{-2})$ system takes the form
\begin{align}\label{-2-2}
\partial_{zz}u_0-f(u_0)=g_0d_0.
\end{align}
\item The $O(\ve^{k-2})(k\geq 1)$ system takes the form
\begin{align}\label{innerk-2}
\partial_{zz}u_k-Df(u_0)u_k=&\mathcal{D}_k+g_kd_0.
\end{align}
\begin{itemize}[label=$\vartriangle$]
\item
For $k=1,$ we have
\begin{align}
\mathcal{D}_1=&\partial_zu_0(\partial_td_0-\Delta d_0)-2\nabla d_0\cdot\nabla\partial_z u_0 +g_0(d_1-z).\label{D1}
\end{align}
\item For $k=2$, one has
\begin{align}
\mathcal{D}_2=&\partial_zu_0(\partial_td_{1}-\Delta d_{1})+\partial_zu_{1}(\partial_td_{0}-\Delta d_{0})-2\nabla d_{1}\cdot \nabla \partial_zu_0-2\nabla d_{0}\cdot \nabla \partial_zu_{1}\nonumber\\
&+f^{(1)}(u_0,u_{1})+d_2g_0+(d_1-z)g_{1}+\partial_{t}u_{0}-\Delta u_{0}+H_0^+\eta_M^++H_0^-\eta_M^-.\label{D2}
\end{align}
\item For $k\geq 3$, we deduce that
\begin{align}
\mathcal{D}_k=&\partial_zu_0(\partial_td_{k-1}-\Delta d_{k-1})+\partial_zu_{k-1}(\partial_td_{0}-\Delta d_{0})-2\nabla d_{k-1}\cdot \nabla \partial_zu_0\nonumber\\
&-2\nabla d_{0}\cdot \nabla \partial_zu_{k-1}+f^{(k-1)}(u_0, \cdots, u_{k-1})+(d_kg_0+d_{k-1}g_1)\nonumber\\
&+(d_1-z)g_{k-1}+\mathcal{A}_{k-2}, \label{Dk}
\end{align}
with
\begin{align*}
\mathcal{A}_{k-2}=&\partial_{t}u_{k-2}-\Delta u_{k-2}+\sum_{\substack{i+j=k-1,\\i,j<k-1}}\partial_zu_i\partial_td_j-\sum_{\substack{i+j=k-1,\\i,j<k-1}}\partial_zu_i\Delta d_j\\
&-2\sum_{\substack{i+j=k-1,\\i,j<k-1}}\nabla d_j\cdot \nabla \partial_zu_i +\sum_{\substack{i+j=k,\\i,j<k-1}}d_ig_j+H_{k-2}^+\eta_M^++H_{k-2}^-\eta_M^-.
\end{align*}
\end{itemize}
\end{itemize}
Set $P^\pm=I-\omega^\pm\otimes\omega^\pm$.
We can choose
\begin{align}\label{hkpm}
H_k^\pm&=-P^\pm\Big((\pa_t-\Delta)u_k^\pm+f^{k+1}(u_0^\pm,\cdots,u_{k+1}^\pm)\Big).
\end{align}
By \eqref{outuk-2} and \eqref{de2}, we have $H^{\pm,\ve}=0$ for $(x,t)\in \G(\delta)\cap U_\pm.$

We can extend $\rho_{k+1}^\pm$ to
$\Gamma(\delta)$ by \eqref{outersolve}.
Since $u_0^\pm\cdot u_1^\pm=0$,
\eqref{fk-1} gives, for $k\ge1$,
\[
H_k^\pm
=
-2\rho_0^\pm(b^2-a^2)^2\rho_{k+1}^\pm u_1^\pm
-P^\pm\left[
(\partial_t-\Delta)u_k^\pm
+f_4^{(k+1)}(u_0^\pm,\ldots,u_k^\pm)
\right].
\]
By \eqref{outersolve}, we know that $\rho_{k+1}^\pm$ depends only on $u_0^\pm,\cdots,u_k^\pm.$ Together with Lemma \ref{fex}, we arrive at $H_k^\pm$ depends only on $u_0^\pm,\ldots,u_k^\pm$. For $k=0$, $H_0^\pm
=-\rho_0^\pm P^\pm(\partial_t\omega^\pm-\Delta\omega^\pm)$
depends only on $u_0^\pm$. In conclusion, we derive that $H_k^\pm$ depends only on $u_0^\pm,\ldots,u_k^\pm$.

Our aim is to find a family of functions $u_k,d_k$ and $g_k$ for $k\geq 0$ to construct an approximate solution. For convenience, we take for $k\geq 0,$
$$\mathcal{V}^k=(u_k^\pm(x,t),u_k(z,x,t),d_k(x,t),g_k(z,x,t)).$$

\subsection{The leading order system}
To solve the leading order system \eqref{-2-2}, we choose $u_0$ as in \eqref{solution:profile}:
\begin{align}\label{u0}
u_0(z,x,t)=\bar{\omega}\big(z,x,t\big) \rho_0(z).
\end{align}
From \eqref{eq:0remainder}, we may take
\begin{align}
g_0(z,x,t)=\frac{1}{d_0}\big(2\partial_z\bar{\omega}(z, x, t)\rho_0'+\partial_z^2\bar{\omega}(z, x, t)\rho_0\big).\label{g0}
\end{align}
Then \eqref{-2-2} is satisfied for $(x,t)\in\Gamma(\delta)$.
In particular, we have by \eqref{basic computation} that
\begin{align}\label{g0omega}
g_0(z,x,t)\cdot\bar{\omega}=\frac{\rho_0}{d_0}\partial_z^2\bar{\omega}\cdot\bar{\omega}=-\frac{\rho_0}{d_0}\kappa_0^2(x,t)(\eta'(z))^2,
\end{align}
and $(g_0\cdot\bar{\omega})|_{\Gamma}=0$ as $\kappa_0=O(d_0)$. In addition, we have by \eqref{basic computation} that
\begin{align}
g_0(z,x,t)\cdot\xi_{\alpha}&= \frac{1}{d_0\rho_0}\pa_z(\rho_0^2\partial_z\bar{\om})\cdot\xi_\alpha
 =\frac{1}{d_0\rho_0}\pa_z(\rho_0^2\partial_z\bar{\om}\cdot\xi_\alpha\big)-\frac{\rho_0}{d_0}\partial_z\bar{\om}\cdot\partial_z\xi_\alpha\notag\\
 &=\frac{1}{d_0\rho_0}\pa_z(\rho_0^2\partial_z\bar{\om}\cdot\xi_\alpha\big)
 =-\frac{\kappa_\alpha}{d_0\rho_0}\pa_z(\rho_0^2\eta'\big).
 \label{g0xibeta}
\end{align}
\subsection{The analysis of the linearized system for next orders}
Now we analyze the equation \eqref{innerk-2} for $k\ge 1$.

For
\begin{align*}
  u(z,x,t)=\rho(z,x,t)\bar{\omega}(z,x,t)+\sum_\al \sigma_\al(z,x,t) \xi_\al(z,x,t),
\end{align*}
we have
\begin{align*}
\partial_{zz}u-Df(u_0)u=&\Big(\big(\partial_z^2\rho-f_A(\rho_0)\rho\big)\bar\omega  +\sum_\al\big(\partial_z^2\sigma_{\al}-f_B(\rho_0)\sigma_{\al}\big)\xi_\alpha\Big)\notag\\
&+\Big(2\partial_z\rho\partial_z\bar\omega+\rho\partial_z^2\bar\omega+\sum_\al 2\partial_z\sigma_{\al}
\partial_z\xi_\alpha+\sigma_{\al}\partial_z^2\xi_\alpha\Big),
\end{align*}
where
\begin{align}
f_A(s)= & (s^2-a^2)(s^2-b^2)(2s^2-a^2-b^2)
+2\big((2s^2-a^2-b^2)^2+2(s^2-a^2)(s^2-b^2)\big)s^2 ,\label{fA'}\\
f_B(s)= & (s^2-a^2)(s^2-b^2)(2s^2-a^2-b^2).\label{fB'}
\end{align}
Using \eqref{1D:profile}, \eqref{fA'} and \eqref{fB'}, we obtain the identities
\begin{align}
  f_A(\rho_0(z))=\frac{\rho_0'''}{\rho_0'}\ \text{and}\ f_B(\rho_0(z))=\frac{\rho_0''}{\rho_0},\label{fAfB}
\end{align}
We introduce
\begin{align}
  \mathcal{L}_*u=\Big(\big(\partial_z^2\rho-f_A(\rho_0)\rho\big)\bar\omega  +\sum_\al\big(\partial_z^2\sigma_{\al}-f_B(\rho_0)\sigma_{\al}\big)\xi_\alpha\Big).
\end{align}
Then $\pa_{zz}u-Df(u_0)u=\mathcal{L}_*u$ for $(x,t)\in\Gamma$. Thus, we can define for $(x,t)\in\Gamma(\delta)$:
\begin{align}\label{def:Du} \nonumber
\mathcal{D}_{*}u:=&~\frac{1}{d_0}\Big(\partial_{zz}u-\partial_z^2\rho\bar\omega -\sum_\al\partial_z^2\sigma_{\al}\xi_\alpha\Big)\\
=&~~\frac{1}{d_0}\Big(2\partial_z\bar\omega(z,x,t)\rho'+\partial_z^2\bar\omega(z,x,t)\rho +2\sum_\al \partial_z\xi_\alpha(z,x,t)\sigma_{\alpha}'+\sum_\al\partial_z^2\xi_\alpha(z,x,t)\sigma_{\alpha}\Big).
\end{align}
Then
\begin{align}
  \pa_{zz}u-Df(u_0)u=  \mathcal{L}_*u+d_0\mathcal{D}_{*}u.
\end{align}
Moreover, by \eqref{basic computation}, we have
\begin{align}
\mathcal{D}_{*}u\cdot\bar\omega= & \frac{1}{d_0}\Big(-\kappa_0^2(\eta')^2\rho+\sum_\al(2\kappa_\alpha\eta'\sigma_{\alpha}'+\kappa_\alpha\eta''\sigma_{\alpha})\Big),\notag\\
\mathcal{D}_{*}u\cdot\xi_\beta=& - \frac{\kappa_\beta}{d_0}\Big(2\eta'\rho'+\eta''\rho+\sum_\al\kappa_\alpha(\eta')^2\sigma_{\alpha}\Big).\label{d*}
\end{align}

Now, we construct
\begin{align}\label{gk}
    g_k(z,x,t)=&\tilde{g}_k(z,x,t)+\mathcal{D}_{*}u_k,
\end{align}
where
\begin{align*}
\tilde{g}_k(z,x,t)=  L_{k,0}(x,t)\eta'(z)\bar\omega+\sum_\al L_{k,\al}(x,t)\rho_0^{-1}\eta'(z)\xi_{\alpha}+\sum_\al h_{k,\alpha}(x,t)\Big(2\rho_0'\eta'+\rho_0\eta'' \Big)\xi_{\alpha}.
\end{align*}
Here $L_{k,0}(x,t), L_{k,\al}(x,t), h_{k,\alpha}(x,t)$ are functions  determined later.
Then the system \eqref{innerk-2} is reduced to
\begin{align}\label{uk-simple}
\mathcal{L}_{*}u_{k}=\mathcal{D}_{k}+d_0\tilde{g}_k.
\end{align}
If we introduce the decompositions for $u_k$:
\begin{align}\label{ukdecomp}
u_k(z,x,t)&~=\rho_k(z,x,t)\bar\omega(z,x,t)+\sum_{\alpha}\sigma_{k,\al}(z,x,t)\xi_\alpha(z,x,t),
\end{align}
the system \eqref{uk-simple} can be decomposed as:
\begin{align}
\partial_z^2\rho_k-f_A(\rho_0)\rho_k=&~d_0L_{k,0}\eta'+\mathcal{D}_k\cdot\omega,\label{rhol}\\
\Big(\partial_z^2-f_B(\rho_0)\Big)\Big(\sigma_{k, \beta}-h_{k,\beta}(x,t)d_0\rho_0\eta\Big)=&~d_0L_{k, \beta}\rho_0^{-1}\eta'+\mathcal{D}_k\cdot\xi_{\beta},\label{sigmal}
\end{align}
with the boundary conditions
\begin{align*}
\rho_k(\pm\infty,x,t)=\rho_k^\pm(x,t),\quad
\sigma_{k,\beta}(\pm\infty,x,t)=\sigma_{k,\beta}^\pm(x,t).
 \end{align*}
Here, we have used the fact that
\begin{align*}
  \Big(\partial_z^2-f_B(\rho_0)\Big)(\rho_0(z)\eta(z))=2\rho_0'\eta'+\rho_0\eta''.
\end{align*}

In light of \eqref{rhol} and \eqref{sigmal}, it is important to give the solvability conditions for the following systems
\begin{align}
\partial_z^2v_1-f_A(\rho_0(z)) v_1=h_1(z,x,t), \label{ODE:A}\\
\partial_z^2v_2-f_B(\rho_0(z)) v_2=h_2(z,x,t). \label{ODE:B}
\end{align}
By \eqref{fAfB}, we can rewrite \eqref{ODE:A} and \eqref{ODE:B} as
\begin{align}
\partial_z^2v_1-f_A(\rho_0(z)) v_1&=\frac{1}{\rho_0'}\partial_z\Big((\rho_0')^2\partial_z\Big(\frac{v_1}{\rho_0'}\Big)\Big),\label{v1eq}\\
\partial_z^2v_2-f_B(\rho_0(z)) v_2&=\frac{1}{\rho_0}\partial_z\Big(\rho_0^2\partial_z\Big(\frac{v_2}{\rho_0}\Big)\Big).\label{v2eq}
\end{align}
\begin{lemma}\label{v1product}
    Let $v_1$ be the bounded solution to \eqref{ODE:A}, then we have
   \begin{align}\label{energy:v1}
-2\int_{\mathbb{R}}\rho_0'(z)v_1'(z)dz=\int_{\mathbb{R}}z\rho_0'(z)h_1(z)dz.
\end{align}
\end{lemma}
\begin{proof}
By \eqref{v1eq},  the equation \eqref{ODE:A} can be equivalently written as
\begin{align*}
  \frac{1}{\rho_0'}\partial_z\Big((\rho_0')^2\partial_z\Big(\frac{v_1}{\rho_0'}\Big)\Big)=h_1(z).
\end{align*}
Multiplying both sides by $z\rho_0'$ to get
\begin{align*}
\int_{\mathbb{R}}z\rho_0'(z)h_1(z)\ud z=&\int_{\mathbb{R}}z\partial_z\Big((\rho_0')^2\partial_z\Big(\frac{v_1}{\rho_0'}\Big)\Big)\ud z\\
=&\int_{\mathbb{R}}z\partial_z(\rho_0'v_1'-\rho_0''v_1)\ud z\\
=&-\int_{\mathbb{R}}(\rho_0'v_1'-\rho_0''v_1)\ud z\\
=&-2\int_{\mathbb{R}}\rho_0'v_1'\ud z.
\end{align*}
\end{proof}
Let $$\theta_1(s)=\rho_0'(s),\  \theta_2(s)=\rho_0(s) \text{ and } \alpha_0\in (0,\min \{\sqrt{2}(b^2-a^2)b,\sqrt{2}(b^2-a^2)a\}].$$
Define the spaces as:
\begin{align*}
\mathcal{S}_{J,L,M}(\alpha_0,k)=\{&f(\cdot,x,t)\in C^j(\mathbb{R}):\ \forall\  (j,l,m)\in[0,J]\times [0,L]\times[0,M],\\
&|\partial_z^j\partial_x^l\partial_t^m(f(z,x,t)-f^\pm(x,t))|\lesssim |z|^ke^{-\alpha_0|z|},\ \text{as}\ z\to \pm\infty\}.
\end{align*}
The following lemma gives the compatibility conditions for \eqref{ODE:A} and \eqref{ODE:B}.
\begin{lemma}\label{lem:s1}
(Compatibility conditions) Assume $h_i(\cdot,x,t)\in \mathcal{S}_{J,L,M}(\alpha_0,k)$ for $i=1,2$ such that
\begin{align}\label{solvablecondition}
  \int_{\mathbb{R}}h_i(z,x,t)\theta_i(z)\ud z =0 \quad (i=1,2),
\end{align}
and
\begin{align}\label{BoundaryCond}
~~h_2^\pm(x,t)=0.
\end{align}
Then \eqref{ODE:A} or \eqref{ODE:B} has a  unique bounded solution
\begin{equation}\label{u1u2star}
v^*_{i}(\cdot, x,t)\in \mathcal{S}_{J+2,L,M}(\alpha_0,k+1),
\end{equation}
satisfying
\begin{align*}
v_1^*(0,x,t)=0,\ v_2^*(-\infty,x,t)=0.
\end{align*}
Precisely, the solution can be written as
\begin{align}
v^*_{1}(z,x,t)&=-\theta_1(z)\int_{0}^{z}\theta_1^{-2}(\varsigma)\int_\varsigma^{+\infty}h_1(\tau,x,t)\theta_1(\tau)\ud \tau \ud \varsigma,\\
v^*_{2}(z,x,t)&=\theta_2(z)\int^z_{-\infty}\theta_2^{-2}(\varsigma)\int_{-\infty}^{\varsigma}h_2(\tau,x,t)\theta_2(\tau)\ud \tau \ud \varsigma,\label{v2solution}
\end{align}
and we have
\begin{align*}
 v_1^{*}(-\infty, x,t)=-\frac{h_1(-\infty, x,t)}{2(a^2-b^2)^2a^2}, \quad v_1^{*}(+\infty, x,t)=-\frac{h_1(+\infty, x,t)}{2(a^2-b^2)^2b^2}.
\end{align*}
Moreover, all bounded solutions of \eqref{ODE:A} and \eqref{ODE:B} are given by
\begin{align*}
  v_i(z,x,t)=  v_i^*(z,x,t)+q_i(x,t)\theta_i(z),\quad \text{for}\ i=1,2.
\end{align*}
\proof
The proof is provided in Appendix \ref{prooflemma3.1}.
\qed
\end{lemma}
\subsection{Deriving equations for the $O(\ve^{-1})$ system}
Setting $k=1$ in \eqref{rhol} and \eqref{sigmal}, we obtain
\begin{equation}\label{rho1equation}
\left\{
\begin{array}{ll}
 \partial_z^2\rho_1-f_A(\rho_0)\rho_1=L_{1, 0}\eta'd_0+\mathcal{D}_1\cdot\bar{\omega}, \\
 \rho_1(\pm\infty, x, t)=0,
\end{array}
\right.
\end{equation}
and
\begin{equation}\label{sigma1eq}
\left\{
\begin{array}{ll}
\big(\partial_z^2-f_B(\rho_0)\big)\big(\sigma_{1, \beta}-h_{1,\beta}d_0\rho_0\eta\big)=L_{1, \beta}\rho_0^{-1}\eta'd_0+\mathcal{D}_1\cdot\xi_\beta, \\
\sigma_{1, \beta}(\pm\infty, x, t)=\sigma_{1, \beta}^\pm,
\end{array}
\right.
\end{equation}
where
\begin{align*}
\mathcal{D}_1=&\partial_zu_0(\partial_td_0-\Delta d_0)-2\nabla d_0\cdot\nabla\partial_z u_0 +g_0(d_1-z).
\end{align*}
We can directly check that the right hand side of above equations exponentially tends to its values
at $\pm\infty$. On $\G,$ for \eqref{rho1equation}, according to the compatibility conditions \eqref{solvablecondition},
we have
\begin{align}\label{inner:m1-cancel1}
0=&\langle \mathcal{D}_1\cdot \bar{\omega},\rho_0'\rangle.
\end{align}
By \eqref{g0omega}, we have
\begin{align}\label{d1omega}
   \mathcal{D}_1\cdot\bar{\omega}= \rho_0'(\partial_td_0-\Delta d_0)-2\rho_0\nabla d_0\cdot\nabla\pa_z
 \bar\omega \cdot\bar\omega-\frac{\rho_0}{d_0}\kappa_0^2(x,t)(\eta'(z))^2(d_1-z)
\end{align}
It follows from $\partial_\nu\bar{\omega}\bot\bar{\omega}$ that
\begin{align}
  \partial_\nu(\partial_z\bar{\omega})=\partial_z\partial_\nu\bar{\omega}\bot\bar\omega,\quad
  \partial_\nu(\partial_z^2\bar{\omega})=\partial_z^2\partial_\nu\bar{\omega}\bot\bar\omega\quad \text{on}\G.\label{bot} \end{align}
Then  we obtain
\begin{align}\label{d0}
  \partial_td_0-\Delta d_0=0,\qquad \text{on } \Gamma,
\end{align}
which implies that $\Gamma$ evolves according to the mean curvature flow. In $\Gamma(\delta)$, by \eqref{nablad0nabladk}, we have the closed system about $d_0$:
\begin{equation}\label{d0delta}
\left\{
\begin{array}{ll}
\partial_t d_0-\Delta d_0=0, &\text{on} \ \Gamma;\\
d_0=0, &\text{on}\ \Gamma;\\
|\nabla d_0|^2=1, &\text{in}\ \Gamma(\delta).
\end{array}
\right.
\end{equation}

On the other hand, for \eqref{sigma1eq},  from the fact \eqref{pazxigamma},  the compatibility condition \eqref{solvablecondition}  leads to
\begin{align*}
0=&\langle \mathcal{D}_1\cdot\xi_\beta,\rho_0\rangle\\
=&\langle -2\rho_0'\nabla d_0\cdot\nabla\bar\omega\cdot\xi_\beta-2\rho_0\nabla d_0\cdot\nabla\partial_z\bar\omega\cdot\xi_\beta +g_0\cdot\xi_{\beta}(d_1-z),\rho_0\rangle_\G\notag\\
=&\langle -2\rho_0'\pa_\nu\bar\omega\cdot\xi_\beta-2\rho_0\partial_{\nu z}\bar\omega\cdot\xi_\beta +(\pa_z^2\pa_\nu\bar{\omega}\rho_0+2\pa_{\nu z}\bar{\omega}\pa_z\rho_0)\cdot\xi_{\beta}(d_1-z),\rho_0\rangle\notag\\
=&-\langle\partial_z(\rho_0^2\pa_\nu\bar\omega)\cdot\xi_\beta,1\rangle-\langle\rho_0^2,\partial_{\nu z}\bar\omega\cdot\xi_\beta \rangle +\langle \pa_z(\partial_z\pa_\nu\bar{\omega}\rho_0^2)\cdot\xi_{\beta},d_1-z\rangle\\
=&-\rho_0^2\pa_\nu\bar{\omega}\cdot \xi_\beta|_{-\infty}^{+\infty}+\langle\pa_z\big(\pa_{z\nu}\bar{\omega}\rho_0^2\cdot \xi_\beta(d_1-z)\big),1\rangle\\
=&-\rho_0^2\pa_\nu\bar{\omega}\cdot \xi_\beta|_{-\infty}^{+\infty},
\end{align*}
which implies that on $\G$, $\omega^\pm$ satisfies a jump condition
 \begin{align}\label{jumpcondition}
b^2\pa_\nu\omega^+=a^2\pa_\nu\omega^-.
\end{align}

In $\G(\delta),$ the equation \eqref{rho1equation} has a bounded solution if and only if
\begin{align}
\langle L_{1, 0}\eta'd_0+\mathcal{D}_1\cdot\bar{\omega}, \rho_0'(z)\rangle=0.
\end{align}
Together with \eqref{d1omega}, we have
\begin{align*}
L_{1, 0}d_0\langle\eta', \rho_0'\rangle+\Big\langle\rho_0'(\partial_td_0-\Delta d_0)-2\rho_0\nabla d_0\cdot\nabla\pa_z
 \bar\omega \cdot\bar\omega-\frac{\rho_0}{d_0}\kappa_0^2(x,t)(\eta'(z))^2(d_1-z), \rho_0'\Big\rangle=0,
\end{align*}
 Thus we can obtain
\begin{align}\nonumber
L_{1,0}=&
\Big(d_0\langle \rho_0',\eta'\rangle\Big)^{-1}\Big((\Delta d_0-\pa_td_0)e+2\langle \rho_0(\nabla d_0\cdot\nabla)\pa_z
 \bar\omega \cdot\bar\omega,\rho_0'\rangle+\frac{\kappa_0^2}{d_0}\langle\rho_0(\eta')^2(d_1-z),\rho_0' \rangle\Big)\\
=&\Big(d_0\langle\rho_0',\eta'\rangle\Big)^{-1}\Big((\Delta d_0-\pa_td_0)e+2\langle\rho_0(\nabla d_0\cdot\nabla)\pa_z
 \bar\omega \cdot\bar\omega,\rho_0'\rangle-\frac{\kappa_0^2}{d_0}\langle \rho_0(\eta')^2z,\rho_0'\rangle\nonumber\\ \nonumber
 &\qquad\qquad \qquad \qquad +\frac{\kappa_0^2d_1}{d_0}\langle\rho_0(\eta')^2,\rho_0'\rangle\Big)\\
 =&\frac{d_1\kappa_0^2}{d_0^2}c_1^*+L_{1,0}^*,\label{L10}
\end{align}
where $c_1^*=\langle \rho_0(\eta')^2,\rho_0'\rangle/\langle\rho_0',\eta'\rangle$ and $L_{1,0}^*=L_{1,0}^*(d_0, \rho_0,\bar\omega, \eta)$.
Then by solving \eqref{rho1equation}, we get
\begin{align}\nonumber
\rho_{1}(z,x,t)&=-\rho_0'(z)\int_{0}^{z}\rho_0'(\varsigma)^{-2}\int_\varsigma^{+\infty}(L_{1,0}\eta'd_0+\mathcal{D}_1\cdot\bar\omega)\rho_0'(\tau)\ud \tau \ud \varsigma\\
&=-\frac{\kappa_0^2d_1}{d_0}\rho_0'(z)\int_{0}^{z}\rho_0'(\varsigma)^{-2}\int_\varsigma^{+\infty}
\big(c_1^*\eta'-\rho_0(\eta')^2\big)\rho_0'\ud s \ud \varsigma+\rho_1^*(z,x,t),\label{rho1}
\end{align}
where $\rho_1^*$ depends only on $d_0, \rho_0,\bar\omega$ and $\eta$.

For \eqref{sigma1eq}, according to the compatibility conditions \eqref{solvablecondition},  one has
\begin{align*}
\langle L_{1,\beta}\rho_0^{-1}\eta' d_0+\mathcal{D}_1\cdot\xi_\beta,\rho_0\rangle=0.
\end{align*}
Direct calculation gives us
\begin{align}\label{d1xibeta}
\mathcal{D}_1\cdot\xi_\beta=\rho_0\partial_z\bar\omega\cdot\xi_\beta(\partial_td_0-\Delta d_0)
-2\rho_0'\nabla d_0\cdot\nabla\bar\omega\cdot\xi_\beta-2\rho_0\nabla d_0\cdot\nabla\partial_z\bar\omega\cdot\xi_\beta +g_0\cdot\xi_{\beta}(d_1-z).
\end{align}
Using \eqref{g0xibeta}, we have $\int_{\mathbb{R}}\rho_0g_0\cdot\xi_{\beta}\ud z=0$ and
\begin{align*}
    \langle z g_0\cdot \xi_\beta ,\rho_0 \rangle=\Big\langle z\frac{1}{d_0\rho_0}\pa_z(\rho_0^2\pa_z\bar{\omega}\cdot\xi_\beta),\rho_0\Big\rangle=-\frac{1}{d_0}\langle \pa_z\bar{\omega}\cdot\xi_\beta,\rho_0^2\rangle.
\end{align*}
Thus,
\begin{align}\nonumber
L_{1,\beta}=&
-\frac{1}{d_0}\Big(\langle \partial_z\bar\omega\cdot\xi_\beta(\partial_td_0-\Delta d_0),\rho_0^2\rangle
-2\langle \rho_0'\nabla d_0\cdot\nabla\bar\omega\cdot\xi_\beta,\rho_0\rangle\\
&-2\langle \nabla d_0\cdot\nabla\partial_z\bar\omega\cdot\xi_\beta,\rho_0^2\rangle
+\frac{1}{d_0}\langle \partial_z\bar{\om}\cdot\xi_\beta,\rho_0^2\rangle\Big)\notag\\
=&L_{1,\beta}(d_0,\rho_0,\bar\omega,\xi_\beta),
\label{L1beta}
\end{align}
By \eqref{g0xibeta} and \eqref{v2eq}, one has
\begin{align*}
 -\frac{\kappa_\beta d_1}{d_0}\Big( \partial_z^2(\eta\rho_0) -f_B(\rho_0)\eta\rho_0\Big)=-\frac{\kappa_\beta d_1}{d_0}\frac{1}{\rho_0}\partial_z(\rho_0^2\eta')=g_0\cdot\xi_\beta d_1.
\end{align*}
Then there exists $e_{1,\beta}(x,t)$ such that
 \begin{align}
 \sigma_{1,\beta}(z,x,t) &=-\frac{\kappa_\beta d_1}{d_0}\eta\rho_0+\sigma_{1,\beta}^*(d_0,\rho_0,\bar\omega,\xi_\beta)+h_{1,\beta}(x,t)d_0\rho_0\eta+e_{1,\beta}(x,t)\rho_0(z),\label{sigma1beta}
 \end{align}
 where $\sigma_{1,\beta}^*$ solves
 \begin{align}\label{sigma1beta*}
     (\pa_z^2-f_B(\rho_0))&\sigma_{1,\beta}^*
     =L_{1,\beta}\rho_0^{-1}\eta' d_0+\big(\partial_zu_0(\partial_td_0-\Delta d_0)-2\nabla d_0\cdot\nabla\partial_z u_0 -zg_0\big)\cdot\xi_\beta.
 \end{align}
On $\G$, then we have
 \begin{align}
  \sigma_{1,\beta}(z,x,t)=&-\frac{\kappa_\beta d_1}{d_0}\eta(z)\rho_0(z)+\sigma^*_{1,\beta}(z,x,t)+
 e_{1,\beta}(x,t)\rho_0(z).\label{sigma1ongamma}
 \end{align}
 Taking $z\rightarrow \pm\infty$ in \eqref{sigma1beta}, we have
 \begin{align*}
    \sigma_{1,\beta}^+&=-\frac{b\kappa_\beta d_1}{d_0}+\sigma_{1,\beta}^{*+}+bh_{1,\beta}d_0+be_{1,\beta},\notag\\
    \sigma_{1,\beta}^-&=\sigma_{1,\beta}^{*-}+ae_{1,\beta},
\end{align*}
 which yields
\begin{equation}\label{e1betabefore}
\begin{aligned}
 &e_{1,\beta}=\frac{1}{a}(\sigma_{1,\beta}^--\sigma_{1,\beta}^{*-}),\\
 &h_{1,\beta}(x,t)=
 \begin{cases}
 \frac{a\sigma_{1,\beta}^+-b\sigma_{1,\beta}^-+b\sigma_{1,\beta}^{*-}-a\sigma_{1,\beta}^{*+}+ab\kappa_\beta d_1/d_0}{ab d_0},\quad \qquad \qquad\qquad\qquad \qquad d_0\neq 0,\\
 \frac{1}{ab}\Big(a\pa_\nu\sigma_{1,\beta}^+-b\pa_\nu\sigma_{1,\beta}^-+b\pa_\nu\sigma_{1,\beta}^{*-}-a\pa_\nu\sigma_{1,\beta}^{*+}\Big)+\pa_\nu\Big(\frac{ \kappa_\beta d_1}{d_0}\Big),\quad d_0=0.
 \end{cases}
 \end{aligned}
\end{equation}

On $\G$, substituting  \eqref{g0omega}, \eqref{g0xibeta},  \eqref{bot} and  \eqref{d0} into \eqref{d1omega} and \eqref{d1xibeta}, we have
\begin{align}
&\mathcal{D}_1\cdot\bar{\om}=0,\\
 &\mathcal{D}_1\cdot\xi_\beta=-2\rho_0'\pa_\nu\bar\omega\cdot\xi_\beta-2\rho_0\partial_{\nu z}\bar\omega\cdot\xi_\beta+(\pa_z^2\pa_\nu\bar{\omega}\rho_0+2\pa_{\nu z}\bar{\omega}\rho_0')\cdot \xi_\beta(d_1-z).\label{d1betagamma}
\end{align}
 We deduce from the  definition of $\rho_1$ in \eqref{rho1} that
\begin{align}\label{rho1gamma}
\rho_1|_\G=0.
\end{align}
Moreover, it is straightforward to conclude that  $\mathcal{D}_1$, $L_{1,0}$,   $\rho_1$ and $\sigma_{1,\beta}$ are all linear in $d_1$ by their definitions in $\G(\delta)$. While $L_{1,0}^*$, $L_{1,\beta}$ and $\sigma_{1,\beta}^*$ depends only on $\mathcal{V}^0$.
\subsection{Deriving equations for the $O(1)$ system}
Setting $k=2$ in \eqref{rhol} and \eqref{sigmal}, we obtain
\begin{align}
\partial_z^2\rho_2-f_A(\rho_0)\rho_2=&L_{2,0}\eta'd_0+\mathcal{D}_2\cdot\bar\omega,\label{innerrho2}\\
\big(\partial_z^2-f_B(\rho_0)\big)\big(\sigma_{2, \beta}-h_{2,\beta}d_0\rho_0\eta\big)=&L_{2, \beta}\rho_0^{-1}\eta'd_0+\mathcal{D}_2\cdot\xi_{\beta},\label{innersigma2}
\end{align}
with the boundary conditions
\begin{align*}
\rho_2(\pm\infty,x,t)=\rho_2^\pm(x,t),\ \sigma_{2,\beta}(\pm\infty,x,t)=\sigma_{2,\beta}^\pm(x,t),
\end{align*}
where
\begin{align*}
\mathcal{D}_2=&\partial_zu_0(\partial_td_{1}-\Delta d_{1})+\partial_zu_{1}(\partial_td_{0}-\Delta d_{0})-2\nabla d_{1}\cdot\nabla \partial_zu_0-2\nabla d_{0}\cdot\nabla \partial_zu_{1}\nonumber\\
&+f^{(1)}(u_0,u_{1})+d_2g_0+(d_1-z)g_{1}+\partial_{t}u_{0}-\Delta u_{0}+H_0^+\eta_M^++H_0
^-\eta_M^-.
\end{align*}
\begin{lemma}\label{led1}
On $\G$,
    $d_1$ satisfies
    \begin{align}\label{d1evolution}
    \pa_t d_1-\Delta d_1=f_{1,1}d_1+f_{1,2},
\end{align}
where
\begin{align*}
   ef_{1,1}=&-\sum_\al\pa_\nu\kappa_\al\langle \eta'',\rho_0'\sigma_{1,\al}^*-\rho_0\pa_z\sigma_{1,\al}^*\rangle+\f32\sum_{\al=1}^{n-1} |\pa_\nu\kappa_\al|^2\langle (\rho_0)^2,(\eta')^2\rangle-L_{1,0}^*\langle \eta',\rho_0'\rangle,\\
   ef_{1,2}=&-\langle \pa_tu_0-\Delta u_0,\pa_zu_0\rangle-\sum_\al\pa_\nu\kappa_\al\langle z\eta''+2\eta',\rho_0\pa_z\sigma_{1,\al}^* -\rho_0'\sigma_{1,\al}^*\rangle\\
    &-\pa_\nu\langle z\rho_0',(\pa_t d_0-\Delta d_0)\rho_0'-2\rho_0\nabla d_0\cdot\nabla\pa_z\bar{\omega}\cdot\bar{\omega}\rangle-\sum_{\al=1}^{n-1} |\pa_\nu\kappa_\al|^2\langle z\rho_0',z\rho_0(\eta')^2 \rangle
\end{align*}
only depends on $\mathcal{V}^0$.
\end{lemma}
\begin{proof}
    We have put the proof in Appendix \ref{A.1}.
\end{proof}

Once $d_1$ is determined on $\G$, we can extend it to the neighborhood $\G(\delta)$ by the ODE:
\begin{align}\label{d0d1}
\nabla d_0\cdot\nabla d_1=0.
\end{align}
To make \eqref{innersigma2} solvable on $\G$, we need the condition:
\begin{equation}\label{d2betarho0}
\langle \mathcal{D}_2\cdot\xi_\beta,\rho_0\rangle=0.
\end{equation}
After a routine calculation about  \eqref{d2betarho0}, we have the following lemma.
\begin{lemma}\label{lem:boundary}
On $\G$, $\sigma_{1,\beta}^\pm$ satisfies mixed boundary conditions:
\begin{align}
b\pa_\nu \sigma_{1,\beta}^+&-a\pa_\nu\sigma_{1,\beta}^-+b\sum_{\al\neq \beta}\sigma_{1,\al}^{+}\pa_\nu\xi_\al^+\cdot\xi_\beta^+-a\sum_{\al\neq \beta}\sigma_{1,\al}^{-}\pa_\nu\xi_\al^-\cdot\xi_\beta^-=H_1(\sigma_{1,\al}^*,d_1),\label{sigma1betaboundary1}
\end{align}
and
 \begin{equation}\label{e1beta1}
a\sigma_{1, \beta}^+(x, t)-b\sigma_{1, \beta}^-(x, t)=a\sigma_{1,\beta}^{*+}-b\sigma_{1,\beta}^{*-}-ab d_1\pa_\nu\kappa_\beta,
 \end{equation}
 where
 \begin{align*}
H_1(\sigma_{1,\al}^*,d_1)
=&a\pa_\nu\sigma_{1,\beta}^{*-}-b\pa_\nu\sigma_{1, \beta}^{*+}+\sum_{\al\neq \beta}b^2d_1\pa_\nu\kappa_\al\pa_\nu\xi_\al^+\cdot\xi_\beta^+\notag\\
&+2\pa_\nu\langle\sigma_{1, \beta}^*, \rho_0'\rangle+d_1\pa_\nu\Big(\frac{\kappa_\beta}{d_0}\Big)\langle\eta',\rho_0^2 \rangle-2d_1\sum_{\al\neq\beta}\pa_\nu\kappa_\al\langle\eta\pa_\nu\xi_\al\cdot \xi_\beta,\rho_0\rho_0' \rangle\notag\\
&-b\sum_{\al\neq \beta}\sigma_{1, \al}^{*+}\pa_\nu\xi_\al^+\cdot\xi_\beta^++a\sum_{\al\neq \beta}\sigma_{1, \al}^{*-}\pa_\nu\xi_\al^-\cdot\xi_\beta^-+2\sum_{\al\neq \beta}\langle\sigma_{1, \al}^*\pa_\nu\xi_\al\cdot\xi_\beta,\rho_0'\rangle\notag\\
&+L_{1, \beta}\langle d_1-z,\eta'(z)\rangle-2\langle(\nabla d_{1}\cdot\nabla )\partial_zu_0\cdot\xi_\beta,\rho_0\rangle+\langle (\pa_t u_0-\Delta u_0)\cdot\xi_\beta ,\rho_0 \rangle.
 \end{align*}
\end{lemma}
\begin{proof}
We have put the proof in Appendix \ref{A.2}.
\end{proof}

\subsection{Deriving equations for the  $O(\ve^{k-1})$ $(k\geq 1)$ system}
Consider \eqref{rhol} and \eqref{sigmal} for $k+1$. By performing an inner product with $\bar{\omega}$ and $\xi_{\beta}$ respectively, we obtain
\begin{align}
\partial_z^2\rho_{k+1}-f_A(\rho_0)\rho_{k+1}=&L_{{k+1},0}\eta'd_0+\mathcal{D}_{k+1}\cdot\bar
\omega,\label{inneromega}\\
\Big(\partial_z^2-f_B(\rho_0)\Big)\Big(\sigma_{k+1,\beta}-h_{{k+1},\beta}(x,t)d_0\rho_0\eta\Big)=&L_{{k+1}, \beta}\rho_0^{-1}\eta'd_0+\mathcal{D}_{k+1}\cdot\xi_\beta,\label{innersigma}
\end{align}
 with the boundary conditions
 \begin{align*}
\rho_{k+1}(\pm\infty,x,t)=\rho_{k+1}^\pm(x,t),\
\sigma_{{k+1},\beta}(\pm\infty,x,t)=\sigma_{{k+1},\beta}^\pm(x,t).
\end{align*}
To make  \eqref{inneromega} solvable in $\G(\delta)$, by the compatibility conditions \eqref{solvablecondition}, one has
\begin{equation}\label{lk0}
L_{k+1,0}=-\frac{\langle\mathcal{D}_{k+1}\cdot\bar{\omega},\rho_0'\rangle}{d_0\langle\rho_0',\eta'\rangle}=\frac{d_{k+1}\kappa_0^2}{d_0^2}c_1^*+L_{k+1,0}^*,
\end{equation}
where $c_1^*=\langle\rho_0\rho_0',(\eta')^2 \rangle/\langle\rho_0',\eta'\rangle$ and $L_{k+1,0}^*$ depends only on $\mathcal{V}^i$ for $i\leq k$.
Similarly, for \eqref{innersigma}, we have
\begin{equation}\label{lkbeta}
L_{k+1,\beta}=-\frac{1}{d_0}\langle\mathcal{D}_{k+1}\cdot\xi_\beta,\rho_0\rangle=L_{k+1,\beta}(\mathcal{V}^i)\text{ for }i\leq k.
\end{equation}
Using
$$
g_0(z,x,t)\cdot\xi_{\beta}
 =-\frac{\kappa_\beta}{d_0\rho_0}\pa_z(\rho_0^2\eta'\big)\text{ and }f_B(\rho_0)=\frac{\rho_0''}{\rho_0},$$
one has $\int_\BR \rho_0g_0\cdot \xi_\beta \ud z=0$ and
\begin{align*}
 -\frac{\kappa_\beta d_{k+1}}{d_0}\Big( \partial_z^2(\eta\rho_0) -f_B(\rho_0)\eta\rho_0\Big)=-\frac{\kappa_\beta d_{k+1}}{d_0}\frac{1}{\rho_0}\partial_z(\rho_0^2\eta')=d_{k+1} g_0\cdot\xi_\beta .
\end{align*}
 Lemma \ref{lem:s1} yields that there exists $e_{k+1,\beta}(x,t)$ such that
\begin{align}
\rho_{k+1}(z,x,t)&=-\rho_0'(z)\int_{0}^{z}\rho_0'(\varsigma)^{-2}\int_\varsigma^{+\infty}(L_{k+1,0}\eta'd_0+\mathcal{D}_{k+1}\cdot\bar{\omega})\rho_0'(\tau)\ud \tau \ud \varsigma\notag\\
&=-\frac{\kappa_0^2d_{k+1}}{d_0}\rho_0'(z)\int_{0}^{z}\rho_0'(\varsigma)^{-2}\int_\varsigma^{+\infty}
\big(c_1^*\eta'-\rho_0(\eta')^2\big)\rho_0'\ud s \ud \varsigma+\rho_{k+1}^*(z,x,t)
\label{rhok}
\end{align}
and
 \begin{align}
 \sigma_{k+1,\beta}(z,x,t)&=-\frac{\kappa_\beta d_{k+1}}{d_0}\eta\rho_0+\sigma_{k+1,\beta}^*(z,x,t)+h_{k+1,\beta}(x,t)d_0\rho_0\eta+e_{k+1,\beta}(x,t)\rho_0,\label{sigmakbeta}
 \end{align}
 where $\sigma_{k+1,\beta}^*(z,x,t)$ solves
 \begin{align*}
     \pa_z^2\sigma_{k+1,\beta}^*-f_B(\rho_0)\sigma_{k+1,\beta}^*=&d_0L_{k+1,\beta}\rho_0^{-1}\eta'+(\mathcal{D}_{k+1}-d_{k+1}g_0)\cdot\xi_\beta.
 \end{align*}
 We have by \eqref{sigmakbeta} that
 \begin{align*}
  \sigma_{k+1,\beta}(z,x,t)|_{\G}=&-\frac{\kappa_\beta d_{k+1}}{d_0}\eta\rho_0+\sigma^*_{k+1,\beta}(z,x,t)+
 e_{k+1,\beta}(x,t)\rho_0.
 \end{align*}
Taking $z\rightarrow \pm \infty$ in \eqref{sigmakbeta}, we obtain
\begin{align*}
    &\sigma_{{k+1},\beta}^+=-\frac{b\kappa_\beta d_{k+1}}{d_0}+\sigma_{k+1,\beta}^{*+}+bh_{k+1,\beta}d_0+be_{k+1,\beta},\notag\\
    &\sigma_{{k+1},\beta}^-=\sigma_{k+1,\beta}^{*-}+ae_{k+1,\beta},
\end{align*}
which yields
\begin{equation}\label{ek+1betabefore}
\begin{aligned}
    &e_{k+1,\beta}=\frac{1}{a}(\sigma_{{k+1},\beta}^--\sigma_{k+1,\beta}^{*-}),\\
     &h_{{k+1},\beta}(x,t)=
     \begin{cases}
     \frac{a\sigma_{{k+1},\beta}^+-b\sigma_{{k+1},\beta}^-+b\sigma_{{k+1},\beta}^{*-}-a\sigma_{{k+1},\beta}^{*+}+ab\kappa_\beta d_{k+1}/d_0}{ab d_0},\qquad  d_0\neq 0,\\
     \frac{1}{ab}\Big(a\pa_\nu\sigma_{{k+1},\beta}^+-b\pa_\nu\sigma_{{k+1},\beta}^-+b\pa_\nu\sigma_{{k+1},\beta}^{*-}-a\pa_\nu\sigma_{{k+1},\beta}^{*+}\Big)\\
     \ +\pa_\nu\Big(\frac{\kappa_\beta d_{k+1}}{d_0}\Big), \qquad d_0=0.
     \end{cases}
\end{aligned}
\end{equation}

  We conclude that $\mathcal{D}_{k+1},\ L_{k+1,0}$ and $\rho_{k+1}$  depend linearly on $d_{k+1}$ in $\G(\delta)$. While $L_{k+1,0}^*,\ L_{k+1,\beta}$ and $\sigma_{k+1,\beta}^*$ depends only on $\mathcal{V}^i$ for $i\leq k$. Moreover, by  \eqref{rhok}, we have that   $\rho_{k+1}|_{\G}$ depends only on $\mathcal{V}^i$ for $i\leq k$.
\subsection{Deriving the jump conditions from $O(\ve^{k})$$(k\geq 1)$ term}
Recall  \eqref{rhol} and \eqref{sigmal} $(k+2)$, we have
\begin{align*}
\partial_z^2\rho_{k+2}-f_A(\rho_0)\rho_{k+2}=&~d_0L_{{k+2},0}\eta'+\mathcal{D}_{k+2}\cdot\omega,\\
\Big(\partial_z^2-f_B(\rho_0)\Big)\Big(\sigma_{{k+2}, \beta}-h_{{k+2},\beta}(x,t)d_0\rho_0\eta\Big)=&~d_0L_{{k+2}, \beta}\rho_0^{-1}\eta'+\mathcal{D}_{k+2}\cdot\xi_{\beta},
\end{align*}
where
\begin{align*}
\mathcal{D}_{k+2}=&\partial_zu_0(\partial_td_{k+1}-\Delta d_{k+1})+\partial_zu_{k+1}(\partial_td_{0}-\Delta d_{0})-2(\nabla d_{k+1}\cdot\nabla) \partial_zu_0\nonumber\\
&-2(\nabla d_{0}\cdot\nabla) \partial_zu_{k+1}+f^{(k+1)}(u_0, \cdots, u_{k+1})+(d_{k+2}g_0+d_{k+1}g_1)\nonumber\\
&+(d_1-z)g_{k+1}+\mathcal{A}_{k}.
\end{align*}
On $\Gamma$, by the compatibility conditions \eqref{solvablecondition},  we obtain
\begin{align}\label{dk+2rho0'}
0=\langle \mathcal{D}_{k+2}\cdot\bar{\omega},\rho_0'\rangle.
\end{align}
\begin{lemma}\label{lemdk}
On $\G$, $d_{k+1}$ satisfies
\begin{align}\label{dkevolution}
\pa_t d_{k+1}-\Delta d_{k+1}=f_{k+1,1}d_{k+1}+f_{k+1,2},
\end{align}
where
\begin{align*}
 ef_{k+1,1}=&-\langle g_1,\pa_zu_0\rangle-\sum_\al\langle g_0,\sigma_{1,\al}'\xi_\al\rangle-\sum_{\al=1}^{n-1} |\pa_\nu\kappa_\al|^2\big(\langle d_1-z,\eta'\eta''(\rho_0)^2+c_1^*\eta'\rho_0'\rangle\\
 &+\langle z\rho_0',c_1^*\eta'-\rho_0(\eta')^2\rangle-2\langle \rho_0\eta',\rho_0\eta'\rangle\big),\\
 ef_{k+1,2}=&-\pa_\nu\langle z\rho_0',\mathcal{D}_{k+1}^r\cdot \bar{\omega}\rangle-\langle \mathcal{A}_k,\pa_z u_0\rangle-d_1L_{k+1,0}^*\langle\eta', \rho_0'\rangle -\langle f_4^{(k+1)}\cdot \bar{\omega},\rho_0'\rangle\\
    &+\sum_\al \Big(\pa_\nu\kappa_\al\langle 2\eta'-(d_1-z)\eta'',\rho_0'\sigma_{k+1,\al}^*-\rho_0\pa_z \sigma_{k+1,\al}^* \rangle-\langle \mathcal{D}_{k+1}^r\cdot\xi_\al,\sigma_{1,\al}'\rangle\Big).
\end{align*}
only depends on $\mathcal{V}^i$ for $i\leq k.$ Here, $\mathcal{D}_{k+1}^r:=\mathcal{D}_{k+1}-d_{k+1}g_0$ does not depend on $d_{k+1}$.
\end{lemma}
\begin{proof}
    We have put the proof in Appendix \ref{A.3}.
\end{proof}
We conclude that \eqref{dkevolution} is linear in $d_{k+1}$. Together with
\begin{equation}\label{d0dk}
\nabla d_0\cdot\nabla d_{k+1}=-\frac{1}{2}\sum_{i=1}^{k}\nabla d_i\cdot\nabla d_{k+1-i},\ \text{for}\ (x,t)\in\Gamma(\delta),
\end{equation}
we can extend $d_{k+1}$ from $\Gamma$ into $\Gamma(\delta)$.

On $\G$, using compatibility condition \eqref{solvablecondition}, we obtain
\begin{equation}\label{eq:Dk+2}
 \langle\mathcal{D}_{k+2}\cdot\xi_\beta, \rho_0 \rangle_\G=0.
 \end{equation}
By the following lemma, the equation \eqref{eq:Dk+2} gives jump conditions for $\sigma_{k+1,\beta}^\pm$ on $\Gamma$.
\begin{lemma}\label{lem:sigmakbetaboundary}
On $\G$, $\sigma_{k+1, \beta}^\pm$ satisfies mixed boundary conditions:
 \begin{align}\label{sigmakbetaboundary}
b\partial_\nu&\sigma^+_{k+1,\beta}-a\partial_\nu\sigma^-_{k+1,\beta}+b\sum_{\al\neq \beta}\sigma_{k+1, \al}^+\pa_\nu\xi_\al^+\cdot\xi_\beta^+-a\sum_{\al\neq \beta}\sigma_{k+1, \al}^-\pa_\nu\xi_\al^-\cdot\xi_\beta^-\\
&=H_{k+1}(\sigma_{k+1,\al}^{*},d_{k+1}),\notag
\end{align}
and
\begin{equation}
a\sigma_{k+1, \beta}^+(x, t)-b\sigma_{k+1, \beta}^-(x, t)=a\sigma_{k+1,\beta}^{*+}-b\sigma_{k+1,\beta}^{*-}-ab d_{k+1}\pa_\nu\kappa_\beta,
 \label{sigmakbetaboundary2}
 \end{equation}
where
 \begin{align*}
H_{k+1}(\sigma_{k+1,\al}^{*},d_{k+1})=&a\pa_\nu\sigma_{{k+1},\beta}^{*-}-b\pa_\nu\sigma_{{k+1}, \beta}^{*+}+\sum_{\al\neq \beta}b^2d_{k+1}\pa_\nu\kappa_\al\pa_\nu\xi_\al^+\cdot\xi_\beta^+\notag\\
&+\sum_{\al\neq \beta}-b\sigma_{{k+1},\al}^{*+}\pa_\nu\xi_\al^+\cdot\xi_\beta^++a\sigma_{{k+1},\al}^{*-}\pa_\nu\xi_\al^-\cdot\xi_\beta^-+2\pa_\nu\langle\sigma_{{k+1}, \beta}^* ,\rho_0'\rangle\\
&+\pa_\nu\Big(\frac{d_{k+1}\kappa_\beta}{d_0}\Big)\langle \eta',\rho_0^2 \rangle+2\sum_{\al\neq \beta}\langle\sigma_{{k+1}, \al}^*\pa_\nu\xi_\al\cdot\xi_\beta,\rho_0'\rangle\\
&+\langle f_2^{(k+1)}(u_0,u_{k+1})\sigma_{1,\beta}+f_4^{(k+1)}(u_0,\cdots,u_{k})\cdot\xi_\beta+ \mathcal{A}_k\cdot\xi_\beta,\rho_0\rangle\notag\\
&-2d_{k+1}\sum_{\al\neq\beta}\pa_\nu\kappa_\al\langle\eta\rho_0'\pa_\nu\xi_\al\cdot \xi_\beta, \rho_0\rangle -2\langle \pa_z\rho_{k+1}\pa_\nu\bar{\omega}\cdot\xi_\beta, \rho_0\rangle\\
&-2\langle\rho_{k+1}\pa_{z\nu}\bar{\omega}\cdot \xi_\beta, \rho_0\rangle +L_{{k+1}, \beta}\langle d_1-z,\eta'(z)\rangle\\
&-\pa_\nu\kappa_\beta\langle (d_1-z)(2\eta'\rho_{k+1}'+\eta''\rho_{k+1}),\rho_0\rangle\\
&+d_{k+1}\langle g_1\cdot\xi_\beta,\rho_0\rangle-2\langle (\nabla d_{k+1}\cdot\nabla)\pa_zu_0\cdot\xi_\beta,\rho_0\rangle.
 \end{align*}
 \end{lemma}
\begin{proof}
We have put the proof in Appendix \ref{A.4}.
\end{proof}
\section{Solving the systems for outer/inner expansion}\label{section4}
In this section, we first solve the systems for the outer/inner expansions, and
then construct the approximate solutions by the gluing method.
\subsection{Solving $\mathcal{V}^0$}
\subsubsection{Solving $u_0^\pm$ in $U_\pm$}
By \eqref{u0pm}, we obtain $u_0^\pm$ in $U_\pm$.
\subsubsection{Solving $u_0$ in $\G(\delta)$}
We can find $u_0$ in $\G(\delta)$ by \eqref{u0}.
\subsubsection{Solving $d_0$ in $\G(\delta)$}
The evolution equation \eqref{d0} for $d_0$ means that the interface $\G$ evolves by mean curvature flow. In $\G(\delta)$, $d_0$ can be determined by the system:
\begin{equation*}
\left\{
\begin{array}{ll}
\pa_t d_0-\Delta d_0=0, &(x, t)\in\Gamma, \\
d_0=0, &(x, t)\in\Gamma, \\
\nabla d_0\cdot \nabla d_0=1, &(x, t)\in\Gamma(\delta).
\end{array}
\right.
\end{equation*}
\subsubsection{Solving $g_0$ in $\G(\delta)$}
We can obtain $g_0$ using \eqref{g0}.
\subsection{Solving $\mathcal{V}^1$}
\subsubsection{Solving $\rho_1^\pm$ in $U_\pm$ and $\rho_1^*,\ \sigma_{1,\beta}^*,\ L_{1,0}^*,\ L_{1,\beta}$ in $\G(\delta)$}
We can write $\rho_1^\pm$, $\rho_1^*,\ \sigma_{1,\beta}^*,\ L_{1,0}^*$ and $L_{1,\beta}$ by \eqref{rho1+}, \eqref{rho1_} and \eqref{L10}-\eqref{sigma1beta*}.
\subsubsection{Solving $d_1$ in $\G(\delta)$}
For $d_1$, by \eqref{d1evolution}, one has
\begin{equation}\label{d1ev}
\begin{cases}
\pa_t d_1-\Delta d_1=f_{1,1}(u_0,\sigma_{1,\al}^*,L_{1,0}^*)d_1+f_{1,2}(u_0,\sigma_{1,\al}^*,L_{1,0}^*),\ \ \ (x, t)\in\Gamma, \\
d_1(x, 0)=0, \ \ \ \ \ \ \ \ \ \ \ \ \ \ \ \ \ \ \ \ \ \ \ \ \ \ \ \ \ \ \ \ \ \ \ \ \ \ \ \ \ \ \ \ \ \ \ \ \ \ \ \ \ \ \ \ \ \ \ x\in\G_0,
\end{cases}
\end{equation}
It is a linear equation for $d_1$. Moreover, $f_{1,i}$ for $i=1,2$ is known. We can determine $d_1$ on $\G$ directly by \eqref{d1ev}.
Then  $d_1$ in $\G(\delta)$ can be determined by the equation
\begin{equation}\label{d1delta}
\nabla d_0\cdot \nabla d_1=0.
\end{equation}
\subsubsection{Solving $\sigma_{1,\beta}^\pm$ in $ U_\pm$}
About $\sigma_{1,\beta}^\pm$, combining parabolic system \eqref{outersigmak}$(k=2)$ in $ U_\pm$ with  boundary conditions \eqref{sigma1betaboundary1} and \eqref{e1beta1}, we obtain
\begin{equation}\label{sigma1betasytem}
\begin{cases}
\begin{array}{ll}
\partial_t\sigma_{1, \beta}^\pm-\Delta\sigma_{1, \beta}^\pm-\sigma_{1, \beta}^\pm\Delta \xi_\beta^\pm\cdot \xi_\beta^\pm\\
\ \ =-\rho_{1}^\pm\partial_t\omega^\pm\cdot\xi_{\beta}^\pm+\rho_{1}^\pm\Delta\omega^\pm\cdot\xi_{\beta}^\pm+2(\nabla\rho_{1}^\pm\cdot\nabla)\omega^\pm\cdot\xi_{\beta}^\pm\\
\ \ \ \ \ \ +\sum_{\alpha\neq \beta}\sigma_{1, \alpha}^\pm\Delta\xi_\alpha^\pm\cdot\xi_{\beta}^\pm+2(\nabla\sigma_{1, \alpha}^\pm\cdot\nabla)\xi_\alpha^\pm\cdot\xi_\beta^\pm\\
\ \ \ \ \ \ -\sum_{\alpha\neq \beta}\sigma _{1, \alpha}^\pm \partial_t\xi_\alpha^\pm \cdot \xi_{\beta}^\pm-f^{(2)}(u_0^\pm, u_1^\pm, u_{2}^\pm)\cdot\xi_\beta^\pm, \ (x,t)\in U_\pm,\\
\sigma_{1, \beta}^\pm(x, t)|_{t=0}=\sigma_{1, \beta}^\pm(x, 0),\ \ \ \ \ \ \ \ \ \ \ \ \ \ \ \ \ \ \ \ \ \ \ \ \ \ \ \ \ \ \ \ \ x\in\Omega^\pm_0,
\end{array}
\end{cases}
\end{equation}
with mixed boundary conditions:
\begin{align}
b\pa_\nu \sigma_{1,\beta}^+&-a\pa_\nu\sigma_{1,\beta}^-
+b\sum_{\al\neq \beta}\sigma_{1,\al}^{+}\pa_\nu\xi_\al^+\cdot\xi_\beta^+-a\sum_{\al\neq \beta}\sigma_{1,\al}^{-}\pa_\nu\xi_\al^-\cdot\xi_\beta^-=H_1(\sigma_{1,\al}^*,d_1),\label{sigma1betaboundary}
\end{align}
and
 \begin{equation}\label{e1beta}
a\sigma_{1, \beta}^+(x, t)-b\sigma_{1, \beta}^-(x, t)=a\sigma_{1,\beta}^{*+}-b\sigma_{1,\beta}^{*-}-ab d_1\pa_\nu \kappa_\beta ,
 \end{equation}
 The right-hand side terms of \eqref{sigma1betasytem} are all independent of $\sigma_{2,\alpha}^\pm$ ($1\leq \alpha\leq n-1$). Combining \eqref{sigma1betasytem}-\eqref{e1beta} and \eqref{rho2equation},  one can then determine  $\sigma_{1,\beta}^\pm|_{ U_\pm}$.
 \subsubsection{Solving $u_1$ and $g_1$ in $\G(\delta)$}
Once  $d_1$ is found, we can obtain $L_{1, 0}$ \eqref{L10}, which further gives $\rho_1$ \eqref{rho1}. Using $\sigma_{1,\beta}^\pm$, \eqref{sigmapm} and \eqref{e1betabefore}, we can obtain $e_{1,\beta}$ and $h_{1,\beta}$. At last, one can get $\sigma_{1,\beta}$ \eqref{sigma1beta}  and $g_1$ \eqref{gk}$(k=1)$.
 \subsection{Solving $\mathcal{V}^{k+1}$}
To close this section, we need to find $(k+1)$th-order functions
$$\mathcal{V}^{k+1}=(u_{k+1}^\pm,u_{k+1},d_{k+1},g_{k+1})$$
applying mathematical induction.
  Assume $\mathcal{V}^i$ has been obtained for all $i \leq k$. We then solve for $\mathcal{V}^{k+1}$ via induction.
\subsubsection{Solving $\rho_{k+1}^{\pm}$ in $U_\pm$ and $\rho_{k+1}^*,\ \sigma_{k+1,\beta}^*,\ L_{k+1,0}^*,\ L_{k+1,\beta}$ in $\G(\delta)$}
We can determine $\rho_{k+1}^\pm$, $\rho_{k+1}^*,\ \sigma_{k+1,\beta}^*,\ L_{k+1,0}^*$ and $L_{k+1,\beta}$ by \eqref{outersolve} and \eqref{lk0}-\eqref{sigmakbeta}.
\subsubsection{Solving $d_{k+1}$ in $\G(\delta)$}
From the evolution equation \eqref{dkevolution}, we have
\begin{align}\label{dk+1ev}
\begin{cases}
\pa_t d_{k+1}-\Delta d_{k+1}=f_{k+1,1}d_{k+1}+f_{k+1,2},\ \ \ \quad (x,t)\in\G,\\
d_{{k+1}}(x, t)|_{t=0}=0, \ \ \ \ \ \ \ \ \ \ \ \ \ \ \ \ \ \ \quad\qquad \qquad x\in\Gamma_0.
\end{cases}
\end{align}
Using the equation
\begin{equation}\label{dk+1delta}
\nabla d_0\cdot\nabla d_{k+1}=-\frac{1}{2}\sum_{i=1}^{k}\nabla d_i\cdot\nabla d_{k+1-i}, \ (x, t)\in\Gamma(\delta),
\end{equation}
one can   find $d_{k+1}$ in $\G(\delta)$.
\subsubsection{Solving $\sigma_{k+1,\beta}^\pm$ in $ U_\pm$}
Combining \eqref{outersigmak} (with $k$ replaced by $k+2$), the mixed boundary conditions \eqref{sigmakbetaboundary} and \eqref{sigmakbetaboundary2}, we obtain  that
 \begin{align}\label{ksystem}
\begin{cases}
\begin{array}{ll}
\partial_t\sigma_{{k+1}, \beta}^\pm-\Delta\sigma_{{k+1}, \beta}^\pm-\sigma_{{k+1}, \beta}^\pm\Delta \xi_\beta^\pm\cdot \xi_\beta^\pm \\
\ \ =-\rho_{{k+1}}^\pm\partial_t\omega^\pm\cdot\xi_{\beta}^\pm+\rho_{{k+1}}^\pm\Delta\omega^\pm\cdot\xi_{\beta}^\pm+2(\nabla\rho_{{k+1}}^\pm\cdot\nabla)\omega^\pm\cdot\xi_{\beta}^\pm  \\
\ \ \ \ \ +\sum_{\alpha\neq \beta}\sigma_{{k+1}, \alpha}^\pm\Delta\xi_\alpha^\pm\cdot\xi_{\beta}^\pm+2(\nabla\sigma_{{k+1}, \alpha}^\pm\cdot\nabla)\xi_\alpha^\pm\cdot\xi_\beta^\pm  \\
\ \ \ \ \ -\sum_{\alpha\neq \beta}\sigma _{{k+1}, \alpha}^\pm \partial_t\xi_\alpha^\pm \cdot \xi_{\beta}^\pm-f^{(k+2)}(u_0^\pm, \cdots, u_{k+2}^\pm)\cdot\xi_\beta^\pm, &(x,t)\in U_\pm, \\
\sigma_{{k+1}, \beta}^\pm(x, t)|_{t=0}=\sigma_{{k+1}, \beta}^\pm(x, 0), &x\in\Omega_0^\pm,
\end{array}
\end{cases}
\end{align}
with mixed boundary conditions
\begin{align}\label{sigmakbetaboundary1}
b\partial_\nu&\sigma^+_{k+1,\beta}-a\partial_\nu\sigma^-_{k+1,\beta}+b\sum_{\al\neq \beta}\sigma_{k+1, \al}^+\pa_\nu\xi_\al^+\cdot\xi_\beta^+-a\sum_{\al\neq \beta}\sigma_{k+1, \al}^-\pa_\nu\xi_\al^-\cdot\xi_\beta^-\\
&=H_{k+1}(\sigma_{k+1,\al}^{*},d_{k+1}),\notag
\end{align}
and
 \begin{equation}
a\sigma_{k+1, \beta}^+(x, t)-b\sigma_{k+1, \beta}^-(x, t)=a\sigma_{k+1,\beta}^{*+}-b\sigma_{k+1,\beta}^{*-}-ab d_{k+1}\pa_\nu\kappa_\beta.
 \label{sigmakbetaboundary21}
 \end{equation}
Then we can determine  $\sigma_{k+1,\beta}^\pm$ in $ U_\pm$. Indeed, by \eqref{fk-1} and \eqref{u0u1cdot}, we obtain
\begin{align*}
f^{(k+2)}(u_0^\pm,\cdots,u_{k+2}^\pm)\cdot \xi_\beta^\pm=f_2^{(k+2)}\sigma_{1,\beta}^\pm+f_4^{(k+2)}\cdot \xi_\beta^\pm.
\end{align*}
We conclude that $f^{(k+2)}(u_0^\pm,\cdots,u_{k+2}^\pm)\cdot \xi_\beta^\pm$ depends only on
$u_i^\pm(0\leq i\leq k+1)$ and $\rho_{k+2}^\pm$. Moreover, we can determine $\rho_{k+2}^\pm$ by \eqref{outersolve}.
\subsubsection{Solving $u_{k+1}$ and  $g_{k+1}$ in $\G(\delta)$} 
\begin{figure}\small
\begin{tikzpicture}[man/.style={rectangle,draw,fill=blue!20},
  woman/.style={rectangle,draw,fill=red!20,rounded corners=.8ex},
 VIP/.style={rectangle,very thick, draw}, fill=white]
\draw [dotted](-3,0.25) to (-3,-1*15.3);

\node (Out) at (0,0){Inner expansion};
\node (Inn) at (-6,0){Outer expansion};

\node[man] (A0) at (-0.4,-1) {$u_0,\ g_0,\  d_0$};
\node[man] (A0-out) at (-7,-1) {$u_0^\pm|_{ U_\pm}$};

\node[man] (F1) at (-0.4,-1*2)  {$\rho_1^*,\ \sigma_{1,\beta}^*,\ L_{1,0}^*,\ L_{1,\beta}$};
\node[man] (d11) at (-0.4,-1*3)  {$d_1|_\Gamma$};
\node[VIP] (qgamma) at (-0.1,-1*4) {Boundary conditions for $\sigma_{1,\beta}^\pm|_\G$};
\node[woman] at (3.6,-4){\eqref{sigma1betaboundary}, \eqref{e1beta}};
\node[man] (F12delta) at (2,-1*3) {$d_1|_{\G(\delta)}$};
\node[man] (inneru1) at (1,-1*5.5)  {$L_{1,0},\ L_{1,\beta},\ \rho_1,$};
\draw [->](2.5,-1*3.3) |- (inneru1);
\node[woman] at (-1.6,-1*5.5){\eqref{L10}-\eqref{L1beta}};
\node[man] (P1i) at (-0.4,-1*6-0.6)   {$e_{1,\beta},$$h_{1,\beta},$$u_1,$$g_1$};
\node[woman] at (-5.5,-1){\eqref{u0pm}};
\node[woman] at (3.2,-1*6-0.6)  { \eqref{rho1},\ \eqref{sigma1beta},\ \eqref{e1betabefore}, \eqref{gk}};
\node[woman] at (3,-1)  {\eqref{u0},\ \eqref{g0},\ \eqref{d0delta}};

\node[woman] at (3,-1*2) {\eqref{L10}-\eqref{sigma1beta}};
\node[woman] at (-1.5,-1*3)  {\eqref{d1ev}};
\node[woman] at (3.5,-1*3)  {\eqref{d1delta}};

\node[man] (M1) at (-7,-1*2)  {$\rho_1^\pm|_{ U_\pm}=0$};
\node[man] (V1) at (-7,-1*6-0.6)  {$\sigma_{1,\beta}^\pm|_{ U_\pm}$};
\node[woman] at (-4.8,-1*2){\eqref{rho1+}, \eqref{rho1_}};
\node[VIP] (Eq1) at (-7,-1*4)  {Equations for $\sigma_{1,\beta}^\pm|_{ U_\pm}$: \eqref{outersigmak} $(k=2)$};
\node (plus1) at (-3,-1*4)  {+};

\node (plus2) at (-3,-1*9-2) {+};
\node (stepk) at (-0.4,-1*7-1)  {\underline{ $\{u_i, ~d_i, ~g_i| i\leq k \}$ are solved}};

\node (stepk) at (-7,-1*7-1){\underline{$\{u_i^\pm|i\le k \}$ are solved}};

\node[man] (F2) at (-0.4,-1*7-2){$\rho_{k+1}^*,\ \sigma_{k+1,\beta}^*,\ L_{k+1,0}^*,\ L_{k+1,\beta}$};
\node[VIP] (qgammak) at (-0.1,-1*9-2)  {Boundary conditions for $\sigma_{k+1,\beta}^\pm|_\G$};
\node[woman] at (3.9,-11){\eqref{sigmakbetaboundary1}, \eqref{sigmakbetaboundary21}};
\node[woman] at(-1.6,-1*8-2){\eqref{dk+1ev}};
\node[man] (d22) at (2,-1*8-2){$d_{k+1}|_{\Gamma(\delta)}$};
\node[man] (F22) at (1.4,-1*10-2){$L_{k+1,0},\ L_{k+1,\beta},$ $\rho_{k+1}$};
\draw [->](2.7,-10.3)to(2.7,-11.7);
\node[woman] at (-1.6,-1*10-2){\eqref{lk0}-\eqref{rhok}};
\node[man] (P2i) at (-0.5,-1*11-2.6)  {$e_{k+1,\beta},h_{k+1,\beta},u_{k+1},g_{k+1}$};
\node[woman] at (3.7,-1*11-2.6)  {\eqref{rhok},\eqref{sigmakbeta},\eqref{ek+1betabefore},\eqref{gk}};

\node[man] (k+1boundary)at (-0.4,-1*8-2)  {$d_{k+1}|_\G$};
\node[woman] at (3.5,-1*8-2)  {\eqref{dk+1delta}};
\node[woman] at (3.4,-1*7-2)  {\eqref{lk0}-\eqref{sigmakbeta}};

\node[man] (Mk) at (-7,-1*7-2)  {$\rho_{k+1}^\pm|_{ U_\pm}$};
\node[man] (Vk) at (-7,-1*11-2.6)  {$\sigma_{k+1,\beta}^\pm|_{ U_\pm}$};
\node[VIP] (Eqk) at (-7,-1*8-3)  {Equations for $\sigma_{k+1,\beta}^\pm|_{ U_\pm}$: \eqref{outersigmak}($k+2$)};
\node[woman] at (-5.5,-1*7-2){\eqref{outersolve}};

\node (stepk) at (0,-1*11-3.5){\underline{$u_{k+1},\ d_{k+1},\ g_{k+1}$ are solved}};
\node (stepk1) at (-7,-1*11-3.5){\underline{$u_{k+1}^\pm$ are solved}};

\draw [->](A0) to  (F1);
\draw [->](F1) to (d11);
\draw [->](d11) to (qgamma);
\draw [->](d11) to (F12delta);
\draw [->](V1) to (P1i);
\draw [->](inneru1) to (P1i);
\draw [->](A0-out) to  (M1);
\draw [->](F22) to  (P2i);
\draw [->](M1) to (Eq1);
\draw [->](V1) to (P1i);
\draw [->](plus1) to [in=0,out=-90,looseness=0.75] (V1.north east);
\draw [->](F2) to (k+1boundary);
\draw [->](k+1boundary) to (qgammak);

\draw [->](k+1boundary) to (d22);

\draw [->] (-7,-1*8-0.2) to (Mk);
\draw [->]  (-0.4,-1*8.2) to (F2);
\draw [->](Mk) to (Eqk);
\draw [->](Vk) to (P2i);
\draw [->](plus2) to [in=0,out=-90,looseness=0.75] (Vk.north east);
\draw [very thick, dotted](0,-1*7.2) to (0, -7.82);
\draw [very thick, dotted](-7,-1*7.2) to (-7, -7.82);

\draw [very thick, dotted](0, -1*14.7) to (0, -1*15.3);
\draw [very thick, dotted](-7, -1*14.7) to (-7, -1*15.3);
\end{tikzpicture}
\caption{The whole procedure to solve the outer and inner expansion systems}
\label{fig:expansion}
\end{figure}
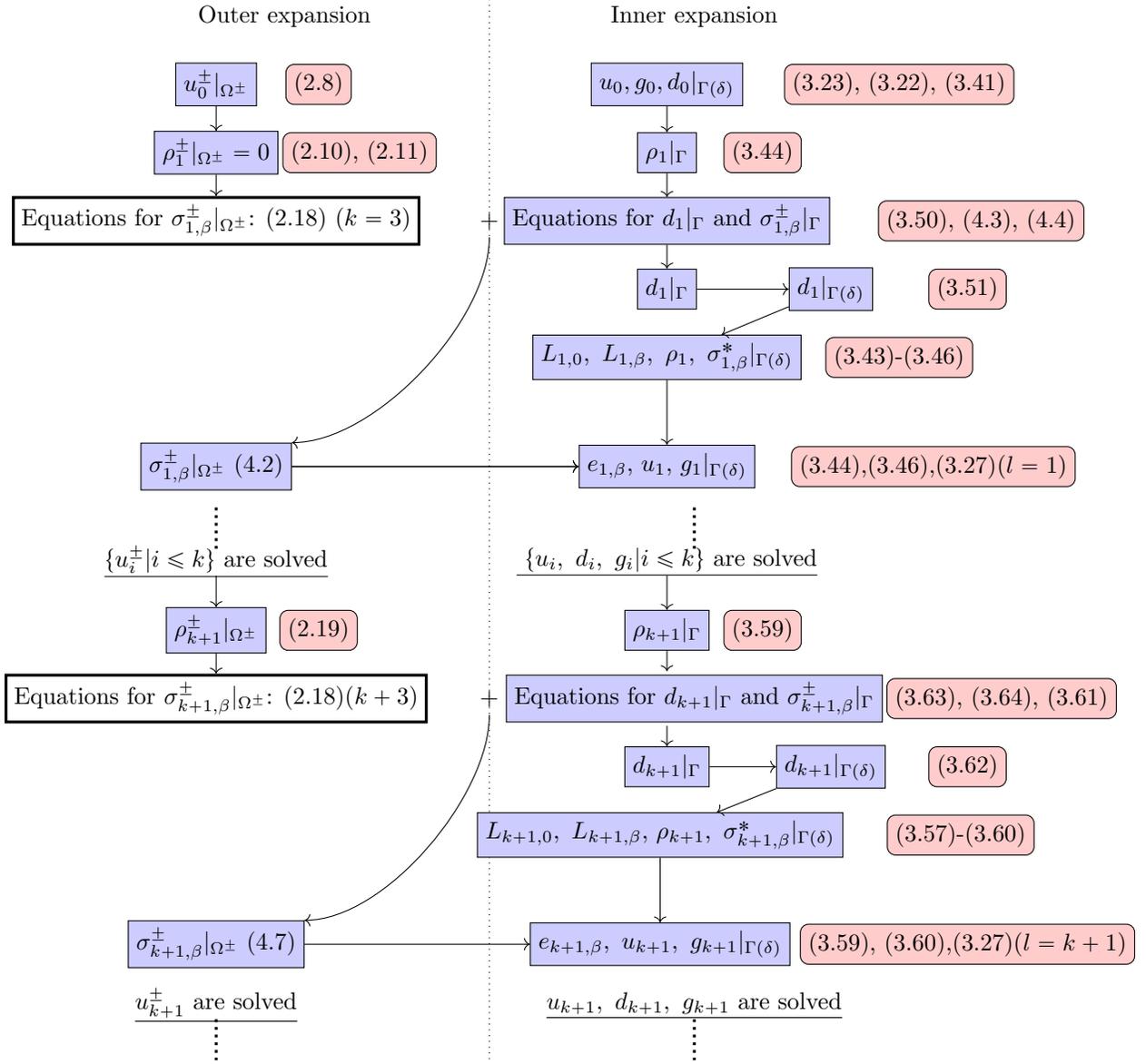

From \eqref{lk0}-\eqref{sigmakbeta}, one gets $L_{k+1, 0}$ and  $\rho_{k+1}$. Using $\sigma_{k+1,\beta}^\pm$, \eqref{sigmak+1pm} and \eqref{ek+1betabefore}, we can obtain $e_{k+1,\beta}$ and $h_{k+1,\beta}$. Then we can write $\sigma_{k+1,\beta}$ \eqref{sigmakbeta}. Finally,  $g_{k+1}$ can be constructed by \eqref{gk}$(l=k+1)$.
 It is not surprising that the step can be promoted. In summary, we can see Figure \ref{fig:expansion}.

\subsection{Constructing  approximate solutions}\label{section4.4}
In this section, we glue together the inner expansion and the outer expansion to obtain the approximate solutions in the whole region $\Omega$. Firstly, we glue together $u_i^+$ and $u_i^-$ by
$$u^K_{o}=\sum_{i=0}^{K}\ve^{i}(u_i^+\chi_{U_+}+u_i^-\chi_{U_-}).$$
where
\begin{equation*}
\chi_{U_\pm}=
\left\{
\begin{array}{ll}
1,&(x,t)\in U_\pm,\\
0,&(x,t)\in \text{other region}.
\end{array}
\right.
\end{equation*}
Then for $(x,t)\in U_\pm$, it holds that
\begin{align*}
\partial_t&u^K_{o}-\Delta u^K_{o}+\frac{1}{\ve^2}f(u^K_{o})\nonumber\\
&=\sum_{k=0}^{K-2}\ve^{k}\big(\partial_tu_k^\pm-\Delta u_k^\pm+D f(u_0^\pm)u_{k+2}^\pm+f^{(k+1)}(u_0^\pm,u_1^\pm,\cdots,u_{k+1}^\pm)\big)+O(\ve^{K-1})\nonumber\\
&=O(\ve^{K-1}).
\end{align*}

As for the inner expansion, we first note
\begin{align*}
&d^K(x,t)=\sum_{k=0}^{K}\ve^kd_k(x,t),\ u_{in}^K(x,t)=\sum_{k=0}^{K}\ve^ku_k\Big(\frac{d^K}{\ve},x,t\Big).
\end{align*}
By \eqref{nablad0nabladk}, we have
\begin{align}\label{nabladK2}
|\nabla d^K|^2=1+\sum_{\substack{1\leq i,j\leq K,\\i+j\geq K+1}}\ve^{i+j}\nabla d^j\cdot\nabla d^i=1+O(\ve^{K+1}).
\end{align}
A routine calculation gives us
\begin{align*}
\partial_tu_{in}^K&-\Delta u_{in}^K+\frac{1}{\ve^2}f(u_{in}^K)\nonumber\\
=&\partial_tu_{in}^K-\Delta u_{in}^K+\frac{1}{\ve^2}f(u_{in}^K)+\Big(\sum_{k=0}^{K}\ve^kd_k-\ve z\Big)\sum_{k=0}^{K}\ve^{k-2}g_k+\sum_{k=0}^{K-2}\ve^k(H_k^+\eta_M^++H_k^-\eta_M^-)\Big|_{z=\frac{d^K}{\ve}}\nonumber\\
=&O(\ve^{K-1}).
\end{align*}

Define
\begin{align}\label{uK}
u^K=u^K_{in}+\Big(1-\xi(\frac{d_0}{\delta})\Big)(u^K_o-u_{in}^K),
\end{align}
where $\xi$ is a smooth cut-off function with
\begin{align}\label{eta-function}
 \text{$\xi(y)=1$ for $|y|\le 1/2$ and $\xi(y)=0$ for $|y|\ge 1$}.
\end{align}
Due to the matching condition \eqref{matching-condition}, we show that
\begin{align*}
\partial_tu^K-\Delta u^K+\ve^{-2}f(u^K)=\mathcal{R}^K,
\end{align*}
where $\mathcal{R}^K=O(\ve^{K-1})$. We have confirmed Theorem \ref{thK-1}.

\section{Spectral lower bound estimate for the linearized operator}\label{spectralestimate}
This section is devoted to proving Theorem \ref{th:uA}, i.e., inequality \eqref{th2ineq} for $u\in H^1(\Omega)$. Obviously, it suffices to consider $\ve$ small enough.
\subsection{Reduction to an inequality on $\Gamma_t^K(\delta/4)$}
In this section, we want to reduce the
inequality \eqref{th2ineq} to $\Gamma_t^K(\delta/4)$, where
\begin{align}\label{def:Gamma-t-K}
\Gamma_t^K(\delta):=\{x:|d^K(x,t)|<\delta\}.
\end{align}
We also recall
\begin{align*}
\Gamma_t(\delta):=\{x:|d_0(x,t)|<\delta\}.
\end{align*}
For each $t\in[0,T]$, we introduce a diffeomorphism from $\Gamma_t^K(\delta/4)$ to $\Gamma_t^K\times [-\delta/4,\delta/4]$, defined as follows
\begin{align}\label{transformation}
(\sigma,r)\in \Gamma_t^K\times(-\delta /4,\delta /4)\to x(\sigma,r)\in\Gamma_t^K(\delta/4),
\end{align}
with
\begin{align}\label{parx}
x(\sigma,0)=\sigma\in\Gamma_t^K\text{ and }\pa_r x(\sigma,r)=\frac{\nabla d^K}{|\nabla d^K|^2}\circ(x(\sigma,r)).
\end{align}
Then we have
$$\frac{d}{dr}\Big(d^K(x(\sigma,r))-r\Big)=0.$$
It follows from $d^K(x(\sigma,0))=0$ that $$d^K(x(\sigma,r),t)=r.$$
Let $$J(\sigma,r)=\det\Big(\frac{\pa x(\sigma,r)}{\pa(\sigma,r)}\Big)$$ be the Jacobian of the transformation. Then one has
$$\ud x=J\ud \sigma \ud r,\ J|_{r=0}=1+O(\ve^{K+1}),\ \text{ and } J(\sigma,r)=1+O(r)+O(\ve^{K+1}).$$

We introduce a group of standard orthogonal  basis in $\BR^n$.
Set $$E_0=\frac{u^K}{|u^K|},$$ where $u^K$ is constructed in \eqref{uK}. We select $E_1,\cdots,E_{n-1}\in \mathbb{S}^{n-1}$ with
\begin{align}\label{eiej}
E_i\cdot E_j=\delta_i^j\ (0\leq i,j\leq n-1),\quad \pa_rE_j=-(\pa_rE_0\cdot E_j)E_0\text{ for }1\leq j\leq n-1.
\end{align}
Indeed, on each local coordinate patch of $\Gamma_t^K$, choose a smooth
orthonormal frame $\{F_1,\ldots,F_{n-1}\}$ of
$E_0(\sigma,0,t)^\perp$. We extend this frame in the $r$-direction
by solving
\[
\partial_rE_j
=-(\partial_rE_0\cdot E_j)E_0,
\qquad
E_j(\sigma,0,t)=F_j(\sigma,t),
\quad 1\le j\le n-1.
\]
Standard linear ODE theory gives a unique smooth extension.
Differentiating the inner products shows that these equations
preserve
\[
E_0\cdot E_j=0,
\qquad
E_i\cdot E_j=\delta_{ij},
\quad 1\le i,j\le n-1.
\]
Thus $\{E_0,\ldots,E_{n-1}\}$ is a smooth orthonormal frame.

Then we decompose $u$ as
\begin{align}\label{decomposition}
u=\sum_{i=0}^{n-1}a_iE_i,
\end{align}
which gives that
\begin{align}\label{g'g''decom}
\frac{1}{\ve^2}&\big(G'(|u^K|^2)|u|^2+2G''(|u^K|^2)(u^K\cdot u)^2\big)\notag\\
&=\frac{1}{\ve^2}\big(\sum_{i=0}^{n-1}a_i^2G'(|u^K|^2)+2a_0^2G''(|u^K|^2)|u^K|^2\big)\notag\\
&=\frac{1}{\ve^2}a_0^2\big(G'(|u^K|^2)+2G''(|u^K|^2)|u^K|^2\big)+\frac{1}{\ve^2}\sum_{i=1}^{n-1}a_i^2G'(|u^K|^2),
\end{align}
where
\begin{align*}
    G'(s)=&(s-a^2)(s-b^2)(2s-a^2-b^2), \\
G''(s)=&(2s-a^2-b^2)^2+2(s-a^2)(s-b^2).
\end{align*}
According to the definition of the approximate solution \eqref{uK}, we have
\begin{align}\label{uK2}
u^K=\hat{u}_0+\ve \hat{u}_1+O(\ve^2),
\end{align}
where
\begin{align}
\hat{u}_0&=\hat{\rho}_0\bar{\omega}+\sum_{i=1}^{n-1}\hat{\sigma}_{0,i}\xi_i,\ \hat{u}_1=\hat{\rho}_1\bar{\omega}+\sum_{i=1}^{n-1}\hat{\sigma}_{1,i}\xi_i,\notag\\
\hat{\rho}_0&=\xi(\frac{d_0}{\delta})\rho_0+\Big(1-\xi(\frac{d_0}{\delta})\Big)b\chi_{U_+}
+\Big(1-\xi(\frac{d_0}{\delta})\Big)a\chi_{U_-},\notag\quad \hat{\sigma}_{0,i}=0,\notag\\
\hat{\rho}_1&=\xi(\frac{d_0}{\delta})\rho_1,\quad
\hat{\sigma}_{1,i}=\xi(\frac{d_0}{\delta})\sigma_{1,i}+\Big(1-\xi(\frac{d_0}{\delta})\Big)\sigma_{1,i}^+\chi_{U_+}+\Big(1-\xi(\frac{d_0}{\delta})\Big)\sigma_{1,i}^-\chi_{U_-}.\label{hat}
\end{align}
By \eqref{uK2}, a direct  computation gives
\begin{align}\label{uA}
|u^K|^2&=\hat{\rho}_0^2+2\ve \hat{\rho}_0\hat{\rho}_1+O(\ve^2), \nonumber\\
G'(|u^K|^2)+2G''(|u^K|^2)|u^K|^2&=G'(|\hat{u}_0|^2)+2G''(|\hat{u}_0|^2)|\hat{u}_0|^2+\ve \hat{\rho}_1f_C(\hat{\rho}_0)+O(\ve^2), \notag\\
G'(|u^K|^2)&=G'(|\hat{u}_0|^2)+\ve \hat{\rho}_1 f_D(\hat{\rho}_0)+O(\ve^2),
\end{align}
where
\begin{align*}
f_C(s)=&6s\big((2s^2-a^2-b^2)^2+4(2s^2-a^2-b^2)s^2+2(s^2-a^2)(s^2-b^2)\big) =f_A'(s), \\
f_D(s)=&2s\big((2s^2-a^2-b^2)^2+2(s^2-a^2)(s^2-b^2)\big)=f_B'(s).
\end{align*}
Together with \eqref{fA'} and \eqref{fB'}, it is easy to find that
\begin{align}\label{g'g''}
G'(|\hat{u}_0|^2)+2G''(|\hat{u}_0|^2)|\hat{u}_0|^2=f_A(\hat{\rho}_0)\ \ \text{and}\ \ G'(|\hat{u}_0|^2)=f_B(\hat{\rho}_0).
\end{align}
\begin{lemma}
For $t\in[0,T]$, $x\in \Omega \backslash \G_t^K(\delta/4)$, we have
 \begin{align}
 \frac{1}{\ve^2}&\big(G'(|u^K|^2)|u|^2+2G''(|u^K|^2)(u^K\cdot u)^2\big)\geq -C|u|^2,
 \end{align}
 where $C$ is independent of $\ve$ and $t\in [0, T]$.
\begin{proof}
Combining  \eqref{g'g''decom}, \eqref{uA} and \eqref{g'g''}, we directly get
\begin{align}
\frac{1}{\ve^2}&\Big(G'(|u^K|^2)|u|^2+2G''(|u^K|^2)(u^K\cdot u)^2\Big)\notag\\
&=\frac{1}{\ve^2}a_0^2\big(f_A(\hat{\rho}_0)+\ve\hat{\rho}_1f_C(\hat{\rho}_0)\big)+\frac{1}{\ve^2}\sum_{i=1}^{n-1}a_i^2\big(f_B(\hat{\rho}_0)+\ve\hat{\rho}_1f_D(\hat{\rho}_0)\big)+O(|u|^2)\notag\\
&=\frac{1}{\ve^2}\Big(a_0^2f_A(\hat{\rho}_0)+\sum_{i=1}^{n-1}a_i^2f_B(\hat{\rho}_0)\Big)+\frac{1}{\ve}\Big(a_0^2\hat{\rho}_1f_C(\hat{\rho}_0)+\sum_{i=1}^{n-1}a_i^2\hat{\rho}_1f_D(\hat{\rho}_0)\Big)+O(|u|^2)\notag\\
&:=I_1+I_2+O(|u|^2).\label{g'g''decom2}
\end{align}
By \eqref{hat}, $f_A(\rho_0(+\infty))=2(a^2-b^2)^2b^2$ and $f_A(\rho_0(-\infty))=2(a^2-b^2)^2a^2$, we obtain that there exists a positive number $C_0$ such that
\begin{equation*}
    f_A(\hat{\rho}_0)\geq (a^2-b^2)^2a^2>0,\ \text{for }|z|\geq C_0.
\end{equation*}
By \eqref{hat}, $f_B(\rho_0(\pm\infty))=0$, we have
\begin{align*}
f_B(\hat{\rho}_0)\geq -Ce^{-\frac{\al_0\delta}{4\ve}}, \ |z|\geq \frac{\delta}{4\ve}.
\end{align*}
Fix $t$. For small $\ve$, one has
\begin{align}\label{fA}
f_A(\hat{\rho}_0)&\geq (a^2-b^2)^2a^2>0,\ \ \ (x,t)\in\{|d^K(x,t)|>\delta \setminus 4\},\\
f_B(\hat{\rho}_0)&\geq
-Ce^{-\frac{\al_0\delta}{4\ve}}, \ \ \ \ \ \ \ \ \ \ \ \ \ (x,t)\in\{|d^K(x,t)|>\delta \setminus 4\}.\label{fB}
\end{align}
Combining \eqref{fA} and \eqref{fB}, for $(x,t)\in \{|d^K(x,t)|>\delta \setminus 4\},$ we have $$I_1\geq -\frac{C}{\ve^2}e^{-\frac{\al_0\delta}{4\ve}}|u|^2\ge -C|u|^2.$$
By the definition of $\hat{\rho}_1$, we can deduce that
\begin{align*}
|\hat{\rho}_1|\leq Ce^{-\frac{\al_0\delta}{4\ve}},\ (x,t)\in \{|d^K(x,t)|>\delta \setminus 4\}.
\end{align*}
Then for $(x,t)\in \{|d^K(x,t)|>\delta \setminus 4\}$, one has
$$I_2\geq -\frac{C}{\ve}e^{-\frac{\al_0\delta}{4\ve}}|u|^2\geq -C|u|^2.$$
We complete the proof.
\end{proof}
\end{lemma}
 Then we only need to consider the inequality in $\G_t^K(\delta/4)$ just as
\begin{align}\label{mainee}
\int_{\Gamma_t^K(\delta/4)}|\nabla u|^2\ud x+\frac{1}{\ve^2}\int_{\Gamma_t^K(\delta/4)}G'(|u^K|^2)|u|^2+2G''(|u^K|^2)(u^K\cdot u)^2\ud x\geq -C\int_{\Gamma_t^K(\delta/4)}|u|^2\ud x.
\end{align}
\begin{Remark}\label{remarkinner}
 We choose $\ve$ sufficiently small  such that
\begin{equation}
\Gamma_t(\delta/8)\subset \Gamma_t^K(\delta/4)\subset \Gamma_t(\delta/2),\ \ \text{for}\ t\in [0,T].\label{subsect:Gamma-delta}
\end{equation}
Then we have $\hat{u}_i=u_i$ for $(x,t)\in \G_t^K(\delta/4)$. We can continue using \eqref{uA} and \eqref{g'g''decom2} with $u_0$ and $u_1$.
\end{Remark}
\subsection{Reduction to inequalities for scalar functions on an interval}
Recall \eqref{transformation} and \eqref{parx}.
Together with \eqref{nabladK2}, we have
\begin{align}\label{nablau2}
|\nabla u|^2\geq\Big|\frac{\nabla d^K}{|\nabla d^K|}\cdot\nabla u\Big|^2=\Big(\frac{\nabla d^K\cdot\nabla u}{|\nabla d^K|^2}\Big)^2|\nabla d^K|^2\geq |\pa_r u|^2-C\ve^{K+1}|\nabla u|^2.
\end{align}
By \eqref{eiej},  we have $\pa_rE_0\cdot \pa_rE_i=0$ for $1\leq i\leq n-1$.   Moreover, by the fact $\partial_r(E_i\cdot E_j)=0$, it is easy to see that
\begin{align}
|\partial_r u|^2&=|\sum_{i=0}^{n-1}(\partial_ra_iE_i+a_i\partial_rE_i)|^2\notag\\
&\geq \sum_{i=0}^{n-1}(\partial_ra_i)^2+2\sum_{i=1}^{n-1}(\partial_ra_0a_iE_0\cdot \partial_rE_i+\pa_ra_ia_0E_i\cdot\pa_rE_0)\notag\\
&=\sum_{i=0}^{n-1}(\partial_ra_i)^2+2\sum_{i=1}^{n-1}(\partial_ra_0a_i-a_0\partial_ra_i)E_0\cdot\partial_rE_i.\label{paru2}
\end{align}
For convenience, note
$$\theta_{1,\ve}=\theta_1(\frac{r}{\ve})=\rho_0'(\frac{r}{\ve})\text{ and }  \theta_{2,\ve}=\theta_2(\frac{r}{\ve})=\rho_0(\frac{r}{\ve}).$$
As a result, by \eqref{g'g''decom2}, \eqref{nablau2}, \eqref{paru2} and Remark \ref{remarkinner}, to prove \eqref{mainee}, we only need to prove
\begin{align}
\int_{\G_t^K}&\int_{-\delta/4}^{\delta/4}\Big((\partial_ra_0)^2+\frac{1}{\ve^2}f_A(\theta_{2, \ve})a_0^2+(\partial_ra_i)^2+\frac{1}{\ve^2}f_B(\theta_{2, \ve})a_i^2\Big)J\ud r\ud \sigma\notag\\
&+\frac{1}{\ve}\int_{\G_t^K}\int_{-\delta/4}^{\delta/4} \big(\rho_1f_C(\theta_{2, \ve})a_0^2+\rho_1 f_D(\theta_{2, \ve})a_i^2\big)J\ud r\ud \sigma \nonumber\\
&+2\sum_{i=1}^{n-1}\int_{\G_t^K}\int_{-\delta/4}^{\delta/4}(\partial_ra_0a_i-\partial_ra_ia_0)E_0\cdot \partial_rE_iJ\ud r\ud \sigma\geq -C\int_{\G_t^K}\int_{-\delta/4}^{\delta/4}(a_0^2+a_i^2)J\ud r\ud \sigma.
\end{align}
Without loss of generality, we assume $\delta/4=1$  and let $I=[-1,1]$. For any $\sigma\in \Gamma_t^K,$ it suffices to prove that
\begin{align}\label{es:main0}
\int_I&\Big((\partial_ra_0)^2+\frac{1}{\ve^2}f_A(\theta_{2, \ve})a_0^2+(\partial_ra_i)^2+\frac{1}{\ve^2}f_B(\theta_{2, \ve})a_i^2\Big)J\ud r\notag\\
&+\underbrace{\frac{1}{\ve}\int_I \big(\rho_1f_C(\theta_{2, \ve})a_0^2+\rho_1 f_D(\theta_{2, \ve})a_i^2\big)J\ud r}_{\text{correction terms}}\notag\\
&+\underbrace{2\sum_{i=1}^{n-1}\int_I(\partial_ra_0a_i-\partial_ra_ia_0)E_0\cdot \partial_rE_iJ\ud r}_{\text{cross terms}}\geq -C\int_I(a_0^2+a_i^2)J\ud r.
\end{align}
We introduce the quadratic forms for $b\in H^1(I)$:
\begin{equation}\label{bilinear}
\mathcal{Q}_i(b)=
\left\{
\begin{array}{ll}
\int_I|\partial_rb|^2+\frac{1}{\ve^2}f_A(\theta_{2,\ve})b^2\ud r,   &i=0;\\
\int_I |\partial_rb|^2+\frac{1}{\ve^2}f_B(\theta_{2,\ve})b^2\ud r,   &i=1.
\end{array}
\right.
\end{equation}
The key step in our estimates will be to bound these two quadratic forms.
\subsection{Estimate for $\mathcal{Q}_i$ for $i=0,1$}\label{chap4.2}
In the sequel, we will assume that $\ve$ is sufficiently small.
 About the  estimate for $\mathcal{Q}_{i}$, we can refer to \cite{C1994,FLWZ,FWZZ2015,FWZZ2018}.
\begin{lemma}\label{le:0}
Let $b_0=\theta_{1,\ve}\bar{b}_0$, for any $\nu_0 >0$, there exists $C(\nu_0)>0$ such that
\begin{align}\label{es:0}
\mathcal{Q}_0(b_0)\geq \Big(\frac{1}{2}-\frac{\nu_0}{2}\Big)\int_I (\theta_{1,\ve}\partial_r \bar{b}_0)^2\ud r-C(\nu_0)e^{-\frac{2\alpha_0}{\ve}}\frac{1}{\ve^2}\int_I b_0^2\ud r.
\end{align}
\proof
Since $\theta_{1,\ve}$ satisfies
$$\ve^2\partial_r^2\theta_{1, \ve}=\theta_{1, \ve}f_A(\theta_{2, \ve}), $$
one has
\begin{align*}
\mathcal{Q}_0(b_0)&=\int_I \big([\partial_r\theta_{1,\ve}\bar{b}_0+\theta_{1,\ve}\partial_r\bar{b}_0]^2+\ve^{-2}f_A(\theta_{2,\ve})\theta^2_{1,\ve}\bar{b}_0^2\big)\ud r\nonumber\\
&=\theta_{1,\ve}\partial_r\theta_{1,\ve}\bar{b}_0^2|_{-1}^1+\int_I \theta^2_{1,\ve}(\partial_r\bar{b}_0)^2\ud r.
\end{align*}
We have by \eqref{rho0expression} that
\begin{align*}
b^2-\rho_0^2&=(b+\rho_0)^2\Big(\frac{\rho_0-a}{\rho_0+a}\Big)^{\frac{b}{a}}e^{-\sqrt{2}b(b^2-a^2) z},\\
\rho_0^2-a^2&=(a+\rho_0)^2\Big(\frac{b-\rho_0}{b+\rho_0}\Big)^{\frac{a}{b}}e^{\sqrt{2}a(b^2-a^2) z}.
\end{align*}
Together with \eqref{pazrho0}, we have
\begin{align}
|\theta_{1,\ve}\partial_r\theta_{1,\ve}(1)|&=\frac{1}{\ve}\sqrt{2}\theta_{2,\ve}|(a^2+b^2-2\theta_{2,\ve}^2)|\theta_{1,\ve}^2(1)\notag\\
&=\frac{1}{2\ve}\sqrt{2}\theta_{2,\ve}|(a^2+b^2-2\theta_{2,\ve}^2)|(\theta_{2,\ve}^2-a^2)^2(\theta_{2,\ve}^2-b^2)^2(1)\notag\\
&\leq \frac{c_1}{\ve}e^{-\frac{2\kappa_+}{\ve}}\Big(1+O(e^{-\al_0/\ve})\Big),\label{theta1par}
\end{align}
and
\begin{align*}
\int_{0}^{1}\theta_{1,\ve}^{-2}\ud r
&=2\int_{0}^{1}
(\theta_{2,\ve}^2-a^2)^{-2}
(\theta_{2,\ve}^2-b^2)^{-2}\ud r\\
&=\frac{2}{(b^2-a^2)^2}
\int_{0}^{1}(b^2-\theta_{2,\ve}^2)^{-2}\ud r
\Big(1+O(e^{-\kappa_+/\ve})\Big)\\
&=\frac{1}{8(b^2-a^2)^2b^4}
\left(\frac{b-a}{b+a}\right)^{-2b/a}
\int_{0}^{1}e^{2\kappa_+r/\ve}\ud r
\Big(1+O(e^{-\kappa_+/\ve})\Big)\\
&=\frac{\ve}{16\kappa_+(b^2-a^2)^2b^4}
\left(\frac{b-a}{b+a}\right)^{-2b/a}
e^{2\kappa_+/\ve}
\Big(1+O(e^{-\kappa_+/\ve})\Big)\\
&=\frac{\ve}{2c_1}e^{\frac{2\kappa_+}{\ve}}
\Big(1+O(e^{-\al_0/\ve})\Big).
\end{align*}
where
$$\kappa_+=\sqrt{2}(b^2-a^2)b,\ c_1=8\sqrt{2}b^5(b^2-a^2)^3\big(\frac{b-a}{b+a}\big)^{\frac{2b}{a}}\text{ and }0<\alpha_0 <\min \{\sqrt{2}(b^2-a^2)b,\sqrt{2}(b^2-a^2)a\}.$$
Notice there exists $a'\in [0,\ve]$ such that
\begin{align*}
\bar{b}_0(a')^2&\leq \ve^{-1}
\int_{0}^{\ve}\bar{b}_0^2\ud r.
\end{align*}
Together with Young's inequality and the Cauchy--Schwarz
inequality, we obtain
\begin{align*}
\bar b_0(1)^2&=
\Big(\int_{a'}^{1}\partial_r\bar{b}_0\ud r+\bar{b}_0(a')\Big)^2\notag\\
&\leq
\left(1+\frac{\nu_0}{2}\right)
\left(\int_{a'}^1\partial_r\bar b_0\ud r\right)^2
+\left(1+\frac{2}{\nu_0}\right)\bar b_0(a')^2\\
&\leq
\left(1+\frac{\nu_0}{2}\right)
\left(\int_0^1\theta_{1,\ve}^{-2}\ud r\right)
\left(\int_0^1\theta_{1,\ve}^2
(\partial_r\bar b_0)^2\ud r\right)
+\frac{C(\nu_0)}{\ve}
\int_0^\ve\theta_{1,\ve}^2\bar b_0^2\ud r\\
&\leq
\left(1+\frac{\nu_0}{2}\right)
\frac{\ve}{2c_1}e^{2\kappa_+/\ve}
\Big(1+O(e^{-\alpha_0/\ve})\Big)
\int_0^1\theta_{1,\ve}^2
(\partial_r\bar b_0)^2\ud r
+\frac{C(\nu_0)}{\ve}
\int_0^\ve\theta_{1,\ve}^2\bar b_0^2\ud r.
\end{align*}
By \eqref{theta1par}, we obtain
\begin{align}\label{es:0fac}
|\theta_{1, \ve}\partial_r\theta_{1, \ve}\bar{b}_0^2(1)|\leq\big(\frac{\nu_0}{2}+\frac{1}{2}\big)\int_0^1(\theta_{1,\ve}\partial_r\bar{b}_0)^2\ud r+C(\nu_0)e^{-\frac{2\alpha_0}{\ve}}\frac{1}{\ve^2}\int_0^\ve\theta_{1,\ve}^2\bar{b}_0^2\ud r.
\end{align}
Similarly, we have
\begin{align*}
   |\theta_{1, \ve}\partial_r\theta_{1, \ve}\bar{b}_0^2(-1)|\leq\big(\frac{\nu_0}{2}+\frac{1}{2}\big)\int_{-1}^0(\theta_{1,\ve}\partial_r\bar{b}_0)^2\ud r+C(\nu_0)e^{-\frac{2\alpha_0}{\ve}}\frac{1}{\ve^2}\int_{-\ve}^0\theta_{1,\ve}^2\bar{b}_0^2\ud r.
\end{align*}
Then we get \eqref{es:0}.
\qed
\end{lemma}
\begin{lemma}\label{le:n-1}
Let $b_1=\theta_{2,\ve}\bar{b}_1$. Then for any $\nu_0>0$, there exists $C(\nu_0)>0$ such that
\begin{equation}\label{q1b1}
\mathcal{Q}_1(b_1) \geq \Big(\frac{1}{2}-\nu_0\Big) \int_I (\theta_{2,\ve}\partial_r \bar{b}_1)^2\ud r-\frac{C(\nu_0)}{\ve} e^{-\frac{\alpha_0}{\ve}}\int_I b_1^2\ud r.
\end{equation}
\proof
By
$$\ve^2\partial_r^2\theta _{2,\ve}=\theta _{2,\ve}f_B(\theta _{2,\ve})\text{ and }\theta_{2,\ve}\partial_r\theta_{2,\ve}\bar{b}_1^2(1)\ge 0,$$
we integrate by parts in $I$:
\begin{align*}
\mathcal{Q}_1(b_1)&=\int_I ([\partial_r\theta_{2,\ve}\bar{b}_1+\theta_{2,\ve}\partial_r\bar{b}_1]^2+\ve^{-2}f_B(\theta_{2,\ve})\theta^2_{2,\ve}\bar{b}_1^2)\ud r\nonumber\\
&=\theta_{2,\ve}\partial_r\theta_{2,\ve}\bar{b}_1^2|_{-1}^1+\int_I \theta^2_{2,\ve}(\partial_r\bar{b}_1)^2\ud r\nonumber\\
&\geq -(\theta_{2,\ve}\partial_r\theta_{2,\ve} \bar{b}_1^2)(-1)+\int_I \theta^2_{2,\ve}(\partial_r\bar{b}_1)^2\ud r.
\end{align*}
We deduce from \eqref{basicestimate} that
\begin{align*}
(\theta_{2,\ve}\partial_r\theta_{2,\ve})(-1)\leq \frac{c_2}{\ve} e^{-\frac{\al_0}{\ve}}\ \text{and}\
\frac{2}{b^2}\leq \int_{-1}^{1}\theta_{2,\ve}^{-2}\ud r \leq \frac{2}{a^2},
\end{align*}
where
$$c_2=2\sqrt{2}a^3(b^2-a^2)\Big(\frac{b-a}{b+a}\Big)^{\frac{a}{b}}.$$
We get by $\theta_{2,\ve}(r)\geq a,\;r\in[0,1]$ that there exists $a'\in [0,1]$ such that
\begin{align*}
\bar{b}_1(a')^2&\leq
\int_{0}^{1}\bar{b}^2_1\ud r\leq \frac{1}{a^2}\int^{1}_0\theta_{2,\ve}^2\bar{b}^2_1\ud r.
\end{align*}
 Choosing  $\ve$ sufficiently small, we obtain
\begin{align*}
\theta_{2, \ve}&\partial_r\theta_{2, \ve}\bar{b}_1^2(-1)\\
\leq& \theta_{2, \ve}\partial_r\theta_{2, \ve}(-1)\Big(\bar{b}_1(a')-\int_{-1}^{a'}\partial_r\bar{b}_1\ud r\Big)^2\\
\leq& \frac{c_2}{\ve} e^{-\frac{\al_0}{\ve}}\Big((1+\nu_1)\int_I(\theta_{2, \ve}\partial_r \bar{b}_1)^2\ud r+C(\nu_1)\int_I\theta _{2, \ve}^2\bar{b}_1^2\ud r\Big)\\
\leq& \Big(\f12+\nu_0\Big)\int_I (\theta_{2, \ve}\partial_r \bar{b}_1)^2\ud r+\frac{C(\nu_0)}{\ve} e^{-\frac{\alpha_0}{\ve}}\int_I (\theta_{2, \ve} \bar{b}_1)^2\ud r.
\end{align*}
The proof of \eqref{q1b1} is complete.
\qed
\end{lemma}
\subsection{Endpoint $L^\infty$ estimates}\label{chap4.3}
\begin{lemma}\label{le:infccontrol}
For any $\nu_0>0$, there exists $C(\nu_0)>0$ such that
\begin{align}\label{eqinf}
\|b\|^2_{L^{\infty}([-1,1])}\leq \nu_0\mathcal{Q}_1(b)+C(\nu_0)\int_I b^2\ud r.
\end{align}
\proof
Let $b=\theta_{2,\ve}\bar{b}$, then by $0<a<\theta_{2,\ve}<b$ and Gagliardo-Nirenberg inequality, it is clear that
\begin{align*}
\|\bar{b}\|_{L^\infty([-1,1])}^2\leq \frac{\nu_0}{4}\int_I\theta_{2,\ve}^2(\partial_r\bar{b})^2\ud r+C\int_I\theta_{2,\ve}^2\bar{b}^2\ud r.
\end{align*}
According to Lemma \ref{le:n-1}, we derive that
\begin{align*}
\int_I \theta^2_{2,\ve}(\partial_r\bar{b})^2\ud r\leq 4\mathcal{Q}_1(b)+\frac{C}{\ve}e^{-\frac{\alpha_0}{\ve}}\int_I b^2\ud r.
\end{align*}
Obviously, we obtain \eqref{eqinf} directly.
\qed
\end{lemma}

\begin{lemma}\label{le:endpoints}
There exists $C>0$ such that
\begin{align}\label{endp}
|b(\pm 1)|^2\leq C\ve\Big(\mathcal{Q}_0(b)+\int_I b^2\ud r\Big).
\end{align}
\proof
Let $b=\theta_{1,\ve}\bar{b}$, by \eqref{theta1par}, it is easy to prove
$$\theta_{1,\ve}^2(\pm 1)=\big(\rho_0'(\frac{r}{\ve})\big)^2(\pm 1)\leq C\ve |(\theta_{1,\ve}\partial_r\theta_{1,\ve})(\pm 1)|.$$
 It is apparent from \eqref{es:0} and \eqref{es:0fac} that \eqref{endp} is proved.
\qed
\end{lemma}
\subsection{Estimate for cross terms}\label{chap4.4}
Now, we estimate the cross terms:
$$\int_I (\partial_ra_0a_i-a_0\partial_ra_i)E_0\cdot\partial_rE_iJ\ud r\text{ for }1\leq i\leq n-1.$$
For convenience, set $a(r)=E_i\cdot \pa_rE_j$ and $b_i(r)=J^{\frac{1}{2}}a_i(r).$
We absorb the Jacobian   $J$ by setting $b_i=J^{1/2}a_i$:
\begin{align*}
(\partial_ra_0a_i-a_0\partial_ra_i)J=&\partial_r(J^{-\frac{1}{2}}b_0)b_iJ^{\frac{1}{2}}-\partial_r(J^{-\frac{1}{2}}b_i)b_0J^{\frac{1}{2}}\\
=&\partial_rJ^{-\frac{1}{2}}b_0b_iJ^{\frac{1}{2}}+\partial_rb_0b_i-J^{\frac{1}{2}}b_0\partial_rJ^{-\frac{1}{2}}b_i-b_0\partial_rb_i\\
=&\partial_rb_0b_i-b_0\partial_rb_i.
\end{align*}
The following Lemma is the key to control the cross terms.
\begin{lemma}\label{par}
There exists a constant $C>0$  independent of $\varepsilon$ and
$t\in[0,T]$ such that
$$|\pa_rE_0|+|\pa_r(\rho_0^2\pa_rE_0)|\leq C\text{ for }(x,t)\in \Gamma_t^K(\delta/4).$$
\begin{proof}
We  take $\sigma'(x,t)$ as the orthogonal projection from $(x,t)$ onto the interface $\Gamma_t$.
In $\Gamma_t^K(\delta/4)$, by
\[
u^K
=
\rho_0\bar\omega
+
\varepsilon
\left(
\rho_1\bar\omega
+
\sum_{\alpha=1}^{n-1}
\sigma_{1,\alpha}\xi_\alpha
\right)
+
O(\varepsilon^2),\quad |u^K|=\rho_0+\ve \rho_1+O(\ve^2),
\]
and
$\bar\omega\cdot\xi_\alpha=0$, we have
\begin{equation}\label{eq:E0-expansion-corrected}
E_0
=\frac{u^K}{|u^K|}=
\bar\omega+\varepsilon R+O(\varepsilon^2),
\qquad
R
:=
\sum_{\alpha=1}^{n-1}
\frac{\sigma_{1,\alpha}}{\rho_0}\xi_\alpha.
\end{equation}
 By \eqref{pazxigamma}, \eqref{v2solution} and \eqref{sigma1eq}, a direct calculation gives us
\begin{align*}
   \pa_zR(z,x,t)|_{\G}=&\sum_\al\pa_z\frac{\sigma_{1,\al}}{\rho_0}\xi_\al \\
   =&\rho_0^{-2}(z)\sum_\al\int_{-\infty}^{z}\mathcal{D}_1\cdot \xi_\al\rho_0(z)\ud z\ \xi_\al\\
   =&-\rho_0^{-2}(z)\sum_\al\int_{-\infty}^{z}\Big(2\rho_0'\pa_\nu\bar{\omega}\cdot\xi_\al+2\rho_0\pa_\nu\pa_z\bar{\omega}\cdot\xi_\al\\
   &-(\pa_z^2\pa_\nu\bar{\omega}\rho_0+2\pa_{\nu}\pa_z\bar{\omega}\pa_z\rho_0)\cdot \xi_\al (d_1-z)\Big)\rho_0(z)\ud z\ \xi_\al\\
   =&-\rho_0^{-2}(z)\sum_\al\Big(\rho_0^2\pa_\nu\bar{\omega}\cdot \xi_\al|_{-\infty}^{z}+\int_{-\infty}^z\rho_0^2\partial_{\nu}\pa_z\bar\omega\cdot\xi_\al \ud z\\
   &-\int_{-\infty}^z\pa_z(\pa_\nu\pa_{z}\bar{\omega}\rho_0^2)\cdot \xi_\al (d_1-z)\ud z\Big)\ \xi_\al\\
=&-\rho_0^{-2}(z)\sum_\al\Big(\rho_0^2\pa_\nu\bar{\omega}\cdot \xi_\al|_{-\infty}^{z}-\pa_z\pa_\nu \bar{\omega}\cdot\xi_\al \rho_0^2(d_1-z)|_{-\infty}^z\Big)\ \xi_\al\\
=&-\rho_0^{-2}(z)\sum_\al\Big(\rho_0^2(\pa_\nu\bar\omega-\pa_z\pa_\nu\bar{\omega}(d_1-z))\cdot \xi_\al-a^2\pa_\nu\omega^-\cdot \xi_\al^-\Big)\xi_\al\\
=&-\pa_\nu\bar{\omega}+\pa_z\pa_\nu\bar{\omega}(d_1-z)+a^2\rho_0^{-2}\pa_\nu\omega^-.
\end{align*}
We have by $\nabla d_0\cdot\nabla d_1=0$ that $d_1(x,t)=d_1(\sigma'(x,t),t)$, which yields
\begin{equation}\label{eq:correct-z-relation}
\frac{d_0}{\varepsilon}+d_1(\sigma'(x,t),t)-z=O(\varepsilon),\quad \pa_r\Big(\frac{d_0}{\varepsilon}+d_1(\sigma'(x,t),t)-z\Big)=O(\varepsilon).
\end{equation}
By \eqref{innerexpansion} and \eqref{parx}, we have
$$\frac{\nabla d^K}{|\nabla d^K|^2}\cdot\nabla=(1+O(\ve^{K+1}))\nabla d^K\cdot \nabla=\nabla d_0\cdot \nabla +O(\ve).$$
Using Taylor expansion at $\Gamma$, together with \eqref{pazxigamma}, \eqref{decay:bar-omega} and the matching conditions, we obtain
\begin{align*}
&\frac{\nabla d^K}{|\nabla d^K|^2}
\cdot\nabla_x\bar\omega(z,x,t)
=
\left.\partial_\nu\bar\omega\right|_\Gamma
+O(d_0)+O(\varepsilon),\\
&\partial_z\bar\omega(z,x,t)=d_0
\left.\partial_z\partial_\nu\bar\omega
\right|_\Gamma+O\left(d_0^2e^{-\alpha_0|z|}\right),\\
&\varepsilon
\frac{\nabla d^K}{|\nabla d^K|^2}
\cdot\nabla_xR(z,x,t)
=
O(\varepsilon),\\
&\partial_zR(z,x,t)
=
\left.\partial_zR\right|_\Gamma
+
O\left(d_0e^{-\alpha_0|z|}\right).
\end{align*}
Here all the quantities restricted to $\Gamma$ are evaluated at
$(z,\sigma'(x,t),t)$. It follows from \eqref{eq:correct-z-relation} that
\[
|d_0|
\leq
C\varepsilon(1+|z|),
\]
and hence
\[
\frac{d_0^2}{\varepsilon}e^{-\alpha_0|z|}
+
|d_0|e^{-\alpha_0|z|}
\leq
C\varepsilon.
\]
Together with \eqref{eq:E0-expansion-corrected}, this gives
\begin{align*}
\partial_rE_0
=&
\frac{\nabla d^K}{|\nabla d^K|^2}
\cdot\nabla_x\bar\omega+
\frac{1}{\varepsilon}\partial_z\bar\omega+\varepsilon\frac{\nabla d^K}{|\nabla d^K|^2}\cdot\nabla_xR+\partial_zR+O(\varepsilon)\\
=&\partial_\nu\bar\omega|_\Gamma+\frac{d_0}{\varepsilon}
\partial_z\partial_\nu\bar\omega|_\Gamma+\partial_zR|_\Gamma+O(d_0)+O(\varepsilon)\\
=&\partial_\nu\bar\omega|_\Gamma+\frac{d_0}{\varepsilon}
\partial_z\partial_\nu\bar\omega|_\Gamma-
\partial_\nu\bar\omega|_\Gamma+
\left(d_1(\sigma',t)-z\right)
\partial_z\partial_\nu\bar\omega|_\Gamma
+\frac{a^2}{\rho_0^2}
\partial_\nu\omega^-(\sigma',t)+O(d_0)+O(\varepsilon)
\\
=&
\left(\frac{d_0}{\varepsilon}+d_1(\sigma'(x,t),t)-z\right)\left.\partial_z\partial_\nu\bar\omega\right|_\Gamma+\frac{a^2}{\rho_0^2}
\partial_\nu\omega^-(\sigma'(x,t),t)+
O(d_0)+O(\varepsilon).
\end{align*}
By \eqref{def:eta1}-\eqref{decay:bar-omega},
\[|\partial_z\partial_\nu\bar\omega|_\G+|\partial_z^2\partial_\nu\bar\omega|_\Gamma+|\partial_z\rho_0|\leq
Ce^{-\alpha_0|z|}.\]
Consequently,
\[
\begin{aligned}
\Bigg|
\partial_r
\Bigg[
\rho_0^2
\left(
\frac{d_0}{\varepsilon}
+
d_1(\sigma'(x,t),t)-z
\right)
\left.
\partial_z\partial_\nu\bar\omega
\right|_\Gamma
\Bigg]
\Bigg|
\leq C.
\end{aligned}
\]
Using
\[
\partial_rd_0=1+O(\varepsilon),
\qquad
\partial_rz=\frac{1}{\varepsilon},
\qquad
\frac{|d_0|}{\varepsilon}e^{-\alpha_0|z|}
\leq
C(1+|z|)e^{-\alpha_0|z|}
\leq C,
\]
we have $
\left|
\partial_r\bigl(\rho_0^2O(d_0)\bigr)
\right|
\leq C.$
Similarly, we have $
\left|\partial_r\bigl(\rho_0^2O(\varepsilon)\bigr)
\right|\leq C.$
Since $d_0(\sigma'(x,t),t)=0,$ one has $\nabla d_0\cdot \pa_r\sigma'=0$ which yields
\begin{align*}
\left|\partial_r\left[a^2\partial_\nu\omega^-(\sigma'(x,t),t)\right]\right|=a^2|\nabla_\G \pa_\nu \omega^- \pa_r\sigma'|\leq C.
\end{align*}
Combining the above estimates, we conclude that
\begin{align*}
|\pa_rE_0|+\left|\partial_r\left(\rho_0^2\partial_rE_0\right)\right|\leq C.
\end{align*}
This completes the proof.
\end{proof}
\end{lemma}
\begin{lemma}\label{cross0ac}
For $i=1\cdots,n-1$. Then for any $\nu_0>0$, there exists $C(\nu_0)>0$ such that
\begin{align*}
\Big|\int_I (\partial_rb_0 b_i-b_0\partial_rb_i)E_0\cdot \pa_rE_i\ud r\Big|\leq \nu_0\big(\mathcal{Q}_0(b_0)+\mathcal{Q}_1(b_i)\big)+C(\nu_0)\int_I (b_0^2+b_i^2)\ud r.
\end{align*}
\proof
Letting $b_0=\theta_{1,\ve} \bar{b}_0$ and $b_i=\theta_{2,\ve}\bar{b}_i$, we have
\begin{align*}
\int_I &(\partial_rb_0 b_i-b_0 \partial_rb_i)E_0\cdot \pa_rE_i\ud r\\
=&\int_I (\pa_r\theta_{1, \ve}\theta_{2, \ve}-\theta_{1, \ve}\pa_r\theta_{2, \ve})\bar{b}_0\bar{b}_iE_0\cdot\pa_rE_i\ud r\\
&+\int_I \theta_{1, \ve}\theta_{2, \ve}(\pa_r\bar{b}_0\bar{b}_i-\bar{b}_0\pa_r\bar{b}_i)E_0\cdot \pa_rE_i\ud r.
\end{align*}
By Lemma \ref{par},  for $i\geq 1$, we have
\begin{align*}
    E_0\cdot \pa_r E_i=-\pa_r E_0\cdot E_i=O(1).
\end{align*}
By Lemmas \ref{le:0} and \ref{le:n-1}, it is easy to obtain that
\begin{align}\label{cross1}
\Big|\int_I \theta_{1,\ve}\theta_{2,\ve}(\pa_r\bar{b}_0\bar{b}_i-\bar{b}_0\pa_r\bar{b}_i)E_0\cdot \pa_rE_i\ud r\Big|\leq \frac{\nu_0}{3}(\mathcal{Q}_0(b_0)+\mathcal{Q}_1(b_i))+C\int_I b_0^2+b_i^2\ud r.
\end{align}
By \eqref{eiej}, we obtain $\pa_r E_0\cdot\pa_rE_i=0$ for $i=1,\cdots,n-1.$
We get by $\theta_{1,\ve}=\ve\pa_r\theta_{2,\ve}$ and integration by parts that
\begin{align*}
\int_I (\pa_r\theta_{1,\ve}&\theta_{2,\ve}-\theta_{1,\ve}\pa_r\theta_{2,\ve})\bar{b}_0\bar{b}_i\pa_rE_0\cdot E_i\ud r\\
=&\int_I \ve \theta_{2,\ve}^2 \pa_r\Big(\frac{\pa_r\theta_{2,\ve}}{\theta_{2,\ve}}\Big)\bar{b}_0\bar{b}_i\pa_rE_0\cdot E_i\ud r\\
=&-\ve\int_I \pa_r\theta_{2,\ve}\pa_r\bar{b}_0\bar{b}_i\theta_{2,\ve}\pa_rE_0\cdot E_i\ud r-\ve\int_I \pa_r\theta_{2,\ve}\bar{b}_0\pa_r\bar{b}_i\theta_{2,\ve}\pa_rE_0\cdot E_i\ud r\\
&-\ve\int_I \frac{\pa_r\theta_{2,\ve}}{\theta_{2,\ve}}\bar{b}_0\bar{b}_i\pa_r(\theta_{2,\ve}^2\pa_rE_0)\cdot E_i\ud r+\ve \theta_{2,\ve}\pa_r\theta_{2,\ve}\bar{b}_0\bar{b}_i\pa_rE_0\cdot E_i|_{-1}^1:=\sum_{i=1}^4I_i.
\end{align*}
It is easy to obtain that
\begin{align}\label{124}
|I_1|+|I_2|\leq \frac{\nu_0}{3}(\mathcal{Q}_0(b_0)+\mathcal{Q}_1(b_i))+C\int_I b_0^2+b_i^2\ud r
\end{align}
using Lemmas \ref{le:0} and \ref{le:n-1}. We deduce from Lemma \ref{par} that
\begin{align}\label{3}
|I_3|\leq C\int_I b_0^2+b_i^2\ud r.
\end{align}
By Lemma \ref{le:infccontrol}, Lemma \ref{le:endpoints} and Young's inequality, we have
\begin{align}\label{4}
    |I_4|=\big|b_0b_i\pa_rE_0\cdot E_i|_{-1}^1\big|\leq \frac{\nu_0}{3}(\mathcal{Q}_0(b_0)+\mathcal{Q}_1(b_i))+C\int_I b_0^2+b_i^2\ud r.
\end{align}
Combining \eqref{cross1}-\eqref{4}, we complete the proof.
\qed
\end{lemma}
\subsection{Estimate for correction terms}\label{chap4.5}
It remains to estimate
$$\int_I \frac{1}{\ve}\rho_1f_C(\theta_{2,\ve})b_0^2\ud r\ \ \text{and}\ \int_I \frac{1}{\ve}\rho_1f_D(\theta_{2,\ve})b_i^2\ud r\ \text{for}\ i\geq 1.$$
\begin{lemma}\label{es:rho1}
There exists a constant $C>0$ such that
\begin{align}\label{rho1fc}
\big|\int_I \frac{1}{\ve}\rho_1f_C(\theta_{2,\ve})b_0^2\ud r\big|\leq C\int_I b_0^2\ud r,
\end{align}
and
\begin{align}\label{rho1fd}
\big|\int_I \frac{1}{\ve}\rho_1f_D(\theta_{2,\ve})b_i^2\ud r\big|\leq C\int_I b_i^2\ud r,\text{ for } i\geq 1.
\end{align}
\proof
We have by \eqref{u1u2star} and \eqref{rho1} that $\rho_1(z,x,t)$ decays exponentially as $z\to \pm\infty$. Using \eqref{rho1gamma}, we have
\begin{align*}
\Big|\frac{1}{\ve}\rho_1(\frac{r}{\ve},x,t)\Big|\leq C\frac{|d_0|}{\ve}e^{-\frac{\alpha_0|r|}{\ve}}\leq C.
\end{align*}
The bound follows for \eqref{rho1fc}. Similarly, we can get \eqref{rho1fd}.
\qed
\end{lemma}
Then by using Lemmas \ref{le:infccontrol} and \ref{le:endpoints} to control $|b_i(\pm 1)|$, we obtain
\begin{align*}
\int_I |\partial_r a_i|^2J\ud r&=\int_I |\partial_r(J^{-\frac{1}{2}}b_i)|^2J\ud r\\
&=\int_I |\partial_rb_i|^2+\big[|\partial_r(J^{-\frac{1}{2}})|^2J-\partial_r(\partial_rJ^{-\frac{1}{2}}J^{\frac{1}{2}})\big]b_i^2\ud r+\partial_rJ^{-\frac{1}{2}}J^{\frac{1}{2}}b_i^2|_{-1}^1\\
&\geq \int_I |\partial_rb_i|^2\ud r-\frac{1}{4(n-1)}\mathcal{Q}_{j(i)}(b_i)-C\int_I b_i^2\ud r,
\end{align*}
where
\begin{equation*}
   j(i)= \begin{cases}
        0,\ i=0,\\
        1,\ i\in\{1,\cdots,n-1\}
    \end{cases}
\end{equation*}
It is enough to estimate \eqref{es:main0}. Theorem \ref{th:uA} is derived directly.

\section{Uniform error estimate}
Let $u^K$ be an approximate solution constructed in Theorem \ref{thK-1} and $u^\ve$ be a solution of \eqref{equation:main1}. Define
$$u_R=\frac{u^\ve-u^K}{\ve^k}.$$
Then  $u_R$ satisfies
\begin{align}\label{maineq}
\partial_tu_R=\Delta u_R-\frac{1}{\ve^2}\big(G'(|u^K|^2)u_R+2G''(|u^K|^2)u^K\cdot u_Ru^K\big)+\ve^{k-2}
\hat{H}(u^K,u_R)+\mathcal{R}_K,
\end{align}
where $\mathcal{R}_K=-\ve^{-k}\mathcal{R}^K$,
\begin{align*}
G'(s)&=(s-a^2)(s-b^2)(2s-a^2-b^2),\\
G''(s)&=(2s-a^2-b^2)^2+2(s-a^2)(s-b^2),\\
\hat{H}(u^K,u_R)&=f_1(u^K,u_R)+\ve^{k}f_2(u^K,u_R)+\cdots+\ve^{5k}f_6(u_R).
\end{align*}
Here
\[
f_\ell(u^K,u_R)
=-\frac{1}{(\ell+1)!}
D^{\ell+1}f(u^K)
[\underbrace{u_R,\ldots,u_R}_{\ell+1\text{ factors}}],
\qquad 1\leq\ell\leq6.
\]
Thus $f_\ell$ is homogeneous of degree $\ell+1$
in $u_R$, and $f_6$ is independent of $u^K$. Direct calculation gives
\[
\|\partial^j u^K\|_{L^\infty}
\leq C\varepsilon^{-j},
\qquad
\|\partial^j \mathcal{R}_K\|_{L^2}
\leq C\varepsilon^{K-k-1-j},
\qquad
0\leq j\leq[\frac m2]+1:=\tau.
\]

Set
$$\mathcal{E}(u)=\sum_{i=0}^{\tau}\ve^{6i}\int_{\Omega}\|\partial^iu\|^2\ud x.$$
Then
\[
\|u_R\|_{H^q}\leq C\varepsilon^{-3q}\mathcal{E}^{1/2},
\qquad 0\leq q\leq \tau,
\qquad
\|u_R\|_{L^\infty}\leq C\varepsilon^{-3\tau}\mathcal{E}^{1/2}.
\]

Choose $K\geq k+1$ and $k=3\tau+3$.
Multiplying (\ref{maineq}) by $u_R$ and integrating over $\Omega$, we find that
\begin{align*}
\frac{1}{2}\frac{d}{dt}\int_{\Omega}|u_R|^2\ud x=&-\int_{\Omega}|\nabla u_R|^2 \ud x-\frac{1}{\ve^2}\int_{\Omega}G'(|u^K|^2)|u_R|^2+2G''(|u^K|^2)(u^K\cdot u_R)^2\ud x\\
&+\ve^{k-2}\int_{\Omega}\hat{H}(u^K,u_R)\cdot u_R\ud x+\int_\Omega\mathcal{R}_K\cdot u_R\ud x.
\end{align*}
With the help of Theorem \ref{th:uA}, standard energy estimates yield that
\begin{align*}
\frac{d}{dt}\int_{\Omega}|u_R|^2\ud x&\leq C(1+\ve\mathcal{E}^{\frac{1}{2}}+\ve^4\mathcal{E}+\ve^7\mathcal{E}^{\frac{3}{2}}+\ve^{10}\mathcal{E}^2+\ve^{13}\mathcal{E}^{\frac{5}{2}}+\ve^{16}\mathcal{E}^3)\int_{\Omega}|u_R|^2\ud x+C\nonumber\\
&\leq C(1+\mathcal{E}+\ve \mathcal{E}^{\frac{3}{2}}+\ve^{13}\mathcal{E}^{\frac{7}{2}}+\ve^{16}\mathcal{E}^4).
\end{align*}
Applying $\partial^i(0\leq i\leq \tau)$ to the equation \eqref{maineq}, it is apparent that
\begin{align*}
\partial_t\partial^iu_R=&\Delta \partial^iu_R-\frac{1}{\ve^2}\partial^i(G'(|u^K|^2)u_R+2G''(|u^K|^2)u^K\cdot u_Ru^K)+\ve^{k-2}\partial^i\hat{H}(u^K,u_R)+\partial^i\mathcal{R}_K,
\end{align*}
where $\partial^i\mathcal{R}_K=O(\ve^{K-1-k-i})$.  Theorem \ref{th:uA} yields that
\begin{align}\label{error1}
\int_\Omega|\nabla \partial^iu_R|^2+\frac{1}{\ve^2}(G'(|u^K|^2)|\partial^iu_R|^2+2G''(|u^K|^2)|u^K\cdot \partial^iu_R|^2)\ud x\geq -C\int_\Omega|\partial^iu_R|^2\ud x.
\end{align}
The Leibniz rule  gives
\[
\|\partial^i f_\ell(u^K,u_R)\|_{L^2}
\leq C\|u_R\|_{L^\infty}^{\ell}
\sum_{j=0}^{i}\varepsilon^{-j}\|u_R\|_{H^{i-j}},
\qquad
0\leq i\leq \tau,\quad 1\leq\ell\leq6.
\]
In particular, since $k=3\tau+3$,
\begin{align}
\varepsilon^{k-2}
\|\partial^i f_1(u^K,u_R)\|_{L^2}
&\leq
C\varepsilon^{k-2}\|u_R\|_{L^\infty}
\sum_{j=0}^{i}\varepsilon^{-j}\|u_R\|_{H^{i-j}}\leq
C\varepsilon^{k-2-3\tau-3i}\mathcal{E}
\sum_{j=0}^{i}\varepsilon^{2j}\leq C\varepsilon^{1-3i}\mathcal{E}.
\end{align}
Similarly, for $2\leq\ell\leq6$,
\[
\varepsilon^{\ell k-2}
\|\partial^i f_\ell(u^K,u_R)\|_{L^2}
\leq
C\varepsilon^{\ell k-2-3\tau\ell-3i}\mathcal{E}^{(\ell+1)/2}
=
C\varepsilon^{3\ell-2-3i}\mathcal{E}^{(\ell+1)/2}.
\]

 So the $L^2-$norm of $$\frac{1}{\ve^2}\partial^i\big(G'(|u^K|^2)u_R+2G''(|u^K|^2)u^K\cdot u_Ru^K\big)$$ except for the term involved by the estimate \eqref{error1} can be controlled by
\begin{align}\label{error2}
\ve^{-2-j}\|\partial^qu_R\|_{L^2}\leq \ve^{-3i+3q}\|\partial^qu_R\|_{L^2}\leq\ve^{-3i}\mathcal{E}^{\frac{1}{2}}(u_R),
\end{align}
where $j+q=i$ and $q\leq i-1$.
It follows from \eqref{error1}-\eqref{error2} that
\begin{align}
\frac{d}{dt}\int_{\Omega}\ve^{6i}\|\partial^iu_R\|^2\ud x &\leq C(\mathcal{E}^{\frac{1}{2}}+\ve \mathcal{E}+\ve^{4}\mathcal{E}^{\frac{3}{2}}+\ve^7\mathcal{E}^2+\ve^{10}\mathcal{E}^{\frac{5}{2}}+\ve^{13}\mathcal{E}^3+\ve^{16}\mathcal{E}^{\frac{7}{2}})\ve^{3i}\|\partial^iu_R\|_{L^2}+C\nonumber\\
&\leq C(1+\mathcal{E}+\ve \mathcal{E}^{\frac{3}{2}}+\ve^{13}\mathcal{E}^{\frac{7}{2}}+\ve^{16}\mathcal{E}^4).
\end{align}
Summing $i$ from 0 to $[\frac{m}{2}]+1$, we have
\begin{align}\label{errormain}
\frac{d}{dt}\mathcal{E}(u_R)\leq C(1+\mathcal{E}+\ve \mathcal{E}^{\frac{3}{2}}+\ve^{13}\mathcal{E}^{\frac{7}{2}}+\ve^{16}\mathcal{E}^4).
\end{align}
After we have these estimates, we prove Theorem \ref{th:error} by a direct continuation argument.
\newpage
\appendix
\section{The basic information for $\omega$}
\begin{lemma}\label{lem:frame-product}
  Assume that \(\omega(s)\;(0\le s\le 1)\) is a geodesic on the unit sphere with arc-length parameter \(s\), and \(\xi(s)\;(0\le s\le 1)\) is a tangent vector field on the unit sphere obtained by parallel transporting \(\xi(0)\) along \(\omega(s)\). Then the following functions are independent of $s$:
 \begin{align*}
   |\xi'(s)|,\quad \omega'(s)\cdot\xi(s), \quad  \omega(s)\cdot\xi'(s), \quad \omega'(s)\cdot \xi'(s), \quad \omega''(s)\cdot \xi(s).
 \end{align*}
 Moreover, if \(\tilde{\xi}(s)\;(0\le s\le 1)\) is another tangent vector field on the unit sphere obtained by parallel transporting \(\tilde\xi(0)\) along \(\omega(s)\) with $\tilde\xi(0)\cdot\xi(0)=0$, then $\partial_\tau \xi(\tau)\cdot\tilde{\xi}(\tau)\equiv 0$.
 \end{lemma}
 \begin{proof}
From the assumptions, we have \(|\omega'(s)|=1\), \(\omega(s)\perp \omega'(s)\), and \(\xi(s)\perp \omega(s)\) with the covariant derivative condition \(\frac{D}{Ds}\xi=0\).

In the embedding \(\mathbb{R}^{n}\) of the sphere, the covariant derivative is given by the projection
\[
0=\frac{D}{Ds}\xi = \xi'(s) - \bigl(\xi'(s)\cdot \omega(s)\bigr)\omega(s),
\]
so \(\xi'(s)\) is parallel to \(\omega(s)\). Hence there exists a function \(\kappa(s)\) such that \(\xi'(s)=\kappa(s)\omega(s)\).

Differentiating \(\xi(s)\cdot \omega(s)=0\) gives
\[
\xi'(s)\cdot \omega(s)+\xi(s)\cdot \omega'(s)=0 \quad\Longrightarrow\quad \kappa(s) = -\xi(s)\cdot \omega'(s).
\]

\begin{enumerate}
\item \(|\xi'(s)| = |\kappa(s)| = |\xi(s)\cdot \omega'(s)|\).
Since \(\xi(s)\) and \(\omega'(s)\) are both parallel vector fields along the curve (the geodesic condition ensures \(\omega'(s)\) is parallel), their inner product \(\xi(s)\cdot \omega'(s)\) is constant. Hence \(|\xi'(s)|\) is independent of \(s\).

\item \(\omega(s)\cdot \xi'(s) = \omega(s)\cdot(\kappa(s)\omega(s)) = \kappa(s)\), which is constant by (1).

\item \(\omega'(s)\cdot \xi(s) = \xi(s)\cdot \omega'(s) = -\kappa(s)\), also constant.

\item \(\omega'(s)\cdot \xi'(s) = \omega'(s)\cdot(\kappa(s)\omega(s)) = \kappa(s)(\omega'(s)\cdot \omega(s)) = 0\), which is trivially constant.

\item \(\omega''(s)\cdot \xi(s) =-\omega'(s)\cdot \xi'(s)=0\) which is also constant.
\end{enumerate}
\end{proof}
\section{Proofs for the inner expansion}
\subsection{The proof of Lemma \ref{rerho}}\label{rho0zpro}
\begin{proof}
Since $\rho_0''=f(\rho_0)$, $\rho_0(-\infty)=a$,
$\rho_0(+\infty)=b$, and $f(a)=f(b)=0$, we have
\[
\lim_{z\to\pm\infty}\rho_0''(z)=0.
\]
By Taylor's formula with integral remainder,
\[
\rho_0(z+1)-\rho_0(z)
=\rho_0'(z)+\int_0^1(1-s)\rho_0''(z+s)\,ds.
\]
Hence,
\[
|\rho_0'(z)|
\leq |\rho_0(z+1)-\rho_0(z)|
+\frac12\sup_{s\in[0,1]}|\rho_0''(z+s)|
\longrightarrow 0
\qquad\text{as }z\to\pm\infty.
\]
Therefore, $\rho_0'(\pm\infty)=0$.
    Multiplying $\pa_\tau^2\rho_0(\tau)=f(\rho_0)$ by $\rho_0'(\tau)$ and integrating the equation from $-\infty$ to $z$ gives that
\begin{align}\label{rho0'2}
(\rho_0'(z))^2=2\int_{-\infty}^{z}f(\rho_0)\rho_0'\ud \tau=\f12(\rho_0^2-a^2)^2(b^2-\rho_0^2)^2. \end{align}
If there exists  $z_0\in \mathbb{R}$ such that $\rho_0(z_0)=a$, then we have $\rho_0'(z_0)=0$ and $f(\rho_0(z_0))=0$ which yields $\rho_0(z)=a$. This contradicts with $\rho_0(+\infty)=b$. Hence we have $a<\rho_0<b$, together with \eqref{rho0'2}, we have $\rho_0'>0$.
Then we have achieved \eqref{pazrho0}.

A direct calculation gives \eqref{rho0expression} and \eqref{dexpression} by \eqref{pazrho0}.
Let $A=b^2-\rho_0^2>0$ and $B=\rho_0^2-a^2$. With the help of \eqref{pazrho0}, we have
\begin{align*}
\lim_{z\to +\infty}\frac{(A')^2}{A^2}&=\lim_{z\to +\infty}\frac{4\rho_0^2(\rho_0')^2}{(b^2-\rho_0^2)^2}\\
&=\lim_{z\to +\infty} \frac{2\rho_0^2(\rho_0^2-a^2)^2(b^2-\rho_0^2)^2}{(b^2-\rho_0^2)^2}\\
&=2b^2(b^2-a^2)^2,
\end{align*}
and
\begin{align*}
\lim_{z\to -\infty}\frac{(B')^2}{B^2}=2a^2(b^2-a^2)^2.
\end{align*}
Then we have $A\le Ce^{-\al_0z}$ for $z>0$ and $B\le Ce^{\al_0z}$ for $z<0$ which yields $\rho_0'\le CAB\le Ce^{-\al_0|z|}$.
 Differentiating \eqref{pazrho0} with respect to $z$ gives the desired estimate \eqref{basicestimate}.
\end{proof}
\subsection{The proof of Lemma \ref{lem:s1}}\label{prooflemma3.1}
\proof We just consider the case of $i=1$. The case $i=2$ can be proved similarly.
Let $v_1$ be the solution of \eqref{ODE:A}. By \eqref{v1eq}, we have
\begin{align*}
\int_{\mathbb{R}}h_1(z,x,t)\theta_1(z)\ud z&=\int_{\mathbb{R}}(\partial_{z}^2v_1-f_A(\rho_0)v_1)\rho_0'\ud z=\int_{\mathbb{R}} \pa_z\Big((\rho_0')^2\pa_z(\frac{v_1}{\rho_0'})\Big) \ud z=0.
\end{align*}
Let $v_1=\theta_{1}\Phi$. By \eqref{v1eq}, we have
\begin{align*}
\partial_z^2 v_1-f_A(\rho_0(z)) v_1=\frac{1}{\rho_0'}\partial_z\Big((\rho_0')^2\partial_z\Big(\frac{v_1}{\rho_0'}\Big)\Big)=\frac{1}{\rho_0'}\partial_z\big((\rho_0')^2\partial_z\Phi\big).
\end{align*}
Thus the equation \eqref{ODE:A} is equivalent to
$$\partial_z\big((\rho_0')^2\partial_z\Phi\big)=h_1\rho_0'.$$
Integrating the above equation with respect to $z$ from $\varsigma$ to $+\infty$, it follows that
\begin{align*}
\partial_\varsigma\Phi=-\frac{1}{(\rho_0')^2(\varsigma)}\int_{\varsigma}^{+\infty}\rho_0' h_1(\tau,x,t)\ud \tau.
\end{align*}
By integrating the above equation with respect to $\varsigma$ from $0$ to $z$, we obtain
\begin{align*}
 \Phi(z)=-\int_{0}^{z}\frac{1}{(\rho_0')^2(\varsigma)}\int_{\varsigma}^{+\infty}\rho_0' h_1(\tau,x,t)\ud \tau \ud \varsigma+q_1(x,t).
\end{align*}
Using $v_1=\theta_1\Phi,$ we have
$$ v_1=-\rho_0'(z)\int_{0}^{z}\frac{1}{(\rho_0')^2(\varsigma)}\int_{\varsigma}^{+\infty}\rho_0' h_1(\tau,x,t)\ud \tau \ud \varsigma+\rho_0'(z)q_1(x, t)=v_1^*+\rho_0'(z)q_1(x, t). $$
We infer that $\int_{0}^{z}\frac{1}{(\rho_0')^2(\varsigma)}\int_{\varsigma}^{+\infty}(\rho_0' h_1)(\tau)\ud \tau \ud \varsigma$ is meaningful because of the condition $$\int_{\mathbb{R}}h_1(z,x,t)\rho_0'\ud z=0.$$ A routine computation gives rise to $$v^*_{1}(\cdot, x,t)\in \mathcal{S}_{J+2,L,M}(\alpha_0,k+1).$$
Taking $z\to\pm\infty$ in \eqref{ODE:A}, we gain $v_1^*(\pm\infty,x,t)$.
\qed
\subsection{The proof of Lemma \ref{led1}}\label{A.1}
\begin{proof}
Since $H_0^\pm=0$ in $U_\pm$ and the smoothness of $H_0^\pm$, it follows that $H_0^\pm|_\Gamma=0.$
Therefore,
\[
\left.
\bigl(H_0^+\eta_M^++H_0^-\eta_M^-\bigr)
\right|_\Gamma=0,
\]
On $\Gamma$, by \eqref{solvablecondition}, we have $\langle\mathcal{D}_2\cdot\bar{\omega},\rho_0'\rangle=\langle\mathcal{D}_2,\pa_zu_0\rangle=0$ which determines $d_1$ on $\Gamma$. First, we consider
\begin{align*}
&  \langle \mathcal{D}_1, \partial_z u_1\rangle+\langle\mathcal{D}_2, \partial_z u_0 \rangle  \\
&= \langle\partial_zu_0(\partial_td_0-\Delta d_0)-2\nabla d_0\cdot\nabla\partial_z u_0 +g_0(d_1-z), \partial_z u_1 \rangle\\
&\quad + \big\langle  \partial_zu_0(\partial_td_{1}-\Delta d_{1})+\partial_zu_{1}(\partial_td_{0}-\Delta d_{0})-2\nabla d_{1}\cdot \nabla \partial_zu_0-2\nabla d_{0}\cdot \nabla \partial_zu_{1}\nonumber\\
&\quad +f^{(1)}(u_0,u_{1})+d_2g_0+(d_1-z)g_{1}+\partial_{t}u_{0}-\Delta u_{0}, \partial_z u_0\big\rangle \\
&= (\partial_td_{1}-\Delta d_{1})\langle  \partial_zu_0,~\partial_zu_0\rangle -2\langle\nabla d_0\cdot\nabla\partial_z u_0 , \partial_z u_1 \rangle\\
&\quad -\nabla d_{1}\cdot\nabla \langle \partial_zu_0, \partial_z u_0\rangle -2\langle \nabla d_{0}\cdot\nabla \partial_zu_{1}, \partial_z u_0\rangle\nonumber\\
&\quad +\langle f^{(1)}(u_0,u_{1})+\partial_{t}u_{0}-\Delta u_{0}, \partial_z u_0\big\rangle +\langle g_0(d_1-z), \partial_z u_1 \rangle
+d_2\langle g_0,~\partial_z u_0\rangle+\langle (d_1-z)g_{1},~\partial_z u_0\big\rangle\\
&= (\partial_td_{1}-\Delta d_{1})\langle  \partial_zu_0,~\partial_zu_0\rangle -2\nabla d_0\cdot\nabla\langle\partial_z u_0 , \partial_z u_1 \rangle+\langle \partial_{t}u_{0}-\Delta u_{0}, \partial_z u_0\rangle\\
&\quad +\langle f^{(1)}(u_0,u_{1}), \partial_z u_0\big\rangle +\big(\langle g_0(d_1-z), \partial_z u_1 \rangle
+\langle (d_1-z)g_{1},~\partial_z u_0\big\rangle\big)\\
&:=\sum_{i=1}^5I_i.
\end{align*}
Here, we have used that on $\Gamma$:
$\partial_td_0-\Delta d_0=0,\ g_0\cdot \bar\omega =0.$

For $I_1$, one has
\begin{align*}
    I_1= e (\pa_t d_1-\Delta d_1),
\end{align*}
where $e=\langle \rho_0',\rho_0'\rangle.$
By \begin{align*}
    \bar{\omega}\cdot\xi_\al=0,\ \pa_z\bar{\om}\cdot\xi_\al=-\bar{\omega}\cdot \pa_z\xi_\al=-\kappa_\al (x,t)\eta',\ \pa_z\bar{\omega}\cdot\pa_z\xi_\al=0,\ |\pa_z\bar{\omega}|^2=\kappa_0^2(x,t)(\eta')^2,
\end{align*}
one has
\begin{align*}
    \langle\pa_zu_0,\pa_zu_1\rangle=&\langle \pa_z(\rho_0\bar{\omega}),\pa_z(\rho_1\bar{\omega}+\sum_\al \sigma_{1,\al}\xi_\al)\rangle\\
    =&\langle\rho_0'\bar{\omega}+\rho_0\pa_z\bar{\omega},\rho_1'\bar{\omega}+\rho_1\pa_z\bar{\omega}\rangle+\sum_\al\langle\rho_0'\bar{\omega}+\rho_0\pa_z\bar{\omega},\sigma_{1,\al}'\xi_\al +\sigma_{1,\al}\partial_z\xi_\al\rangle\\
    =&\langle \rho_0',\rho_1'\rangle +\langle \rho_0 \pa_z\bar{\omega},\rho_1\pa_z\bar{\omega}\rangle +\sum_\al\langle \rho_0\pa_z\bar{\omega},\sigma_{1,\al}'\xi_\al\rangle +\langle \rho_0'\bar{\omega},\sigma_{1,\al}\pa_z\xi_\al\rangle\\
    =&\langle \rho_0',\rho_1'\rangle +\kappa_0^2\langle |\eta'|^2\rho_0,\rho_1\rangle +\sum_\al\kappa_\al (-\langle \rho_0\eta',\sigma_{1,\al}'\rangle +\langle \rho_0'\eta',\sigma_{1,\al}\rangle).
\end{align*}
 Multiplying \eqref{rho1equation} by $z\rho_0'$, and using Lemma \ref{v1product}, we have
\begin{align*}
&\langle\rho_0',\rho_1'\rangle=-\f12\langle z\rho_0',L_{1,0}\eta' d_0+\mathcal{D}_1\cdot \bar{\omega}\rangle= -\f12d_0L_{1,0}\langle z\rho_0',\eta'\rangle-\f12\langle z\rho_0',\mathcal{D}_1\cdot \bar{\omega}\rangle.
\end{align*}
By
\begin{align}
 &\sigma_{1,\al}=-\frac{\kappa_\al d_1}{d_0}\eta(z)\rho_0(z)+\sigma^*_{1,\al}(z,x,t)+h_{1,\al}d_0\rho_0\eta+
 e_{1,\al}(x,t)\rho_0(z),\label{sigma1alg}\\
 &\pa_z\sigma_{1,\al}=-\frac{\kappa_\al d_1}{d_0}(\eta\rho_0)'+\pa_z\sigma_{1,\al}^*+h_{1,\al}d_0(\rho_0\eta)'+e_{1,\al}\rho_0',\label{pazsigma1alg}\\
&L_{1,0}=\frac{d_1\kappa_0^2}{d_0^2}c_1^*+L_{1,0}^*,\label{l10ga}
\end{align}
one has
\begin{align*}
    &\langle \rho_0\eta',\sigma_{1,\al}'\rangle=-\frac{\kappa_\al d_1}{d_0}\langle\rho_0\eta',(\eta\rho_0)'\rangle+\langle\rho_0\eta',\pa_z\sigma_{1,\al}^*\rangle+h_{1,\al}d_0\langle \rho_0\eta',(\rho_0\eta)' \rangle+e_{1,\al}\langle\rho_0\eta',\rho_0'\rangle,\\
    &\langle\rho_0'\eta',\sigma_{1,\al}\rangle=-\frac{\kappa_\al d_1}{d_0}\langle\rho_0'\eta',\eta\rho_0\rangle+\langle\rho_0'\eta',\sigma_{1,\al}^*\rangle+h_{1,\al}d_0\langle \rho_0'\eta',\rho_0\eta\rangle+e_{1,\al}\langle\rho_0'\eta',\rho_0\rangle.
\end{align*}
Using \eqref{kappa0kappaal}, we have
\begin{align*}
    -I_2=&\frac{2}{d_0}\langle\pa_z u_0,\pa_zu_1\rangle\\
    =&-\Big(\frac{d_1\kappa_0^2}{d_0^2}c_1^*+L_{1,0}^*\Big)\langle z\rho_0',\eta'\rangle+2d_1\sum_\al \frac{\kappa_\al^2}{d_0^2}\langle \rho_0\eta',\rho_0\eta'\rangle\\
    &-\frac{1}{d_0}\langle z\rho_0',(\pa_t d_0-\Delta d_0)\rho_0'-2\rho_0\nabla d_0\cdot\nabla\pa_z\bar{\omega}\cdot\bar{\omega}-(d_1-z)\frac{\kappa_0^2}{d_0}\rho_0(\eta')^2\rangle\\
    &+2\sum_\al \frac{\kappa_\al}{d_0}\langle \eta',\rho_0'\sigma_{1,\al}^*-\rho_0\pa_z \sigma_{1,\al}^* \rangle\\
    =&-d_1 \frac{\kappa_0^2}{d_0^2}\Big(\langle z\rho_0',c_1^*\eta'-\rho_0(\eta')^2\rangle-2\langle \rho_0\eta',\rho_0\eta'\rangle\Big)-L_{1,0}^*\langle z\rho_0',\eta'\rangle\\
    &-\pa_\nu\langle z\rho_0',(\pa_t d_0-\Delta d_0)\rho_0'-2\rho_0\nabla d_0\cdot\nabla\pa_z\bar{\omega}\cdot\bar{\omega}+z\frac{\kappa_0^2}{d_0}\rho_0(\eta')^2\rangle\\
    &+2\sum_\al \frac{\kappa_\al}{d_0}\langle \eta',\rho_0'\sigma_{1,\al}^*-\rho_0\pa_z \sigma_{1,\al}^* \rangle.
\end{align*}

For $u=\rho \bar\omega+\sum_\al\sigma_{\alpha}\xi_\alpha$, by \eqref{d*},
 we have
\begin{align}\label{d*gamma}
   \mathcal{D}_{*}u\cdot\bar\omega|_\G= \sum_\al\frac{\kappa_\al}{d_0}\Big(2\eta'\sigma_{\alpha}'+\eta''\sigma_{\alpha}\Big),\quad
\mathcal{D}_{*}u\cdot\xi_\beta|_\G= - \frac{\kappa_\beta}{d_0}\Big(2\eta'\rho'+\eta''\rho\Big).
\end{align}
For $I_5$, by
\begin{align}
  &g_0=\mathcal{D}_{*}u_0,\notag\\
  &g_1=\mathcal{D}_{*}u_1+L_{1,0}(x,t)\eta'(z)\bar\omega+\sum_\al L_{1,\al}(x,t)\rho_0^{-1}\eta'(z)\xi_{\alpha}+\sum_\al h_{1,\alpha}(x,t)\Big(2\rho_0'\eta'+\rho_0\eta'' \Big)\xi_{\alpha},\label{g0g1}\\
  &\pa_z\bar{\omega}=\pa_z\xi_\al|_\G=0,\notag
\end{align}
one has
\begin{align*}
I_5=&\langle (d_1-z)\mathcal{D}_* u_0,\rho_1'\bar{\omega}+\sum_\al\pa_z\sigma_{1,\al}\xi_\al\rangle+\langle(d_1-z)\mathcal{D}_*u_1,\rho_0'\bar{\omega}\rangle+L_{1,0}\langle (d_1-z)\eta'\bar\omega,\rho_0'\bar{\omega}\rangle\\
=&-\sum_\al\frac{\kappa_\al}{d_0}\langle(d_1-z),(2\rho_0'\eta'+\eta''\rho_0)\pa_z\sigma_{1,\al}\rangle+\sum_\al\frac{\kappa_\al}{d_0}\langle d_1-z,\rho_0'(2\eta'\pa_z\sigma_{1,\al}+\eta''\sigma_{1,\al})\rangle\\
&+L_{1,0}\langle(d_1-z)\eta',\rho_0'\rangle\\
=&\sum_\al\frac{\kappa_\al}{d_0}\langle d_1-z,\eta''(-\rho_0\pa_z\sigma_{1,\al}+\rho_0'\sigma_{1,\al})\rangle+L_{1,0}\langle(d_1-z)\eta',\rho_0'\rangle.
\end{align*}

Using  \eqref{sigma1alg}, \eqref{pazsigma1alg} and \eqref{l10ga},  a direct calculation gives
\begin{align*}
    -\langle (d_1-z)&\eta''\rho_0,\pa_z\sigma_{1,\al}\rangle_\G=-\langle (d_1-z)\eta''\rho_0,-\frac{\kappa_\al d_1}{d_0}(\eta\rho_0)'+\pa_z\sigma_{1,\al}^*+e_{1,\al}\rho_0'\rangle\\
    =&\frac{\kappa_\al }{d_0}d_1^2 \langle \eta''\rho_0,(\eta\rho_0)'\rangle-d_1\langle \eta''\rho_0,\pa_z\sigma_{1,\al}^*+e_{1,\al}\rho_0'\rangle\\
    &-\frac{\kappa_\al}{d_0}d_1\langle  z\eta''\rho_0,(\eta\rho_0)'\rangle+\langle z\eta''\rho_0,\pa_z\sigma_{1,\al}^*+e_{1,\al}\rho_0'\rangle,
\end{align*}
\begin{align*}
    \langle (d_1-z)&\eta''\rho_0',\sigma_{1,\al}\rangle_\G=\langle (d_1-z)\eta''\rho_0',  -\frac{\kappa_\al d_1}{d_0}\eta(z)\rho_0(z)+\sigma^*_{1,\al}(z,x,t)+
 e_{1,\al}(x,t)\rho_0(z) \rangle\\
 =&-\frac{\kappa_\al }{d_0}d_1^2\langle \eta'' \rho_0',\eta\rho_0 \rangle+d_1\langle \eta''\rho_0',  \sigma^*_{1,\al}(z,x,t)+
 e_{1,\al}(x,t)\rho_0(z) \rangle\\
 &+d_1\frac{\kappa_\al}{d_0}\langle z\eta''\rho_0',  \eta(z)\rho_0(z)\rangle-\langle z\eta''\rho_0',  \sigma^*_{1,\al}(z,x,t)+
 e_{1,\al}(x,t)\rho_0(z) \rangle,
\end{align*}
and
\begin{align*}
L_{1,0}\langle(d_1-z)\eta',\rho_0'\rangle
=&d_1^2\frac{\kappa_0^2}{d_0^2}c_1^*\langle \eta',\rho_0'\rangle-d_1\frac{\kappa_0^2}{d_0^2}c_1^*\langle z\eta',\rho_0'\rangle+d_1 L_{1,0}^*\langle\eta',\rho_0'\rangle-L_{1,0}^*\langle z\eta',\rho_0'\rangle.
\end{align*}
Combining the above  four equalities, we have
\begin{align*}
I_5=&\sum_\al\frac{\kappa_\al}{d_0}\langle d_1-z,\eta''(-\rho_0\pa_z\sigma_{1,\al}+\rho_0'\sigma_{1,\al})\rangle+L_{1,0}\langle(d_1-z)\eta',\rho_0'\rangle\\
=&d_1^2\Big(\sum_\al\frac{\kappa_\al^2}{d_0^2}\langle \eta''\eta',(\rho_0)^2 \rangle +\frac{\kappa_0^2}{d_0^2}c_1^*\langle \eta',\rho_0'\rangle\Big)\\
&+d_1\Big(\sum_\al\frac{\kappa_\al}{d_0}\langle \eta'',\rho_0'\sigma_{1,\al}^*-\rho_0\pa_z\sigma_{1,\al}^*\rangle-\sum_\al\frac{\kappa_\al^2}{d_0^2}\langle z\eta''\rho_0,\eta'\rho_0\rangle-\frac{\kappa_0^2}{d_0^2}c_1^*\langle z\eta',\rho_0'\rangle+L_{1,0}^*\langle \eta',\rho_0'\rangle\Big)\\
&+\sum_\al\frac{\kappa_\al}{d_0}\langle z\eta'',\rho_0\pa_z\sigma_{1,\al}^* -\rho_0'\sigma_{1,\al}^*\rangle-L_{1,0}^*\langle z\eta',\rho_0'\rangle.
\end{align*}

In conclusion, for \eqref{innerrho2} on $\G$, by compatibility condition \eqref{solvablecondition}, one has
\begin{align*}
    \langle \mathcal{D}_2,\pa_z u_0\rangle|_\G=0.
\end{align*}
We have by \eqref{de2} \eqref{fAfB} \eqref{sigma1eq}, $\rho_1^\pm=0$ and $f_B(\rho_0(\pm\infty))=0$ that
\begin{align*}
   \langle\mathcal{D}_1,\pa_zu_1\rangle-I_4 =&\langle\mathcal{D}_1,\rho_1'\bar{\omega}\rangle+\langle \mathcal{D}_1, \sigma_{1,\al}'\xi_\al\rangle- \langle f^{(1)}(u_0,u_1),\rho_0'\bar{\omega}\rangle\\
    =&\langle\mathcal{D}_1,\rho_1'\bar{\omega}\rangle+\langle \pa_z^2\sigma_{1,\al}-f_B(\rho_0)\sigma_{1,\al}, \sigma_{1,\al}'\rangle\\
    &-\langle \s_{1,\al}^2[(2\rho_0^2-a^2-b^2)^2+2(\rho_0^2-a^2)(\rho_0^2-b^2)]\rho_0,\rho_0'\rangle\\
    =&\langle\mathcal{D}_1,\rho_1'\bar{\omega}\rangle+\f12\langle \pa_zf_B(\rho_0),\sigma_{1,\al}^2\rangle-\f12\langle \pa_zf_B(\rho_0),\sigma_{1,\al}^2\rangle\\
    =&\langle\mathcal{D}_1,\rho_1'\bar{\omega}\rangle=0.
\end{align*}
Then the equation of $d_1$ is determined by
\begin{align*}
   \langle \mathcal{D}_2,\pa_zu_0\rangle= I_1+I_2+I_3+I_5=0.
\end{align*}
By \eqref{kappa0kappaal}, then the coefficient of $d_1^2$ is
\begin{align*}
   \sum_\al\frac{\kappa_\al^2}{d_0^2}&\langle \eta''\eta',(\rho_0)^2 \rangle+\frac{\kappa_0^2}{d_0^2}c_1^*\langle \eta',\rho_0'\rangle\\
   =&\frac{\kappa_0^2}{d_0^2}\Big(\langle \eta''\eta',(\rho_0)^2 \rangle+\frac{\langle \rho_0\rho_0',(\eta')^2\rangle}{\langle \rho_0',\eta'\rangle}\langle \eta',\rho_0'\rangle\Big)\\
   =&\frac{\kappa_0^2}{d_0^2}\Big(\langle \eta''\eta',(\rho_0)^2 \rangle+\langle \rho_0\rho_0',(\eta')^2\rangle\Big)\\
   =&0.
\end{align*}
Combining $I_i$ for $i=1,2,3,5$, we obtain \eqref{d1evolution}.
\end{proof}
\subsection{The proof of Lemma \ref{lem:boundary}}\label{A.2}
\begin{proof}
By the matching condition \eqref{matching-condition}, one has
$$\lim_{z\to \pm\infty}\sigma_{1, \beta}(z, x, t)=\sigma_{1, \beta}^\pm(x, t). $$
Recall \eqref{sigma1alg},
we have
\begin{align}
    \sigma_{1,\beta}^+&=-\frac{b\kappa_\beta d_1}{d_0}+\sigma_{1,\beta}^{*+}+bh_{1,\beta}d_0+be_{1,\beta},\notag\\
    \sigma_{1,\beta}^-&=\sigma_{1,\beta}^{*-}+ae_{1,\beta},\label{sigmapm}
\end{align}
which yields the relation \eqref{e1beta1}.

By  \eqref{d0}, one can directly verify that
\begin{align*}
\langle \mathcal{D}_2\cdot\xi_\beta,\rho_0\rangle|_\G=&\langle \partial_zu_0(\partial_td_{1}-\Delta d_{1})\cdot\xi_\beta,\rho_0\rangle-2\langle \nabla d_{0} \cdot \nabla) \partial_zu_{1}\cdot\xi_\beta, \rho_0\rangle\\&+\langle f^{(1)}(u_0,u_{1})\cdot\xi_\beta,\rho_0\rangle +d_2\langle g_0\cdot\xi_\beta,\rho_0\rangle\\
&+\langle (d_1-z)g_{1}\cdot\xi_\beta,\rho_0\rangle+\big(-2\langle(\nabla d_{1}\cdot\nabla )\partial_zu_0\cdot\xi_\beta,\rho_0\rangle+\langle(\partial_{t}u_{0}-\Delta u_{0})\cdot\xi_\beta,\rho_0\rangle\big)\\
:=&\sum_{i=1}^6I_i
\end{align*}
By $\pa_z u_0\cdot \xi_\beta|_\G=\rho_0'\bar{\omega}\cdot\xi_\beta+\rho_0\pa_z\bar{\omega}\cdot \xi_\beta=0$, we have $I_1=0.$ Using \eqref{rho1gamma}, one gets $u_0\cdot u_1|_\G=\rho_0\rho_1=0$ and $u_0\cdot\xi_\beta=0$. Thus, we have from \eqref{de2} that
$f^{(1)}(u_0, u_1)\cdot\xi_\beta=0$ which implies $I_3=0$. For $I_4$, it yields from \eqref{g0xibeta} that
\begin{align}\label{intrho0g0xibeta}
I_4&=d_2\langle g_0\cdot\xi_\beta,\rho_0\rangle=-d_2\langle\frac{\kappa_\beta}{d_0\rho_0}\pa_z(\rho_0^2\eta'\big),\rho_0\rangle\notag=-d_2\frac{\kappa_\beta}{d_0}\langle\pa_z(\rho_0^2\eta'\big),1\rangle =0.
\end{align}
For the term $I_2$, direct calculations give us
\begin{align}
&(\nabla d_0\cdot \nabla \pa_z)u_1 \cdot\xi_\beta|_\G\notag\\
 &=\partial_\nu\pa_z(\rho_1\bar{\omega})\cdot\xi_\beta+\pa_{z\nu}\sigma_{1, \beta}+\sum_{\al\neq \beta}\pa_{z}\sigma_{1,\al}\pa_\nu \xi_\al \cdot\xi_\beta+\sigma_{1, \al}\pa_{z\nu}\xi_\al \cdot\xi_\beta.
\end{align}
As $\partial_z\bar\omega=0,\ \rho_1=0$ on $\Gamma$, one has
\begin{align}
& \partial_\nu \pa_z(\rho_1\bar{\omega})  \cdot\xi_\beta=(\partial_\nu\partial_z\rho_1\bar\omega+\partial_\nu\rho_1\partial_z\bar\omega+\partial_z\rho_1\partial_\nu\bar\omega+\rho_1\partial_\nu\partial_z\bar\omega)\cdot\xi_\beta=0.
\end{align}
We deduce from \eqref{sigma1beta} and  \eqref{sigmapm} that
\begin{align}
\pa_\nu\sigma_{1,\beta}&=-\pa_\nu\Big(\frac{\kappa_\beta d_1}{d_0}\Big)\eta\rho_0+\pa_\nu\sigma_{1,\beta}^*+h_{1,\beta}\rho_0\eta+\pa_\nu e_{1,\beta}\rho_0 \notag\\
    \pa_\nu\sigma_{1,\beta}^+&=-\pa_\nu\Big(\frac{b\kappa_\beta d_1}{d_0}\Big)+\pa_\nu\sigma_{1,\beta}^{*+}+bh_{1,\beta}+b\pa_\nu e_{1,\beta},\notag\\
    \pa_\nu\sigma_{1,\beta}^-&=\pa_\nu\sigma_{1,\beta}^{*-}+a\pa_\nu e_{1,\beta},\label{panuspm}
\end{align}
Then we have
 \begin{align*}
\langle \pa_{z\nu}\sigma_{1, \beta}, \rho_0\rangle_\G=&\pa_\nu\sigma_{1, \beta}\rho_0|_{-\infty}^{+\infty}-\pa_\nu\langle\sigma_{1, \beta}, \rho_0'\rangle\nonumber\\
=&b\pa_\nu\sigma_{1, \beta}^+-a\pa_\nu \sigma_{1, \beta}^--\pa_\nu\langle\sigma_{1, \beta}, \rho_0'\rangle\nonumber\\
=&b\pa_\nu\sigma_{1, \beta}^+-a\pa_\nu \sigma_{1, \beta}^-+\pa_\nu\Big(d_1\frac{\kappa_\beta }{d_0}\Big)\langle \eta\rho_0,\rho_0'\rangle\\
&-\pa_\nu\langle\sigma_{1, \beta}^*, \rho_0'\rangle- h_{1,\beta}\langle\rho_0\eta,\rho_0'\rangle-\frac{b^2-a^2}{2}\pa_\nu e_{1\beta}\\
=&b\pa_\nu\sigma_{1, \beta}^+-a\pa_\nu \sigma_{1, \beta}^--\frac{d_1}{2}\pa_\nu\Big(\frac{\kappa_\beta }{d_0}\Big)\langle \eta',\rho_0^2\rangle+\frac{b^2}{2}\Big(d_1\pa_\nu\Big(\frac{\kappa_\beta}{d_0}\Big)-h_{1,\beta}\Big)\\
&-\pa_\nu\langle\sigma_{1, \beta}^*, \rho_0'\rangle+\frac{h_{1,\beta}}{2}\langle \eta',\rho_0^2\rangle-\frac{b^2-a^2}{2}\pa_\nu e_{1\beta},
 \end{align*}
 and for $\al\neq \beta$ that
\begin{align}
\langle\pa_z&\sigma_{1, \al}\pa_\nu \xi_\al \cdot\xi_\beta +\sigma_{1, \al}\pa_{z\nu}\xi_\al \cdot\xi_\beta,\rho_0 \rangle_\G\nonumber\\
=&\sigma_{1, \al}\pa_\nu\xi_\al\cdot\xi_\beta\rho_0|_{-\infty}^{+\infty}-\langle \sigma_{1, \al}\pa_\nu\xi_\al\cdot\xi_\beta,\rho_0'\rangle\nonumber\\
=&b\sigma_{1, \al}^+\pa_\nu\xi_\al^+\cdot\xi_\beta^+-a\sigma_{1, \al}^-\pa_\nu\xi_\al^-\cdot\xi_\beta^--\langle\sigma_{1, \al}^*\pa_\nu\xi_\al\cdot\xi_\beta,\rho_0'\rangle\notag\\
&+d_1\frac{\kappa_\al}{d_0}\langle \eta\rho_0\pa_\nu\xi_\al\cdot \xi_\beta,\rho_0'\rangle-{e_{1, \al}(x,t)}\langle\rho_0\pa_\nu\xi_\al\cdot\xi_\beta,\rho_0'\rangle.
\end{align}

For $I_5,$  we yield from
\begin{align*}
  &g_1=\mathcal{D}_{*}u_1+L_{1,0}(x,t)\eta'(z)\bar\omega+\sum_\al L_{1,\al}(x,t)\rho_0^{-1}\eta'(z)\xi_{\alpha}+\sum_\al h_{1,\alpha}(x,t)\Big(2\rho_0'\eta'+\rho_0\eta'' \Big)\xi_{\alpha},\\
&\mathcal{D}_{*}u_i\cdot\xi_\beta|_\G=- \frac{\kappa_\beta}{d_0}\Big(2\eta'\rho_i'+\eta''\rho_i\Big),\quad \rho_1|_\G=0,
\end{align*}
 that
\begin{align*}
I_5=&\langle(d_1-z)g_1\cdot\xi_\beta, \rho_0\rangle\\
=&L_{1, \beta}\langle d_1-z,\eta'(z)\rangle+h_{1,\beta}\langle (d_1-z)(2\rho_0'\eta'+\rho_0\eta''),\rho_0\rangle+\langle (d_1-z)\mathcal{D}_*u_1\cdot\xi_\beta,\rho_0\rangle\\
=&L_{1, \beta}\langle d_1-z,\eta'(z)\rangle+h_{1,\beta}\langle \rho_0^2,\eta' \rangle.
\end{align*}
We have by \eqref{sigmapm} that
\begin{align*}
  b^2h_{1,\beta}+(b^2-a^2)\pa_\nu e_{1, \beta}+b\pa_\nu \sigma_{1,\beta}^{*+}-a\pa_\nu \sigma _{1,\beta}^{*-}-b^2\pa_\nu\Big(\frac{ \kappa_\beta d_1}{d_0}\Big)=b\partial_\nu\sigma^+_{1,\beta}-a\partial_\nu\sigma^-_{1,\beta},
\end{align*}
and
\begin{align*}
    h_{1,\beta}(x,t)=\frac{1}{ab}\Big(a\pa_\nu\sigma_{1,\beta}^+-b\pa_\nu\sigma_{1,\beta}^-+b\pa_\nu\sigma_{1,\beta}^{*-}-a\pa_\nu\sigma_{1,\beta}^{*+}\Big)+\pa_\nu\Big(\frac{ \kappa_\beta d_1}{d_0}\Big).
\end{align*}
Therefore, we have
\begin{align*}
I_2+I_5=&-2b\pa_\nu\sigma_{1, \beta}^++2a\pa_\nu \sigma_{1, \beta}^-+d_1\pa_\nu\Big(\frac{\kappa_\beta }{d_0}\Big)\langle \eta',\rho_0^2\rangle-b^2\Big(d_1\pa_\nu\Big(\frac{\kappa_\beta}{d_0}\Big)-h_{1,\beta}\Big)\\
&+2\pa_\nu\langle\sigma_{1, \beta}^*, \rho_0'\rangle-h_{1,\beta}\langle \eta',\rho_0^2\rangle+{(b^2-a^2)}\pa_\nu e_{1\beta}\\
&-2b\sum_{\al\neq \beta}\sigma_{1, \al}^+\pa_\nu\xi_\al^+\cdot\xi_\beta^++2a\sum_{\al\neq \beta}\sigma_{1, \al}^-\pa_\nu\xi_\al^-\cdot\xi_\beta^-+2\sum_{\al\neq \beta}\langle\sigma_{1, \al}^*\pa_\nu\xi_\al\cdot\xi_\beta,\rho_0'\rangle\notag\\
&-2d_1\sum_{\al\neq\beta}\frac{\kappa_\al}{d_0}\langle \eta\rho_0\pa_\nu\xi_\al\cdot \xi_\beta,\rho_0' \rangle+2\sum_{\al\neq \beta}{e_{1, \al}(x, t)}\langle\pa_\nu\xi_\al\cdot\xi_\beta\rho_0,\rho_0'\rangle\\
&+L_{1, \beta}\langle d_1-z,\eta'(z)\rangle+h_{1,\beta}\langle \rho_0^2,\eta' \rangle\\
=&a\pa_\nu\sigma_{1,\beta}^--b\pa_\nu \sigma_{1,\beta}^++a\pa_\nu\sigma_{1,\beta}^{*-}-b\pa_\nu\sigma_{1, \beta}^{*+}+\sum_{\al\neq \beta}\frac{\kappa_\al d_1 b^2}{d_0}\pa_\nu\xi_\al^+\cdot\xi_\beta^+\\
&+\sum_{\al\neq \beta}-b\sigma_{1,\al}^{*+}\pa_\nu\xi_\al^+\cdot\xi_\beta^++a\sigma_{1,\al}^{*-}\pa_\nu\xi_\al^-\cdot\xi_\beta^-+2\pa_\nu\langle \sigma_{1, \beta}^*, \rho_0'\rangle +d_1\pa_\nu\Big(\frac{\kappa_\beta}{d_0}\Big)\langle \eta',\rho_0^2 \rangle\\
&-b\sum_{\al\neq \beta}\sigma_{1, \al}^+\pa_\nu\xi_\al^+\cdot\xi_\beta^++a\sum_{\al\neq \beta}\sigma_{1, \al}^-\pa_\nu\xi_\al^-\cdot\xi_\beta^-+2\sum_{\al\neq \beta}\langle\sigma_{1, \al}^*\pa_\nu\xi_\al\cdot\xi_\beta,\rho_0'\rangle\notag\\
&-2d_1\sum_{\al\neq\beta}\frac{\kappa_\al}{d_0}\langle \eta\rho_0\pa_\nu\xi_\al\cdot \xi_\beta,\rho_0'\rangle+L_{1, \beta}\langle d_1-z,\eta'(z)\rangle.
\end{align*}
Indeed, by  $\pa_z\xi_\al\cdot\xi_\beta=0$, we have used
\begin{align*}
    2\sum_{\al\neq\beta}e_{1,\al}(x,t)\langle \pa_\nu\xi_\al\cdot \xi_\beta \rho_0,\rho_0' \rangle =&\sum_{\al\neq\beta}e_{1,\al}\Big(\pa_\nu\xi_\al\cdot\xi_\beta \rho_0^2|_{-\infty}^{+\infty}-\langle \pa_{z\nu}\xi_\al \cdot \xi_\beta, \rho_0^2\rangle\Big)\\
    =&\sum_{\al\neq \beta}e_{1,\al}(b^2\pa_\nu\xi_\al^+\cdot\xi_\beta^+-a^2\pa_\nu\xi_\al^-\cdot \xi_\beta^-)-e_{1,\al}\langle \pa_{z\nu}\xi_\al \cdot \xi_\beta, \rho_0^2\rangle\\
    =&\sum_{\al\neq \beta}\Big(b\sigma_{1,\al}^++\frac{\kappa_\al d_1 b^2}{d_0}-b\sigma_{1,\al}^{*+}\Big)\pa_\nu\xi_\al^+\cdot\xi_\beta^+-\Big(a\sigma_{1,\al}^--a\sigma_{1,\al}^{*-}\Big)\pa_\nu\xi_\al^-\cdot\xi_\beta^-.
\end{align*}
Combining $I_i$ for $i=1,\cdots,6$, we complete the proof of \eqref{sigma1betaboundary1}.
\end{proof}
\subsection{The proof of Lemma \ref{lemdk}}\label{A.3}
\begin{proof}
By $\pa_t d_0-\Delta d_0=0$ and $g_0\cdot\bar{\omega}=0$ on $\G$, one has
\begin{align*}
    \langle \mathcal{D}_{k+2},&\pa_z u_0 \rangle_\G+\langle \mathcal{D}_1,\pa_z u_{k+1}\rangle_\G\\
    =&\langle \partial_zu_0(\partial_td_{k+1}-\Delta d_{k+1})+\partial_zu_{k+1}(\partial_td_{0}-\Delta d_{0})-2\nabla d_{k+1}\cdot \nabla \partial_zu_0\nonumber\\
&-2\nabla d_{0}\cdot \nabla \partial_zu_{k+1}+f^{(k+1)}(u_0, \cdots, u_{k+1})+(d_{k+2}g_0+d_{k+1}g_1)\nonumber\\
&+(d_1-z)g_{k+1}+\mathcal{A}_{k},\pa_zu_0\rangle+\langle \partial_zu_0(\partial_td_0-\Delta d_0)-2\nabla d_0\cdot\nabla\partial_z u_0 \\
&+g_0(d_1-z), \pa_zu_{k+1}\rangle\\
=&( \pa_t d_{k+1}-\Delta d_{k+1})\langle \pa_zu_0,\pa_zu_0\rangle-2\langle \nabla d_{k+1}\cdot \nabla\pa_z u_0,\pa_z u_0\rangle-2\pa_\nu\langle\pa_zu_{k+1},\pa_z u_0\rangle\\
&+\big(\langle(d_1-z)g_{k+1}, \pa_z u_0 \rangle+\langle (d_1-z)g_0,\pa_zu_{k+1}\rangle\big)\\
&+\langle f^{(k+1)},\pa_z u_0\rangle  +(\langle \mathcal{A}_k,\pa_z u_0\rangle+d_{k+1}\langle g_1,\pa_z u_0\rangle):=\sum_{i=1}^{6}I_i.
\end{align*}
For $I_1$, one has
$$I_1=e(\pa_t d_{k+1}-\Delta d_{k+1}).$$
Direct calculation gives us
\begin{align*}
I_2=-\nabla d_{k+1}\cdot \nabla\langle \pa_z u_0,\pa_zu_0\rangle=-\nabla d_{k+1}\cdot\nabla (\rho_0'(z))^2=0.
\end{align*}
By \begin{align*}
    \bar{\omega}\cdot\xi_\al=0,\ \pa_z\bar{\om}\cdot\xi_\al=-\bar{\omega}\cdot \pa_z\xi_\al=-\kappa_\al (x,t)\eta',\ \pa_z\bar{\omega}\cdot\pa_z\xi_\al=0,\ |\pa_z\bar{\omega}|^2=\kappa_0^2(x,t)(\eta')^2,
\end{align*}
one has
\begin{align*}
    \langle\pa_zu_0,\pa_zu_{k+1}\rangle=&\langle \pa_z(\rho_0\bar{\omega}),\pa_z(\rho_{k+1}\bar{\omega}+\sum_\al \sigma_{{k+1},\al}\xi_\al)\rangle\\
    =&\langle\rho_0'\bar{\omega}+\rho_0\pa_z\bar{\omega},\rho_{k+1}'\bar{\omega}+\rho_{k+1}\pa_z\bar{\omega}\rangle+\sum_\al\langle\rho_0'\bar{\omega}+\rho_0\pa_z\bar{\omega},\sigma_{{k+1},\al}'\xi_\al +\sigma_{{k+1},\al}\partial_z\xi_\al\rangle\\
    =&\langle \rho_0',\rho_{k+1}'\rangle +\langle \rho_0 \pa_z\bar{\omega},\rho_{k+1}\pa_z\bar{\omega}\rangle +\sum_\al\langle \rho_0\pa_z\bar{\omega},\sigma_{{k+1},\al}'\xi_\al\rangle +\langle \rho_0'\bar{\omega},\sigma_{{k+1},\al}\pa_z\xi_\al\rangle\\
    =&\langle \rho_0',\rho_{k+1}'\rangle +\kappa_0^2\langle |\eta'|^2\rho_0,\rho_{k+1}\rangle +\sum_\al\kappa_\al (-\langle \rho_0\eta',\sigma_{{k+1},\al}'\rangle +\langle \rho_0'\eta',\sigma_{{k+1},\al}\rangle).
\end{align*}
Multiplying  $\pa_z^2\rho_{k+1}-f_A(\rho_0)\rho_{k+1}=L_{k+1,0}\eta'd_0+\mathcal{D}_{k+1}\cdot \bar{\omega}$ by $z\rho_0'$, integrating over $\BR$ and using Lemma \ref{v1product}, we have
\begin{align*}
&\langle\rho_0',\rho_{k+1}'\rangle= -\f12d_0L_{{k+1},0}\langle z\rho_0',\eta'\rangle-\f12\langle z\rho_0',\mathcal{D}_{k+1}\cdot \bar{\omega}\rangle,
\end{align*}
where $\mathcal{D}_{k+1}\cdot\bar
\omega|_\G$ depends only on $\mathcal{V}^i$ for $i\leq k$.
By
\begin{align}
 &\sigma_{{k+1},\al}=-\frac{\kappa_\al d_{k+1}}{d_0}\eta(z)\rho_0(z)+\sigma^*_{{k+1},\al}+h_{k+1,\al}d_0\rho_0\eta+
 e_{{k+1},\al}(x,t)\rho_0(z),\label{sigma{k+1}alg}\\
 &\pa_z\sigma_{{k+1},\al}=-\frac{\kappa_\al d_{k+1}}{d_0}(\eta\rho_0)'+\pa_z\sigma_{{k+1},\al}^*+h_{k+1,\al}d_0(\rho_0\eta)'+e_{{k+1},\al}\rho_0',\label{pazsigma{k+1}alg}\\
&L_{{k+1},0}=\frac{d_{k+1}\kappa_0^2}{d_0^2}c_{1}^*+L_{{k+1},0}^*,\label{l{k+1}0ga}
\end{align}
one has
\begin{align*}
    &\langle \rho_0\eta',\sigma_{{k+1},\al}'\rangle=-\frac{\kappa_\al d_{k+1}}{d_0}\langle\rho_0\eta',(\eta\rho_0)'\rangle+\langle\rho_0\eta',\pa_z\sigma_{{k+1},\al}^*\rangle+h_{k+1,\al}d_0\langle\rho_0\eta',(\rho_0\eta)'\rangle\\
    &\qquad\qquad\qquad\ \ \ +e_{{k+1},\al}\langle\rho_0\eta',\rho_0'\rangle,\\
    &\langle\rho_0'\eta',\sigma_{{k+1},\al}\rangle=-\frac{\kappa_\al d_{k+1}}{d_0}\langle\rho_0'\eta',\eta\rho_0\rangle+\langle\rho_0'\eta',\sigma_{{k+1},\al}^*\rangle+h_{k+1,\al}d_0\langle \rho_0'\eta',\rho_0\eta\rangle\\
    &\qquad\qquad\qquad\ \ \ +e_{{k+1},\al}\langle\rho_0'\eta',\rho_0\rangle.
\end{align*}
Then we have
\begin{align*}
    -I_3=&2\pa_\nu\langle\pa_z u_0,\pa_zu_{k+1}\rangle\\
    =&-\Big(\frac{d_{k+1}\kappa_0^2}{d_0^2}c_{1}^*+L_{{k+1},0}^*\Big)\langle z\rho_0',\eta'\rangle-\pa_\nu\langle z\rho_0',\mathcal{D}_{k+1}\cdot \bar{\omega}\rangle\\
    &+2d_{k+1}\sum_\al\frac{\kappa_\al ^2}{d_0^2}\langle \rho_0\eta',\rho_0\eta'\rangle+2\sum_\al \frac{\kappa_\al}{d_0}\langle \eta',\rho_0'\sigma_{k+1,\al}^*-\rho_0\pa_z \sigma_{k+1,\al}^*\rangle\\
    =&-d_{k+1}\frac{\kappa_0^2}{d_0^2}\Big(\langle z\rho_0',c_1^*\eta'-\rho_0(\eta')^2\rangle-2\langle \rho_0\eta',\rho_0\eta'\rangle\Big)-L_{{k+1},0}^* \langle z\rho_0',\eta'\rangle\\
    &-\pa_\nu\langle z\rho_0',\mathcal{D}_{k+1}^r\cdot \bar{\omega}\rangle+2\sum_\al \frac{\kappa_\al}{d_0}\langle \eta',\rho_0'\sigma_{k+1,\al}^*-\rho_0\pa_z \sigma_{k+1,\al}^*\rangle.
\end{align*}
Recall $\mathcal{D}_{k+1}=\mathcal{D}_{k+1}^r+d_{k+1}g_0$.
 Then by \eqref{g0omega}, we have
\begin{align*}
    \pa_\nu\langle z\rho_0',\mathcal{D}_{k+1}\cdot\bar{\omega}\rangle&=\pa_\nu\langle z\rho_0',d_{k+1}g_0\cdot\bar{\omega}\rangle+\pa_\nu\langle z\rho_0',\mathcal{D}_{k+1}^r\cdot\bar{\omega}\rangle\\
    &=-\pa_\nu\Big(\frac{d_{k+1}\kappa_0^2}{d_0}\Big)\langle z\rho_0',\rho_0(\eta')^2\rangle+\pa_\nu\langle z\rho_0',\mathcal{D}_{k+1}^r\cdot\bar{\omega}\rangle.
\end{align*}
By \eqref{d*gamma}, \eqref{g0g1} and
\begin{align*}
 & g_{k+1}=\mathcal{D}_{*}u_{k+1}+L_{k+1,0}(x,t)\eta'\bar\omega+\sum_\al L_{k+1,\al}(x,t)\rho_0^{-1}\eta'\xi_{\alpha}+\sum_\al h_{k+1,\alpha}(x,t)\Big(2\rho_0'\eta'+\rho_0\eta'' \Big)\xi_{\alpha},
\end{align*}
one has
\begin{align*}
    I_4=&\langle(d_1-z)g_{k+1}, \pa_z u_0 \rangle+\langle (d_1-z)g_0,\pa_zu_{k+1}\rangle\\
    =&\langle(d_1-z)g_{k+1}, \rho_0'\bar{\om} \rangle+\langle (d_1-z)g_0,\pa_z\sigma_{k+1,\al}\xi_\al\rangle\\
    =&\langle(d_1-z)\mathcal{D}_*u_{k+1}, \rho_0'\bar{\om} \rangle+L_{k+1,0}\langle(d_1-z)\eta', \rho_0'\rangle+\langle (d_1-z)\mathcal{D}_*u_0,\pa_z\sigma_{k+1,\al}\xi_\al\rangle\\
    =&(\langle(d_1-z)\mathcal{D}_*u_{k+1}, \rho_0'\bar{\om} \rangle+\langle (d_1-z)\mathcal{D}_*u_0,\pa_z\sigma_{k+1,\al}\xi_\al\rangle)+L_{k+1,0}\langle(d_1-z)\eta', \rho_0'\rangle\\
    =&\sum_\al\frac{\kappa_\al}{d_0}\big(\langle(d_1-z)(2\eta'\sigma_{k+1,\al}'+\eta''\sigma_{k+1,\al}), \rho_0' \rangle-\langle (d_1-z)(2\eta'\rho_0'+\eta''\rho_0),\pa_z\sigma_{k+1,\al}\rangle\big)\\
    &+L_{k+1,0}\langle(d_1-z)\eta', \rho_0'\rangle\\
    =&\sum_\al \frac{\kappa_\al}{d_0}\langle d_1-z, \eta''(-\rho_0\pa_z\sigma_{k+1,\al}+\rho_0'\sigma_{k+1,\al})+L_{k+1,0}\langle(d_1-z)\eta', \rho_0'\rangle.
\end{align*}
By \eqref{sigma{k+1}alg}-\eqref{l{k+1}0ga}, one has
\begin{align*}
    -\langle (d_1-z)&\eta''\rho_0,\pa_z\sigma_{{k+1},\al}\rangle_\G=-\langle (d_1-z)\eta''\rho_0,-\frac{\kappa_\al d_{k+1}}{d_0}(\eta\rho_0)'+\pa_z\sigma_{{k+1},\al}^*+e_{{k+1},\al}\rho_0'\rangle\\
    =&\frac{\kappa_\al }{d_0}d_{k+1} \langle d_1-z,\eta''\rho_0(\eta\rho_0)'\rangle-\langle (d_1-z)\eta''\rho_0,\pa_z\sigma_{{k+1},\al}^*+e_{{k+1},\al}\rho_0'\rangle,
\end{align*}
\begin{align*}
    \langle (d_1-z)&\eta''\rho_0',\sigma_{{k+1},\al}\rangle_\G=\langle (d_1-z)\eta''\rho_0',  -\frac{\kappa_\al d_{k+1}}{d_0}\eta\rho_0+\sigma^*_{{k+1},\al}+
 e_{{k+1},\al}(x,t)\rho_0 \rangle\\
 =&-\frac{\kappa_\al }{d_0}d_{k+1}\langle (d_1-z), \eta'' \rho_0'\eta\rho_0 \rangle+\langle (d_1-z)\eta''\rho_0',  \sigma^*_{{k+1},\al}+
 e_{{k+1},\al}\rho_0 \rangle,
\end{align*}
and
\begin{align*}
L_{k+1,0}\langle(d_1-z)\eta', \rho_0'\rangle=&d_{k+1}\frac{\kappa_0^2}{d_0^2}c_1^*\langle (d_1-z)\eta',\rho_0'\rangle+L_{k+1,0}^*\langle(d_1-z)\eta', \rho_0'\rangle.
\end{align*}
Combining the above  four equalities, we have
\begin{align*}
I_4=&\sum_\al\frac{\kappa_\al}{d_0}\langle d_1-z,\eta''(-\rho_0\pa_z\sigma_{{k+1},\al}+\rho_0'\sigma_{{k+1},\al})\rangle+L_{{k+1},0}\langle(d_1-z)\eta',\rho_0'\rangle\\
=&\frac{\kappa_0^2}{d_0^2}d_{k+1}\langle d_1-z, \eta''\eta'\rho_0^2+c_1^*\eta'\rho_0'\rangle+\sum_\al\pa_\nu\kappa_\al\langle (d_1-z)\eta'',\rho_0'\sigma_{k+1,\al}^*-\rho_0\pa_z\sigma_{k+1,\al}^*\rangle\\
&+L_{k+1,0}^*\langle(d_1-z)\eta', \rho_0'\rangle.
\end{align*}
Using  \eqref{fk-1} and $\rho_1|_\G=0$, one has
 \begin{align*}
     f^{(k+1)}\cdot \bar{\omega}|_\G=2\sum_\al \sigma_{1,\al}\sigma_{k+1,\al}[(2\rho_0^2-a^2-b^2)^2+2(\rho_0^2-a^2)(\rho_0^2-b^2)]\rho_0+f_4^{(k+1)}\cdot \bar{\omega}.
 \end{align*}
We have by  \eqref{fAfB}, $\rho_1^\pm=0$, $f_B(\rho_0(\pm\infty))=0$, $\mathcal{D}_1\cdot\bar{\omega}|_\G= 0$ and
\begin{align*}
    (\pa_z^2-f_B(\rho_0))\sigma_{1,\al}=\mathcal{D}_{1}\cdot\xi_\al,\ \ (\pa_z^2-f_B(\rho_0))\sigma_{k+1,\al}=\mathcal{D}_{k+1}\cdot\xi_\al,\text{ on }\G,
\end{align*}
that
\begin{align*}
  \langle \mathcal{D}_1,\pa_z u_{k+1}\rangle_\G-I_5=&\langle\mathcal{D}_1,\rho_{k+1}'\bar{\omega}\rangle+\sum_\al \langle \mathcal{D}_1, \sigma_{k+1,\al}'\xi_\al\rangle- \langle f^{(k+1)},\rho_0'\bar{\omega}\rangle\\
    =&\langle\mathcal{D}_1,\rho_{k+1}'\bar{\omega}\rangle+\sum_\al\langle \pa_z^2\sigma_{1,\al}-f_B(\rho_0)\sigma_{1,\al}, \sigma_{k+1,\al}'\rangle\\
    &-2 \sum_\al\langle \s_{1,\al}\s_{k+1,\al}[(2\rho_0^2-a^2-b^2)^2+2(\rho_0^2-a^2)(\rho_0^2-b^2)]\rho_0,\rho_0'\rangle\\
    &-\langle f_4^{(k+1)},\rho_0'\bar{\omega}\rangle\\
    =&\langle\mathcal{D}_1,\rho_{k+1}'\bar{\omega}\rangle+\sum_\al\langle f_B(\rho_0)\s_{k+1,\al},\s_{1,\al}'\rangle-\langle \pa_z^2\s_{k+1,\al},\s_{1,\al}'\rangle\\
    &+\langle f_B'(\rho_0),\sigma_{1,\al}\s_{k+1,\al}\rangle-\langle f_B'(\rho_0),\sigma_{1,\al}\sigma_{k+1,\al}\rangle-\langle f_4^{(k+1)},\rho_0'\bar{\omega}\rangle\\
    =&-\sum_\al\langle \mathcal{D}_{k+1},\sigma_{1,\al}'\xi_\al\rangle -\langle f_4^{(k+1)},\rho_0'\bar{\omega}\rangle\\
    =&-\sum_\al\big(\langle  \mathcal{D}_{k+1}^r,\sigma_{1,\al}'\xi_\al\rangle+d_{k+1}\langle g_0\cdot\xi_\al, \sigma_{1,\al}'\rangle\big)-\langle f_4^{(k+1)},\rho_0'\bar{\omega}\rangle
\end{align*}
Together with \eqref{dk+2rho0'}, we obtain
\begin{align*}
    \langle \mathcal{D}_{k+2},\pa_z u_0\rangle_\G =&\sum_{i=1}^4I_i+I_6+\sum_\al\big(\langle  \mathcal{D}_{k+1}^r,\sigma_{1,\al}'\xi_\al\rangle+d_{k+1}\langle g_0\cdot\xi_\al, \sigma_{1,\al}'\rangle\big) +\langle f_4^{(k+1)},\rho_0'\bar{\omega}\rangle\\
    =&e(\pa_t d_{k+1}-\Delta d_{k+1})+d_{k+1}\Big( \langle g_1,\pa_zu_0\rangle+\sum_\al\langle g_0, \sigma_{1,\al}'\xi_\al\rangle\\
    &+\frac{\kappa_0^2}{d_0^2}\big(\langle d_1-z,\eta'\eta''(\rho_0)^2+c_1^*\eta'\rho_0'\rangle
    +\langle z\rho_0',c_1^*\eta'-\rho_0(\eta')^2\rangle-2\langle \rho_0\eta',\rho_0\eta'\rangle\big)\Big)\\
    &+L_{{k+1},0}^* \langle z\rho_0',\eta'\rangle+\pa_\nu\langle z\rho_0',\mathcal{D}_{k+1}^r\cdot \bar{\omega}\rangle-2\sum_\al \frac{\kappa_\al}{d_0}\langle \eta',\rho_0'\sigma_{k+1,\al}^*-\rho_0\pa_z \sigma_{k+1,\al}^* \rangle\\
    &+\langle \mathcal{A}_k,\pa_z u_0\rangle+L_{k+1,0}^*\langle(d_1-z)\eta', \rho_0'\rangle+\sum_\al\frac{\kappa_\al}{d_0}\langle (d_1-z)\eta'',\rho_0'\sigma_{k+1,\al}^*-\rho_0\pa_z\sigma_{k+1,\al}^*\rangle\\
    &+\sum_\al\langle \mathcal{D}_{k+1}^r\cdot\xi_\al,\sigma_{1,\al}'\rangle +\langle f_4^{(k+1)}\cdot \bar{\omega},\rho_0'\rangle.
\end{align*}
\end{proof}
\subsection{The proof of Lemma \ref{lem:sigmakbetaboundary}}\label{A.4}

\begin{proof}
By $$\lim_{z\to \pm\infty}\sigma_{{k+1},\beta}(z,x,t)=\sigma_{{k+1},\beta}^\pm(x,t),$$
and
\begin{align*}
    \sigma_{{k+1},\beta}(z,x,t)|_{\G}=-\frac{\kappa_\beta d_{k+1}}{d_0}\eta(z)\rho_0(z)+\sigma^*_{{k+1},\beta}(z,x,t)+
 e_{{k+1},\beta}(x,t)\rho_0(z),
\end{align*}
we have
\begin{align}
    &\sigma_{{k+1},\beta}^+=-\frac{b\kappa_\beta d_{k+1}}{d_0}+\sigma_{k+1,\beta}^{*+}+bh_{k+1,\beta}d_0+be_{k+1,\beta},\notag\\
    &\sigma_{{k+1},\beta}^-=\sigma_{k+1,\beta}^{*-}+ae_{k+1,\beta}.\label{sigmak+1pm}
\end{align}
It is direct to obtain \eqref{sigmakbetaboundary2}.
By $\pa_t d_0-\Delta d_0=0$ and $\pa_z u_0\cdot\xi_\beta=0$ on $\G$, we can write
\begin{align*}
\langle\mathcal{D}_{k+2}\cdot\xi_\beta,\rho_0\rangle_\G=&-2\langle(\nabla d_{k+1}\cdot \nabla) \partial_zu_0\cdot\xi_\beta,\rho_0\rangle-2\langle(\nabla d_{0}\cdot\nabla) \partial_zu_{k+1}\cdot\xi_\beta,\rho_0\rangle\\
&+\langle f^{(k+1)}(u_0, \cdots, u_{k+1})\cdot\xi_\beta,\rho_0\rangle+\langle (d_{k+2}g_0+d_{k+1}g_1)\cdot\xi_\beta,\rho_0\rangle\\
&+\langle(d_1-z)g_{k+1}\cdot\xi_\beta,\rho_0\rangle+\langle \mathcal{A}_{k}\cdot\xi_\beta,\rho_0\rangle\\
:=&\sum_{i=1}^6I_i.
\end{align*}
We have by \eqref{g0xibeta} that $I_4=d_{k+1}\langle g_1\cdot\xi_\beta,\rho_0\rangle$. For $I_2$, one has
 \begin{align*}
 (\nabla &d_0\cdot\nabla\pa_z)u_{k+1}\cdot \xi_\beta|_\G\\
 &=\pa_z\rho_{k+1}\pa_\nu\bar{\omega}\cdot\xi_\beta+\rho_{k+1}\pa_{z\nu}\bar{\omega}\cdot \xi_\beta +\pa_{z\nu}\sigma_{{k+1}, \beta}+\sum_{\al\neq \beta}\pa_{z}\sigma_{{k+1},\al}\pa_\nu \xi_\al \cdot\xi_\beta +\sigma_{{k+1},\al}\pa_{z\nu}\xi_\al \cdot\xi_\beta.
 \end{align*}
 We deduce from \eqref{sigmakbeta} and
\begin{align*}
\pa_\nu \sigma_{k+1,\beta}&=-\pa_\nu\Big(\frac{\kappa_\beta d_{k+1}}{d_0}\Big)\eta\rho_0+\pa_\nu  \sigma_{k+1,\beta}^{*}+h_{k+1,\beta}\rho_0\eta+\pa_\nu e_{k+1,\beta}\rho_0, \\
    \pa_\nu\sigma_{k+1,\beta}^+&=-b\pa_\nu\Big(\frac{\kappa_\beta d_{k+1}}{d_0}\Big)+\pa_\nu\sigma_{k+1,\beta}^{*+}+bh_{k+1,\beta}+b\pa_\nu e_{k+1,\beta},\notag\\
    \pa_\nu\sigma_{k+1,\beta}^-&=\pa_\nu\sigma_{k+1,\beta}^{*-}+a\pa_\nu e_{k+1,\beta},
\end{align*}
that
 \begin{align*}
\langle\pa_{z\nu}\sigma_{{k+1}, \beta}, \rho_0\rangle_\G=&\pa_\nu\sigma_{{k+1}, \beta}\rho_0|_{-\infty}^{+\infty}-\pa_\nu\langle\sigma_{{k+1}, \beta}, \rho_0'\rangle\nonumber\\
=&b\pa_\nu\sigma_{{k+1}, \beta}^+-a\pa_\nu \sigma_{{k+1}, \beta}^--\pa_\nu\langle\sigma_{{k+1}, \beta}, \rho_0'\rangle\nonumber\\
=&b\pa_\nu\sigma_{{k+1}, \beta}^+-a\pa_\nu \sigma_{{k+1}, \beta}^-+\pa_\nu \Big(\frac{d_{k+1}\kappa_\beta }{d_0}\Big)\langle\langle \eta\rho_0,\rho_0'\rangle\\
&-\pa_\nu\langle\sigma_{{k+1}, \beta}^*, \rho_0'\rangle- h_{{k+1},\beta}\langle\rho_0\eta,\rho_0'\rangle-\frac{b^2-a^2}{2}\pa_\nu e_{{k+1}\beta}\\
=&b\pa_\nu\sigma_{{k+1}, \beta}^+-a\pa_\nu \sigma_{{k+1}, \beta}^--\pa_\nu\Big(d_{k+1}\frac{\kappa_\beta }{2d_0}\Big)\langle \eta',\rho_0^2\rangle+\frac{b^2}{2}\Big(\pa_\nu\Big(\frac{d_{k+1}\kappa_\beta}{d_0}\Big)-h_{k+1,\beta}\Big)\\
&-\pa_\nu\langle\sigma_{{k+1}, \beta}^*, \rho_0'\rangle+\frac{h_{{k+1},\beta}}{2}\langle \eta',\rho_0^2\rangle-\frac{b^2-a^2}{2}\pa_\nu e_{{k+1}\beta},
 \end{align*}
 and for $\al\neq \beta$ that
\begin{align*}
\langle\pa_z&\sigma_{{k+1}, \al}\pa_\nu \xi_\al \cdot\xi_\beta +\sigma_{{k+1}, \al}\pa_{z\nu}\xi_\al \cdot\xi_\beta,\rho_0 \rangle_\G\nonumber\\
=&\sigma_{{k+1}, \al}\pa_\nu\xi_\al\cdot\xi_\beta\rho_0|_{-\infty}^{+\infty}-\langle\sigma_{{k+1}, \al}\pa_\nu\xi_\al\cdot\xi_\beta,\rho_0'\rangle\nonumber\\
=&b\sigma_{{k+1}, \al}^+\pa_\nu\xi_\al^+\cdot\xi_\beta^+-a\sigma_{{k+1}, \al}^-\pa_\nu\xi_\al^-\cdot\xi_\beta^--\langle\sigma_{{k+1}, \al}^*\pa_\nu\xi_\al\cdot\xi_\beta,\rho_0'\rangle\notag\\
&+d_{k+1}\frac{\kappa_\al}{d_0}\langle\eta\rho_0\pa_\nu\xi_\al\cdot \xi_\beta,\rho_0'\rangle-{e_{{k+1}, \al}(x,t)}\langle\pa_\nu\xi_\al\cdot\xi_\beta\rho_0,\rho_0'\rangle.
\end{align*}
For $I_3$, we have by \eqref{de2}, \eqref{fk-1} and \eqref{rho1gamma} that
\begin{align}
I_3=&\langle f^{(k+1)}(u_0,\cdots,u_{k+1})\cdot \xi_\beta, \rho_0\rangle\nonumber\\
=&\langle f_2^{(k+1)}(u_0,u_{k+1})\sigma_{1,\beta},\rho_0\rangle+\langle f_4^{(k+1)}(u_0,\cdots,u_{k})\cdot\xi_\beta,\rho_0\rangle,\label{fkrho0}
\end{align}
where
$$f_2^{(k+1)}(u_0,u_{k+1})=2\rho_0\rho_{k+1}\big((2|\rho_0|^2-a^2-b^2)^2+2(|\rho
_0|^2-a^2)(|\rho_0|^2-b^2)\big),$$
which is independent of $\sigma_{k+1,\beta}.$ For $I_5,$  we yield from
\begin{align*}
  &g_{k+1}=\mathcal{D}_{*}u_{k+1}+L_{{k+1},0}\eta'(z)\bar\omega+\sum_\al L_{{k+1},\al}\rho_0^{-1}\eta'(z)\xi_{\alpha}+\sum_\al h_{{k+1},\alpha}\Big(2\rho_0'\eta'+\rho_0\eta'' \Big)\xi_{\alpha},\\
&\mathcal{D}_{*}u_i\cdot\xi_\beta|_\G=- \frac{\kappa_\beta}{d_0}\Big(2\eta'\rho_i'+\eta''\rho_i\Big),
\end{align*}
that
\begin{align*}
I_5=&\langle(d_1-z)g_{k+1}\cdot\xi_\beta, \rho_0\rangle\\
=&L_{{k+1}, \beta}\langle d_1-z,\eta'(z)\rangle+h_{{k+1},\beta}\langle (d_1-z)(2\rho_0'\eta'+\rho_0\eta''),\rho_0\rangle+\langle\mathcal{D}_*u_{k+1}\cdot\xi_\beta(d_1-z),\rho_0\rangle\\
=&L_{{k+1}, \beta}\langle d_1-z,\eta'(z)\rangle+h_{{k+1},\beta}\langle \rho_0^2,\eta' \rangle-\frac{\kappa_\beta}{d_0}\langle (d_1-z)(2\eta'\rho_{k+1}'+\eta''\rho_{k+1}),\rho_0\rangle.
\end{align*}
We have by \eqref{sigmak+1pm} that
\begin{align*}
  b^2h_{{k+1},\beta}+(b^2-a^2)\pa_\nu e_{{k+1}, \beta}+b\pa_\nu \sigma_{{k+1},\beta}^{*+}-a\pa_\nu \sigma _{{k+1},\beta}^{*-}-b^2\pa_\nu\Big(\frac{\kappa_\beta d_{k+1}}{d_0}\Big)=b\partial_\nu\sigma^+_{{k+1},\beta}-a\partial_\nu\sigma^-_{{k+1},\beta},
\end{align*}
and
\begin{align*}
    h_{{k+1},\beta}(x,t)=\frac{1}{ab}\Big(a\pa_\nu\sigma_{{k+1},\beta}^+-b\pa_\nu\sigma_{{k+1},\beta}^-+b\pa_\nu\sigma_{{k+1},\beta}^{*-}-a\pa_\nu\sigma_{{k+1},\beta}^{*+}\Big)+\pa_\nu\Big(\frac{\kappa_\beta d_{k+1}}{d_0}\Big).
\end{align*}
Therefore, we have
\begin{align*}
I_2+I_5=&-2b\pa_\nu\sigma_{{k+1}, \beta}^++2a\pa_\nu \sigma_{{k+1}, \beta}^-+\pa_\nu\Big( \frac{d_{k+1} \kappa_\beta }{d_0}\Big)\langle \eta',\rho_0^2\rangle-b^2\Big(\pa_\nu\Big(\frac{d_{k+1}\kappa_\beta}{d_0}\Big)-h_{k+1,\beta}\Big)\\
&+2\pa_\nu\langle\sigma_{{k+1}, \beta}^*, \rho_0'\rangle-h_{{k+1},\beta}\langle \eta',\rho_0^2\rangle+(b^2-a^2)\pa_\nu e_{{k+1}\beta}\\
&-2b\sum_{\al\neq \beta}\sigma_{{k+1}, \al}^+\pa_\nu\xi_\al^+\cdot\xi_\beta^++2a\sum_{\al\neq \beta}\sigma_{{k+1}, \al}^-\pa_\nu\xi_\al^-\cdot\xi_\beta^-+2\sum_{\al\neq \beta}\langle\sigma_{{k+1}, \al}^*\pa_\nu\xi_\al\cdot\xi_\beta,\rho_0'\rangle\notag\\
&-2d_{k+1}\sum_{\al\neq\beta}\frac{\kappa_\al}{d_0}\langle\eta\rho_0\pa_\nu\xi_\al\cdot \xi_\beta,\rho_0'\rangle+2\sum_{\al\neq \beta}{e_{{k+1}, \al}(x, t)}\langle\pa_\nu\xi_\al\cdot\xi_\beta\rho_0,\rho_0'\rangle\\
&-2\langle\pa_z\rho_{k+1}\pa_\nu\bar{\omega}\cdot\xi_\beta ,\rho_0\rangle-2\langle\rho_{k+1}\pa_{z\nu}\bar{\omega}\cdot \xi_\beta,\rho_0\rangle\\
&+L_{{k+1}, \beta}\langle d_1-z,\eta'(z)\rangle+h_{{k+1},\beta}\langle\rho_0^2,\eta' \rangle-\frac{\kappa_\beta}{d_0}\langle (d_1-z)(2\eta'\rho_{k+1}'+\eta''\rho_{k+1}),\rho_0\rangle\\
=&a\pa_\nu\sigma_{{k+1},\beta}^--b\pa_\nu \sigma_{{k+1},\beta}^++a\pa_\nu\sigma_{{k+1},\beta}^{*-}-b\pa_\nu\sigma_{{k+1}, \beta}^{*+}+\sum_{\al\neq \beta}\frac{\kappa_\al d_{k+1} b^2}{d_0}\pa_\nu\xi_\al^+\cdot\xi_\beta^+\\
&+\sum_{\al\neq \beta}-b\sigma_{{k+1},\al}^{*+}\pa_\nu\xi_\al^+\cdot\xi_\beta^++a\sigma_{{k+1},\al}^{*-}\pa_\nu\xi_\al^-\cdot\xi_\beta^-+2\pa_\nu\langle\sigma_{{k+1}, \beta}^* ,\rho_0'\rangle+\pa_\nu \Big(\frac{d_{k+1} \kappa_\beta}{d_0}\Big)\langle \eta',\rho_0^2\rangle\\
&-b\sum_{\al\neq \beta}\sigma_{{k+1}, \al}^+\pa_\nu\xi_\al^+\cdot\xi_\beta^++a\sum_{\al\neq \beta}\sigma_{{k+1}, \al}^-\pa_\nu\xi_\al^-\cdot\xi_\beta^-+2\sum_{\al\neq \beta}\langle\sigma_{{k+1}, \al}^*\pa_\nu\xi_\al\cdot\xi_\beta,\rho_0'\rangle\notag\\
&-2d_{k+1}\sum_{\al\neq\beta}\frac{\kappa_\al}{d_0}\langle\eta\rho_0\pa_\nu\xi_\al\cdot \xi_\beta, \rho_0'\rangle -2\langle \pa_z\rho_{k+1}\pa_\nu\bar{\omega}\cdot\xi_\beta, \rho_0\rangle-2\langle\rho_{k+1}\pa_{z\nu}\bar{\omega}\cdot \xi_\beta, \rho_0\rangle\\
&+L_{{k+1}, \beta}\langle d_1-z,\eta'(z)\rangle-\frac{\kappa_\beta}{d_0}\langle (d_1-z)(2\eta'\rho_{k+1}'+\eta''\rho_{k+1}),\rho_0\rangle.
\end{align*}
Indeed, by $\pa_z\xi_\al\cdot\xi_\beta=0$, we have used
\begin{align*}
    2\sum_{\al\neq\beta}&e_{{k+1},\al}(x,t)\langle \pa_\nu\xi_\al\cdot \xi_\beta, \rho_0\rho_0' \rangle\\=&\sum_{\al\neq\beta}e_{{k+1},\al}\Big(\pa_\nu\xi_\al\cdot\xi_\beta \rho_0^2|_{-\infty}^{+\infty}-\langle \pa_{z\nu}\xi_\al \cdot \xi_\beta, \rho_0^2\rangle\Big)\\
    =&\sum_{\al\neq \beta}e_{{k+1},\al}(b^2\pa_\nu\xi_\al^+\cdot\xi_\beta^+-a^2\pa_\nu\xi_\al^-\cdot \xi_\beta^-)-e_{{k+1},\al}\langle \pa_{z\nu}\xi_\al \cdot \xi_\beta, \rho_0^2\rangle\\
    =&\sum_{\al\neq \beta}\Big(b\sigma_{{k+1},\al}^++\frac{\kappa_\al d_{k+1} b^2}{d_0}-b\sigma_{{k+1},\al}^{*+}\Big)\pa_\nu\xi_\al^+\cdot\xi_\beta^+\\
    &-\Big(a\sigma_{{k+1},\al}^--a\sigma_{{k+1},\al}^{*-}\Big)\pa_\nu\xi_\al^-\cdot\xi_\beta^-.
\end{align*}
Combining $I_1-I_6$, we can get the mixed boundary condition \eqref{sigmakbetaboundary} for $\sigma_{{k+1},\beta}^\pm$ on $\G$.
\end{proof}
\section*{Acknowledgments}
This work is supported by NSF of China under Grant No. 12271476 and 11931010.

\end{document}